\documentclass{article}%
\usepackage{amsfonts}
\usepackage{amsmath}
\usepackage{color}
\usepackage{eurosym}
\usepackage{graphicx}
\usepackage[colorlinks=true,             linkcolor=blue,
citecolor=red,             urlcolor=blue]{hyperref}
\usepackage{tikz}
\usetikzlibrary{positioning,arrows.meta}
\usepackage{longtable}%
\usepackage{amssymb}
\providecommand{\U}[1]{\protect\rule{.1in}{.1in}}
\usetikzlibrary{positioning,arrows.meta}
\allowdisplaybreaks
\newtheorem{theorem}{Theorem}[section]

\newtheorem{corollary}[theorem]{Corollary}

\newtheorem{definition}[theorem]{Definition}

\newtheorem{lemma}[theorem]{Lemma}

\newtheorem{remark}[theorem]{Remark}

\newenvironment{proof}[1][Proof]{\noindent\textbf{#1.} }{\ \rule{0.5em}{0.5em}}
\begin{document}

\title{An $H^{2}$ extension for the flow of the cubic matrix Szeg\H{o} equation via unbounded
Hankel operators}
\author{%
Caixing Gu$^{1}$, Meng Li$^{2}$ and Pan Ma$^{2}$\\[0.6em]
{\small $^{1}$Department of Mathematics}\\
{\small California Polytechnic State University,
San Luis Obispo, USA}\\[0.4em]
{\small $^{2}$School of Mathematics and Statistics, HNP-LAMA}\\
{\small Central South University, Changsha, China}
}
\date{}
\maketitle
\begingroup
\makeatletter
\renewcommand{\@makefntext}[1]{%
  \setlength{\parindent}{0pt}%
  \noindent#1%
}
\makeatother
\renewcommand{\thefootnote}{}

\footnotetext{%
  \textit{E-mail addresses.}
  \href{mailto:cgu@calpoly.edu}{\texttt{cgu@calpoly.edu}}
  (C. Gu);
  \href{mailto:mengli97@126.com}{\texttt{mengli97@126.com}}
  (M. Li);
  \href{mailto:pan.ma@csu.edu.cn}{\texttt{pan.ma@csu.edu.cn}}
  (P. Ma).
}
\begin{abstract}

We construct a global flow for the cubic matrix Szeg\H{o}
equation on the whole Hardy space
$H^2(\mathbb T;\mathbb C^{m\times n})$,
at the regularity of the conserved mass.
The flow agrees with the regular evolution, depends
continuously on the initial datum, and has continuous
$H^2$ trajectories.
It preserves the $H^2$ norm and conserves the quartic
Hamiltonian whenever initially finite.
In the scalar case, it recovers the flow of
G\'{e}rard and Pushnitski.

\end{abstract}

\section{Introduction}
\label{introduction}

Let $n
\in\mathbb{N}_{+}$ be arbitrary fixed positive integers,
$\mathbb{C}^{n\times n}$ denotes the space of all $n\times n$ complex
matrices. We write
\[
L^2_{n\times n}:=L^2(\mathbb T;\mathbb C^{n\times n}),
\qquad
H^2_{n\times n}:=H^2(\mathbb T;\mathbb C^{n\times n}),
\]
where $H^2_{n\times n}$ consists of the functions in
$L^2_{n\times n}$ with vanishing negative Fourier coefficients.
We study the cubic matrix Szeg\H{o} equation on the unit circle $\mathbb T,$
\begin{align}
	\label{cubic matrix}i\partial_{t} U=\mathbb{P}(UU^{\ast} U),\quad U
	\in\mathbb{C}^{n \times n},\quad(t,z)\in\mathbb{R}\times\mathbb{T},
\end{align}
where $\mathbb{P}$ denotes the matrix
Szeg\H{o} projection from $L^2_{n\times n}$ onto
$H^2_{n\times n}. $  We present the operator-theoretic construction for square matrices; the rectangular
case follows by zero padding, as explained in Remark \ref{zero}.

Our aim is to extend the regular evolution to the whole
matrix-valued Hardy space $H^2_{n\times n}$, with
continuous dependence on the initial datum, preservation of the
$H^2_{n\times n}$ norm, and conservation of the Hamiltonian whenever the latter is
initially finite.

We use the
Hilbert--Schmidt matrix norm in $L^2$ and write
$d\lambda(z)=|dz|/(2\pi)$ for the normalized Lebesgue measure on $\mathbb{T}$.
For sufficiently regular data, \eqref{cubic matrix} is Hamiltonian with
respect to
\[
\omega(U,V)=\operatorname{Im}\int_{\mathbb{T}}
\operatorname{tr}(UV^{\ast})\,d\lambda,
\]
with quartic Hamiltonian
\begin{equation*}
	E(U)=\frac14\int_{\mathbb{T}}
	\operatorname{tr}\bigl((UU^{\ast})^2\bigr)\,d\lambda
	=\frac14\|UU^{\ast}\|_{L^2_{n\times n}}^2,
\end{equation*}
defined whenever $UU^{\ast}\in L^2_{n\times n}$.

For $m=n=1$, equation \eqref{cubic matrix} reduces to
\begin{equation}\label{cubic scalar}
	i\partial_t u=\mathbb{P}(|u|^2u).
\end{equation}

G\'{e}rard and Grellier introduced \eqref{cubic scalar} as a model  for completely non-dispersive Hamiltonian dynamics on Hardy space \cite{Gerard06,Gerard08}.  They discovered its Hankel Lax pair structure and
subsequently derived an explicit formula expressing the solution in terms
of the initial datum 
\cite{Gerard10,Gerard12, Gerard12-Anal,Gerard15,Gerard17} . These results
demonstrate that Hankel operators are not merely a convenient analytic tool,
but are intrinsic to the integrable structure of the equation, playing a
central role in its conservation laws and explicit evolution formulas. This operator representation
provides a way to study the evolution beyond the regularity at which the
equation can be treated directly.  

G\'{e}rard and Pushnitski \cite{Gerard23} established a continuous extension of the scalar Szeg\H{o} flow to the whole Hardy space \(H^2(\mathbb T)\), preserving the \(H^2\) norm and conserving the Hamiltonian on the finite-energy class. Their proof combines the domain theory of unbounded Hankel operators with the operator-theoretic structure of the explicit flow formula.

For $s\geq 0$, let
\[
\mathcal{H}_{+}^{s}(\mathbb{T},\mathbb{C}^{n\times n})
:=\mathcal{H}^{s}(\mathbb{T},\mathbb{C}^{n\times n})
\cap H^{2}(\mathbb{T},\mathbb{C}^{n\times n})
\]
denote the matrix-valued analytic Sobolev space. In particular,
$H^2(\mathbb T;\mathbb C^{n\times n})
=\mathcal H_+^0(\mathbb T;\mathbb C^{n\times n})$.
For the matrix equation, Sun \cite{Sun} established global well-posedness
and two equivalent explicit flow formulas on
$\mathcal{H}_+^{1/2}(\mathbb{T};\mathbb{C}^{n\times n})$.
The formulas involve the generally distinct left and right matrix
Hankel operators and their shifted counterparts. Their products satisfy
Lax equations for solutions in $\mathcal{H}_+^s$, $s>1/2$.
Related developments on the real line include traveling waves, explicit
solutions, and low-regularity flows
\cite{Gerard-Pu24,GerardJEMS,Pocovnicu11,Pocovnicu11-2,Sun2};
multiphase and multisoliton solutions were constructed in \cite{mats}.

The matrix product in \eqref{cubic matrix} couples the entries of $U$,
so the evolution cannot in general be obtained by applying the scalar
flow entrywise. For general $U\in H^2_{n\times n}$,
the associated Hankel operators
may also be unbounded.
This leads to the following question :

\medskip\noindent\textbf{Main question.} \emph{Does the cubic matrix Szeg\H{o}
	flow admit a continuous extension to $H^{2}_{n\times n}$, with conservation
	laws preserved?} \medskip

From an operator-theoretic viewpoint, this raises the question of whether
the approach of G\'{e}rard and Pushnitski can be adapted to the noncommutative
matrix setting. From a broader viewpoint, this question fits naturally into
the program emphasized in G\'{e}rard's recent ICM lecture \cite{Gerard26},
where Hardy-space operator models, Lax pair structures, and explicit
formulas are presented as a general framework for the analysis of integrable
nonlinear equations. Our main theorem gives an affirmative answer to the question above.

Let $\mathcal{W}(t)$, $t\in\mathbb{R}$, denote the (nonlinear) flow map
associated with (\ref{cubic matrix}), initially defined for $U_{0}%
\in\mathcal{H}^{s}_{+}( \mathbb{T}, \mathbb{C}^{n\times n})$, $s>1/2$.
Thus, if $U(t,\cdot)$ is the corresponding solution with initial datum
$U(0,z)=U_{0}(z)$, one has
\[
\mathcal{W}(t)U_{0}=U(t,\cdot).
\]

For $U\in H^{2}_{n\times n}$, we use the same notation $\mathcal{W}(t)U$ for
the value at time $t$ of the extended flow with initial datum $U$. 

\begin{theorem}
	[Main Result]\label{main result} The flow map $\mathcal{W}(t)$ admits an
	extension to $H^{2}_{n\times n}$ with the following properties: \newline(i)
	$\Vert\mathcal{W}(t)U\Vert_{H^{2}_{n \times n}}=\Vert U\Vert_{H^{2}_{n \times
			n}} $ for all $U\in H^{2}_{n \times n}$; \newline(ii) if $\Vert U_{k}%
	-U\Vert_{H^{2}_{n \times n}}\to0$, then $\Vert\mathcal{W}(t)U_{k}%
	-\mathcal{W}(t)U\Vert_{H^{2}_{n \times n}}\to0$ as $k \to\infty$;
	\newline(iii) the map $t\mapsto\mathcal{W}(t)U$ is continuous in $H^{2}_{n
		\times n}$ for all $U \in H^{2}_{n \times n}$; \newline(iv) if $UU^{\ast}\in
	L^{2}_{n\times n}$, then
	\begin{align*}
		(\mathcal{W}(t)U)(\mathcal{W}(t)U)^{\ast}\in L^{2}_{n\times n}%
	\end{align*}
	and
	\begin{align*}
		& \| (\mathcal{W}(t)U)(\mathcal{W}(t)U)^{\ast} \|_{L^{2}_{n \times n}} =\|
		UU^{\ast} \|_{L^{2}_{n\times n}},
	\end{align*}
	that is
	\begin{align*}
		E(\mathcal{W}(t)U) =E(U).
	\end{align*}
	
\end{theorem}

\begin{remark}
	For $n=1$, the condition $UU^{\ast}\in L^{2}$ is equivalent to $u\in H^{4}(
	\mathbb{T})$, and part~(iv) reduces to the scalar Hamiltonian conservation law
	of G\'{e}rard and Pushnitski \cite{Gerard23}. Thus, Theorem
	\ref{main result} provides a natural matrix-valued extension of their result.
\end{remark}
\begin{remark}\label{zero}
	For any positive integers $m,n$,
embed $\mathbb{C}^{m\times n}$ into
	$\mathbb{C}^{N\times N}$, $N=\max\{m,n\}$, by appending zero rows or
	columns as needed. The nonlinear term $UU^{\ast}U$ preserves the resulting zero pattern, so the
	embedded subspace is invariant under the regular flow by uniqueness. Since this subspace is closed in
	$H^{2}(\mathbb{T},\mathbb{C}^{N\times N})$,
	regular approximation and
	Theorem \ref{main result}(ii) imply
	that it is also invariant under
	the extended flow.
	Consequently, the conclusions of
	Theorem \ref{main result} extend to
	arbitrary rectangular matrix-valued
	initial data.
\end{remark}

We next describe the main ideas of the proof. The first step is to
establish the basic operator-theoretic properties of unbounded matrix
Hankel operators. For $U\in H^{2}_{n\times n}$, we consider the maximal
realizations of the right and left Hankel operators and prove that they are
closed. %
The commutation relation for the right Hankel operator
implies that its graph is invariant under the
direct sum of the forward and backward shift operators. %
After vectorization,
this leads to a vector-valued invariant
subspace problem, to which the invariant subspace theorem of Gu and
Luo \cite{GuLuo} applies. %
Their result plays a key role in our proof,
particularly in establishing the graph-norm
density of polynomials. %
The adjoint relations for the Hankel operators then allow
us to establish  the positivity and self-adjointness of the products
\[
H_{U}^{(r)}H_{U}^{(l)}  \quad\text{and}\quad H_{U}^{(l)}H_{U}^{(r)},
\]
and to prove strong resolvent convergence of these products when the
symbols converge in $H^{2}$. This strong resolvent convergence is a key
ingredient in proving the continuity of the extended flow.

The $H^{2}$-norm conservation requires a different operator-theoretic
argument. From the unitary groups generated by the Hankel products, we
construct isometries $\Omega^{(r)}$ and $\Omega^{(l)}$ and analyze their Wold
decompositions. The key step is to show that $U$ belongs to the pure
part of the isometry $\Omega^{(l)}$. Once this is established, the Wold
decomposition provides an orthonormal basis for the pure subspace, and
Parseval's identity applied to the explicit formula gives %
\[
\|\mathcal{W}(t)U\|_{H^{2}_{n\times n}}  =  \|U\|_{H^{2}_{n\times n}}.
\]

For the Hamiltonian conservation law, we introduce a resolvent functional
$J(x,U)$ associated with $H_{U}^{(r)}H_{U}^{(l)}$. Its second
derivative at $x=0_{+}$ determines $\|UU^{\ast}\|_{L^{2}_{n\times n}}^{2}$
through an identity involving the zeroth Fourier coefficient $\int_{
\mathbb{T}}UU^{\ast}\,d\lambda$. The Lax pair yields conservation of $J(x,U)$
and of this zeroth Fourier coefficient for $U\in\mathcal{H}_{+}^{s}(
\mathbb{T}, \mathbb{C}^{n\times n})$ with $s>1/2$. Strong resolvent
convergence then extends the conservation of $J(x,U)$ to general $H^{2}$ data,
which yields the conservation of $\|UU^{\ast}\|_{L^{2}_{n\times n}}^{2}$ and
proves part (iv) of the main theorem.

The paper is organized as follows.

In Section \ref{Auxiliary results}, we present several auxiliary theorems that
provide the operator-theoretic foundation for the proof of Theorem
\ref{main result}.

In Section \ref{proof of main result}, we prove Theorem \ref{main result}%
(i)-(iv). The proof of part (i) forms the core of the paper. Its main
ingredients are the domain theory of unbounded anti-linear Hankel operators
(Theorem \ref{Domain dense left and right Hankel}), the explicit formula for
the flow map (Theorem \ref{Sun formula}), %
and the fact that $U$ belongs to the pure part
of the isometry $\Omega^{(l)}$ in its Wold
decomposition (Theorem \ref{complicated}). 
Parts (ii) and (iii) rely on the strong resolvent
convergence of the products $H_{U}^{(r)}H_{U}^{(l)}$ (Theorem
\ref{strong con-H}). Part (iv) is established by combining the conservation laws
in Lemmas \ref{conservation law} and \ref{conservation-Q}, %
yielding the conservation of the Hamiltonian energy.

In Section \ref{Polynomials are dense}, we establish Theorem
\ref{Domain dense}. The key ingredient is the characterization of invariant
subspaces of $S_{E}\oplus S_{F}^{\ast}$ on the vector-valued Hardy space
$H_{E}^{2}\oplus H_{F}^{2}$ due to Gu and Luo \cite{GuLuo}, which replaces
Timotin's scalar result and constitutes one of the principal new ingredients
of the matrix-valued theory.

Section \ref{proof of complicated} is devoted to the proof of Theorem
\ref{complicated}. The central step is to show that $U$ belongs to the pure
part of the isometry $\Omega^{(l)}$. The argument is based on two structural
results: Theorem \ref{Wold reduce} shows that $\mathcal{H}^{(l)}$  maps the
pure part of the Wold decomposition of $\Omega^{(r)}$  into that of
$\Omega^{(l)}$, and Lemma \ref{no a.c.part} which establishes that the unitary
parts of the particular isometries $\Omega^{(r)}$ and $\Omega^{(l)}$ have no
absolutely continuous spectral components.

Finally, in Section  \ref{proof of strong convergence}, we establish the
strong resolvent convergence for the operators $H^{(r)}_{U}H^{(l)}_{U}$
(Theorem  \ref{strong con-H}) which is used in Section
\ref{proof of main result}. The appendix records the operator conventions
needed for \cite{GuLuo} and provides, in the present notation, the
details of the proofs of the Lax pair identities and Sun's explicit
formula.

Several operator-theoretic results used in Section~3 are stated
in Section~2 and proved in Sections~4--6.
Figure~\ref{fig:proof-roadmap} groups the main ingredients
according to the conclusions of Theorem~1.1 that they establish.
The norm identity proved in part~(i) is used again in
parts~(ii) and~(iii), where it upgrades weak convergence
to strong convergence in $H^2$.

\begin{figure}[!htbp]
  \centering
  \makebox[\linewidth][c]{%
    \resizebox{1.10\linewidth}{!}{%
      \begin{tikzpicture}[
        every node/.style={
          align=center,
          minimum width=8mm,
          minimum height=3mm,
          inner sep=1pt,
          font=\footnotesize
        }
      ]

      \node (main) {\textbf{Theorem~1.1}\\
        (Sec.~\ref{introduction})};
      \node[below left=15mm and 28mm of main] (p1)
        {\textbf{Part (i)}\\(Sec.~\ref{1sec})};
      \node[below=15mm of main] (p2)
        {\textbf{Parts (ii)--(iii)}\\(Sec.~\ref{23sec})};
      \node[below right=15mm and 35mm of main] (p3)
        {\textbf{Part (iv)}\\(Sec.~\ref{4sec})};

      \draw (p1) -- (main);
      \draw (p2) -- (main);
      \draw (p3) -- (main);

      \node[below left=22mm and 6mm of p1] (t29)
        {\textbf{Extension of the flow map}\\
         Thm.~\ref{Sun formula} (Sec.~\ref{Extension flow})};
      \node[below right=22mm and 6mm of p1] (t210)
        {\textbf{Pure part of the isometry}\\
         Thm.~\ref{complicated} (Sec.~\ref{pure part sec})};

      \draw (t29) -- (p1);
      \draw (t210) -- (p1);

      \node[below=9mm of p2] (t23a)
        {\textbf{Strong resolvent convergence}\\
         Thm.~\ref{strong con-H} (Sec.~\ref{auxiliary})};

      \draw (t23a) -- (p2);

      \node[below left=22mm and 6mm of p3] (t31a)
        {\textbf{Lax pair}\\
         Thm.~\ref{Sun lax pair} (Sec.~\ref{sun result})};
      \node[below right=22mm and 6mm of p3] (t32a)
        {\textbf{Conservation laws}\\
         Lemmas~\ref{conservation law} and~\ref{conservation-Q}
         (Sec.~\ref{4sec})};

      \draw (t31a) -- (p3);
      \draw (t32a) -- (p3);

      \node[below left=15mm and 10mm of t29] (t14)
        {\textbf{Sun's explicit formula}\\
         Thm.~\ref{Sun explicit formula} (Sec.~\ref{sun result})};
      \node[below right=15mm and 10mm of t29] (t23b)
        {\textbf{Domain properties of}\\
         \textbf{unbounded Hankel operators}\\
         Thm.~\ref{Domain dense left and right Hankel}
         (Sec.~\ref{auxiliary})};

      \draw (t14) -- (t29);
      \draw (t23b) -- (t29);

      \node[below left=15mm and 6mm of t23b] (t21)
        {Thm.~\ref{Domain dense} (Sec.~\ref{auxiliary})};
      \node[below right=15mm and 6mm of t23b] (t22)
        {Thm.~\ref{Domain dense left} (Sec.~\ref{auxiliary})};

      \draw (t21) -- (t23b);
      \draw (t22) -- (t23b);

      \end{tikzpicture}%
    }%
  }\par
  \caption{Main ingredients in the proof of Theorem~1.1.}
  \label{fig:proof-roadmap}
\end{figure}
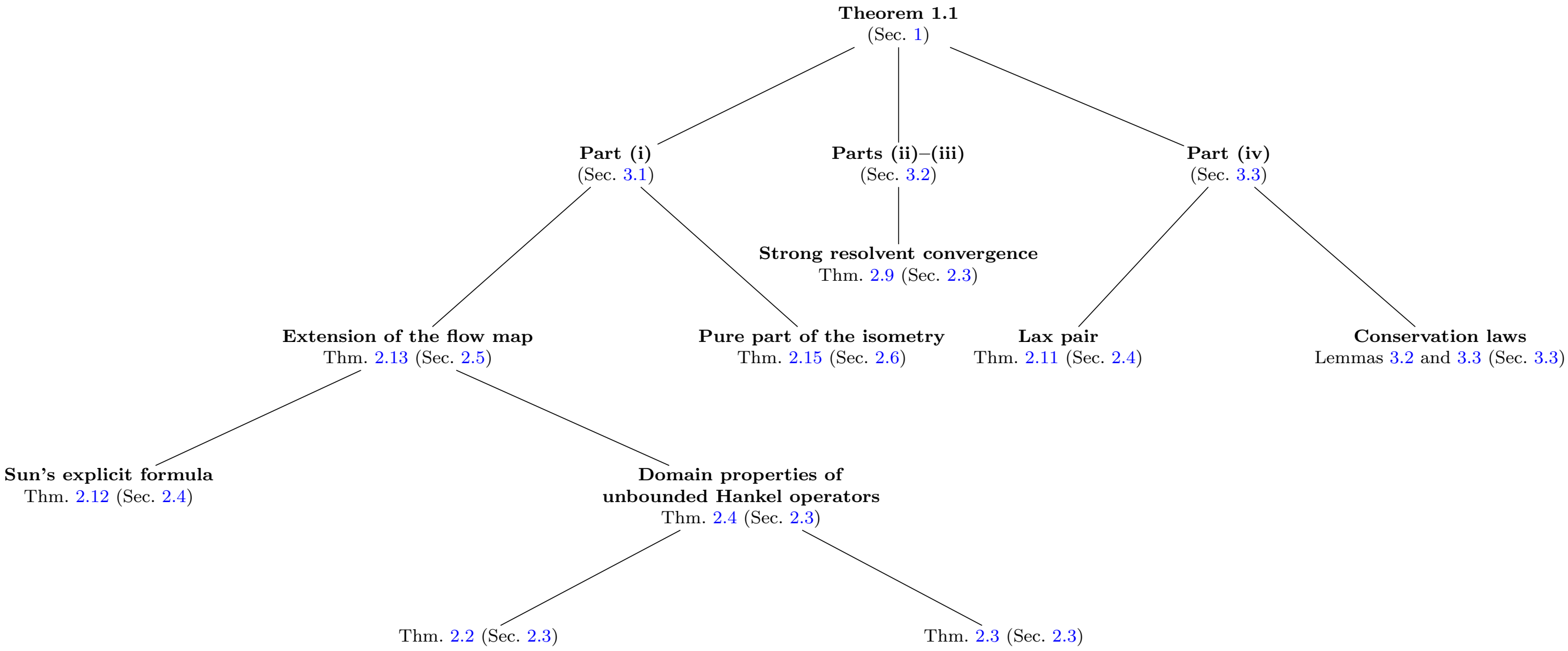

\clearpage
\section{Operator-theoretic framework for the proof of main result}

\label{Auxiliary results}

Since the proof of Theorem \ref{main result}(i) is rather involved, we first
present four auxiliary theorems. The first theorem (Theorem
\ref{Domain dense left and right Hankel}) establishes the basic properties of
the left and right anti-linear Hankel operators, including closedness, density
of polynomials in their graph norms, and their adjoint relations. The second
theorem (Theorem \ref{strong con-H}) proves the strong resolvent convergence
of the operators $H_{U}^{(r)}H_{U}^{(l)}$. The third theorem (Theorem
\ref{Sun formula}) extends the explicit formula of Theorem
\ref{Sun explicit formula} to define the flow map on the whole space
$H^{2}_{n\times n}$. The fourth theorem (Theorem \ref{complicated}) shows that
$U$ belongs to the pure part of the isometry $\Omega^{(l)}$.

\subsection{Matrix-valued Hardy spaces}

\label{Notions and preliminaries}

Let $L^{2}=L^{2}(\mathbb{T})$ be the space of square-integrable functions on
the unit circle $\mathbb{T}$ with respect to the normalized Lebesgue measure.
Let $H^{2}$ be the Hardy space on the open unit disk $\mathbb{D}$. $L^{\infty
}$ and $H^{\infty}$ are the algebras of bounded functions in $L^{2}$ and
$H^{2}$ respectively. Suppose that $E$ and $F$ are separable complex
Hilbert spaces. Let $B(E,F)$ denote the bounded linear operators from $E$ to
$F$, and write $B(E)=B(E,E)$.
$L_{E}^{2}$ and $H_{E}^{2}$ denote $E$-valued $L^{2}$ and $H^{2}$ spaces.
$L_{B(E)}^{\infty}$ and $H_{B(E)}^{\infty}$ are operator-valued $L^{\infty}$
and $H^{\infty}$ algebras.

For $1\le p\le\infty$, we write
\[
L^{p}_{n\times n}:=L^{p}( \mathbb{T}, \mathbb{C}^{n\times n}), \qquad
H^{p}_{n\times n}:=H^{p}( \mathbb{T}, \mathbb{C}^{n\times n}).
\]
The matrix-valued $L^{p}$-space $L^{p}_{n\times n}$ consists of all
(equivalence classes of) measurable functions $U:\mathbb{T}\rightarrow
\mathbb{C}^{n\times n}$. It is well known that $L^{p}_{n\times n}$ is a Banach
space with respect to the $L^{p}$-norm:
\[
\Vert U\Vert_{L^{p}_{n\times n}}:=%
\begin{cases}
\left(  \int_{\mathbb{T}}\Vert U(z)\Vert_{\mathbb{C}^{n\times n}}^{p}%
d\lambda(z)\right)  ^{1/p}, & 1\leq p<\infty,\\
\text{ess sup}_{z\in\mathbb{T}}\Vert U(z)\Vert_{\mathbb{C}^{n\times n}}, &
p=\infty,
\end{cases}
\]
where $d\lambda(z)$ is the normalized Lebesgue measure on $\mathbb{T}$ and the
norm $\Vert\cdot\Vert_{\mathbb{C}^{n\times n}}$ denotes the Hilbert--Schmidt
norm on $\mathbb{C}^{n\times n}$, defined by $\Vert A\Vert_{\mathbb{C}%
^{n\times n}}^{2}=\operatorname{tr}(AA^{\ast})$.

The Hilbert--Schmidt inner product on $\mathbb{C}^{n \times n}$ is defined by
\[
\langle A, B \rangle_{\mathbb{C}^{n \times n}} := \operatorname{tr}(A B^{\ast
}).
\]
The corresponding inner product on $L^{2}_{n\times n}$ is
\[
\langle U, V \rangle_{L^{2}_{n \times n}} := \int_{\mathbb{T}}
\operatorname{tr}(U(z) V(z)^{\ast}) d\lambda(z).
\]

For any $U \in L^{p}_{n \times n} $, its Fourier coefficients are $n \times n
$ complex matrices defined as:
\[
\hat{U}(k) := \int_{\mathbb{T}} U(z) z^{-k} d\lambda(z) \in\mathbb{C}^{n
\times n} \quad\text{for each } k \in\mathbb{Z}.
\]
The matrix-valued Hardy space $H^{p}(\mathbb{T}, \mathbb{C}^{n \times n})$ ($1
\leq p \leq\infty$) is then defined as:
\[
H^{p}_{n \times n} := \left\{  U \in L^{p}_{n \times n} \mid\hat{U}(k) = 0
\text{ for all } k < 0 \right\} .
\]
As a closed subspace of $L^{p}_{n\times n}$, $H^{p}_{n\times n}$ inherits the
Banach space structure from $L^{p}_{n\times n}$. Its norm is the restriction
of the $L^{p} $-norm:
\[
\| U \|_{H^{p}_{n \times n}} := \| U \|_{L^{p}_{n \times n}}.
\]

We next recall the realization of Hardy spaces in terms of holomorphic
functions on the unit disk $\mathbb{D}$.

Let $Hol(\mathbb{D},\mathbb{C}^{n\times n})$ denote the space of all
$\mathbb{C}^{n\times n}$-valued analytic functions, that is, $f\in
Hol(\mathbb{D},\mathbb{C}^{n\times n})$ can be written as%
\[
f(z)=\sum_{l\geq0}f_{l}z^{l},z\in\mathbb{D}\text{ and }f_{l}\in\mathbb{C}%
^{n\times n},
\]
and the above infinite series is convergent in the Hilbert--Schmidt norm on
$\mathbb{C}^{n\times n}$.

The Hardy space on the disk $\mathbb{D}$ is defined by%
\begin{align*}
H^{2}(\mathbb{D},\mathbb{C}^{n\times n}) =\left\{  f\in Hol(\mathbb{D}%
,\mathbb{C}^{n\times n}):\left\Vert f\right\Vert _{H^{2}(\mathbb{D}%
,\mathbb{C}^{n\times n})}<\infty\right\}  ,
\end{align*}
where%
\begin{align*}
\|f\|_{H^{2}(\mathbb{D}, \mathbb{C}^{n\times n})}  & = \sup_{0<r<1} \left(
\int_{ \mathbb{T}} \|f(r\zeta)\|_{ \mathbb{C}^{n\times n}}^{2} \,d\lambda
(\zeta) \right) ^{1/2}\\
& = \sup_{0<r<1} \left(  \int_{ \mathbb{T}} \operatorname{tr} \left(
f(r\zeta)f(r\zeta)^{\ast} \right)  \,d\lambda(\zeta) \right) ^{1/2}.
\end{align*}

By the radial limit theorem, the Hardy spaces $H^{2}(\mathbb{D},
\mathbb{C}^{n\times n})$ and $H^{2}(\mathbb{T}, \mathbb{C}^{n\times n})$ are
isometrically isomorphic. More precisely, every $f \in H^{2}(\mathbb{D},
\mathbb{C}^{n\times n})$ admits almost everywhere non-tangential (in
particular radial) boundary limits on $\mathbb{T}$, which define an element of
$H^{2}(\mathbb{T}, \mathbb{C}^{n\times n})$.

Conversely, every $g \in H^{2}(\mathbb{T}, \mathbb{C}^{n\times n})$ admits a
holomorphic extension to $\mathbb{D}$ given by the Poisson integral, which
belongs to $H^{2}(\mathbb{D}, \mathbb{C}^{n\times n})$.

Under this isometric identification, we do not distinguish between
$H^{2}(\mathbb{D}, \mathbb{C}^{n\times n})$ and $H^{2}(\mathbb{T},
\mathbb{C}^{n\times n})$, and use the notation $H^{2}_{n\times n}$ for either
realization.

The case $p=2$ is of particular importance in this paper. Indeed, the
Hilbert--Schmidt norm equips $L^{2}_{n\times n}$ with a natural Hilbert space
structure. In particular, the $L^{2}$-norm admits an entrywise decomposition.
Such a Hilbert space structure is specific to the case $p=2$. It is easy to
see that a matrix function
\begin{align*}
U=(U_{kj})_{1 \leq k, j \leq n}%
\end{align*}
belongs to $L^{2}_{n \times n}$ if and only if each entry $U_{kj}$ belongs to
$L^{2}(\mathbb{T})$. We set
\begin{align*}
\Vert U\Vert^{2}_{L^{2}_{n \times n}}:=\sum_{k,j=1}^{n}\Vert U_{kj}\Vert
^{2}_{L^{2}(\mathbb{T})}=\int_{\mathbb{T}} \operatorname{tr}(U(z)U(z)^{\ast})
d\lambda(z).
\end{align*}
Similarly,
\[
\Vert U\Vert^{2}_{H^{2}_{n \times n}}:=\sum_{k,j=1}^{n}\Vert U_{kj}\Vert
^{2}_{H^{2}(\mathbb{T})}.
\]

More explicitly, the intersection of the Sobolev space with
$H^{2}(\mathbb{T},\mathbb{C}^{n \times n})$ 
can be characterized as
\begin{align*}
\mathcal{H}^{s}_{+}(\mathbb{T},\mathbb{C}^{n \times n}) 
&  = \{U\in H^{2}(\mathbb{T},\mathbb{C}^{n \times n}): \sum_{l=0}^{\infty}
\operatorname{tr}\left(  \hat{U}(l) \left(  \hat{U}(l)\right)  ^{\ast}
\right)  (1+l^{2})^{s}<\infty\} .
\end{align*}
Then equation (\ref{cubic matrix}) defines a smooth flow on $\mathcal{H}%
^{s}_{+}(\mathbb{T},\mathbb{C}^{n \times n})$ for $s>1/2$, and a continuous
flow on $\mathcal{H}^{1/2}_{+}(\mathbb{T},\mathbb{C}^{n \times n})$.

\subsection{Matrix-valued Hankel and Toeplitz operators}

\label{Matrix-valued intro}

Before discussing the flow on $H^{2}_{n\times n}$, we briefly recall the
motivation for introducing the $\mathrm{BMOA} (\mathbb{T},
\allowbreak\mathbb{C}^{n\times n})$ space. Recall that $U\in\mathrm{BMOA}%
(\mathbb{T}, \mathbb{C}^{n\times n})$ if and only if each entry $U_{kj}$ lies
in the scalar space $\mathrm{BMOA}(\mathbb{T})$.

Let $\mathcal{P}_{n \times n}$ denote the subspace of $H^{2}_{n \times n}$
consisting of matrix-valued analytic polynomials, i.e.
\[
\mathcal{P}_{n \times n}=\{p\in H^{2}_{n \times n}:p(z)=\sum_{i=0}^{m}%
C_{i}z^{i},C_{i}\in\mathbb{C}^{n \times n},m\geq0\}.
\]

Let $U\in H^{2}_{n\times n}$. The (possibly unbounded) right Hankel operator
$\Gamma_{U}^{(r)}$ with the symbol $U$ is initially defined on $\mathcal{P}_{n
\times n}$ by
\begin{equation}
\label{def of right Hankel}\Gamma_{U}^{(r)}(f)=\mathbb{P}(UJf),\quad f
\in\mathcal{P}_{n \times n},
\end{equation}
where $Jf(z)=f(\overline{z})$.

Similarly, the (possibly unbounded) left Hankel operator $\Gamma_{U}^{(l)}$
with the symbol $U$ is initially defined on $\mathcal{P}_{n \times n}$ by
\begin{equation}
\label{def of left Hankel}\Gamma_{U}^{(l)}(g)=\mathbb{P}((Jg) U),g
\in\mathcal{P}_{n \times n}.
\end{equation}

In the scalar case where $u\in H^{2}(\mathbb{T}).$ Then $\Gamma_{u}%
:=\Gamma_{u}^{(r)}=\Gamma_{u}^{(l)}.$ By the Nehari-Fefferman Theorem
\cite[Theorem 1.1.2]{Peller}, the operator $\Gamma_{u}$, initially defined on
$\mathcal{P}$, extends to a bounded operator on $H^{2}$ if and only
if $u\in\mathrm{BMOA}(\mathbb{T}).$ Thus we immediately have the following
result in the matrix-valued case.

\begin{theorem}
\label{Nehari} Let $U\in H^{2}_{n\times n}$. Then $\Gamma_{U}^{(r)},$
initially defined on $\mathcal{P}_{n \times n},$ extends to be a bounded
linear operator from $H^{2}_{n \times n}$ to $H^{2}_{n \times n}$ if and only
if $U\in\mathrm{BMOA}(\mathbb{T},\mathbb{C}^{n\times n})$.
\end{theorem}

This connection between $\mathrm{BMOA}(\mathbb{T},\mathbb{C}^{n\times n})$ and
(possibly unbounded) Hankel operators motivates the study of flows on
$\mathrm{BMOA}(\mathbb{T},\mathbb{C}^{n\times n})$, which will later be
extended to the larger space $H^{2}_{n\times n}$.

We now focus on the case where $U$ does not belong to $\mathrm{BMOA}%
(\mathbb{T},\mathbb{C}^{n\times n})$, i.e. at least one entry $U_{kj}$ is not
in the scalar $\mathrm{BMOA}(\mathbb{T})$ space. Since the right Hankel
operator may be unbounded in this case, it is necessary to specify carefully
how $\Gamma^{(r)}_{U}$ is defined. We define it by its maximal domain:
\begin{align}
\label{Dom of right Hankel}\operatorname{Dom} \Gamma^{(r)}_{U}=\{f \in
H^{2}_{n \times n} : \mathbb{P}(UJf)\in H^{2}_{n \times n}\}.
\end{align}
Here, the product $UJf$ is understood as the pointwise matrix multiplication
on $\mathbb{T}$. The condition $\mathbb{P}(UJf)\in H^{2}_{n \times n}$ is
understood in the distributional sense, namely, that $UJf\in L^{1}_{n \times
n}$ and the nonnegative Fourier coefficients of each matrix entry are square
summable. In this case, we define $\Gamma^{(r)}_{U}(f)=\mathbb{P}(UJf)$.

Throughout this paper, $\Gamma_{U}^{(r)}$ always refers to the maximal
realization. We will show that the polynomials $\mathcal{P}_{n \times n}$ are
dense in $\operatorname{Dom}\Gamma_{U}^{(r)}$ with respect to the graph norm
in Theorem \ref{Domain dense}.

Similarly, we can define $\operatorname{Dom} \Gamma^{(l)}_{U}$ by
(\ref{def of left Hankel}),
\begin{align}
\label{Dom of left Hankel}\operatorname{Dom} \Gamma^{(l)}_{U}=\{g \in H^{2}_{n
\times n}: \mathbb{P}((Jg) U) \in H^{2}_{n \times n}\}.
\end{align}

Next, we will introduce the anti-linear version of the Hankel operator.

Now let $U\in H^{2}_{n\times n}$, so that $H_{U}^{(r)}$ and $H_{U}^{(l)}$ may
be unbounded. Define the conjugation operator $\mathcal{C}$ from $H^{2}_{n
\times n}$ to $H^{2}_{n \times n}$ by
\[
(\mathcal{C}f)(z):=f(\overline{z})^{\ast},f\in H^{2}_{n \times n}.
\]
For convenience, write $f^{\#}:=\mathcal{C}f$.

Set
\begin{equation}
\label{eq:dom-left-anti-hankel}\text{Dom } H^{(r)}_{U}=\{f \in H^{2}_{n \times
n} :\mathbb{P}(Uf^{\ast})\in H^{2}_{n \times n}\}
\end{equation}
and define
\begin{equation}
\label{def-right-H}H^{(r)}_{U} (f):=\Gamma^{(r)}_{U}(\mathcal{C}%
f)=\mathbb{P}(Uf^{\ast}), \quad f \in\text{Dom } H^{(r)}_{U}.
\end{equation}
The graph norm of $H^{(r)}_{U}$ is defined as
\begin{equation}
\label{graph norm-H}\Vert f\Vert_{H^{(r)}_{U}}:= \left(  \Vert f \Vert
^{2}_{H^{2}_{n \times n}}+\Vert H^{(r)}_{U} (f)\Vert^{2}_{H^{2}_{n \times n}%
}\right)  ^{1/2}, \quad f \in\text{Dom } H^{(r)}_{U}.
\end{equation}

In addition, set
\begin{equation}
\label{eq:dom-right-anti-hankel}\text{Dom } H^{(l)}_{U}=\{ g \in H^{2}_{n
\times n}:\mathbb{P}(g^{\ast}U) \in H^{2}_{n \times n}\}
\end{equation}
and define
\begin{equation}
\label{def-left-H}H^{(l)}_{U}(g):=\Gamma^{(l)}_{U}(\mathcal{C}g)=\mathbb{P}%
(g^{\ast}U), \quad g \in\text{Dom } H^{(l)}_{U}.
\end{equation}

Unlike the Hankel operators considered above, all Toeplitz operators
appearing in this paper are bounded. 

Let $V\in L^{\infty}(\mathbb{T},\mathbb{C}^{n\times n})$. The right
Toeplitz operator $T_{V}^{(r)}$ with symbol $V$ is defined by
\[
T_{V}^{(r)}f=\mathbb{P}(Vf),\quad f\in H^{2}_{n\times n}.
\]
The left Toeplitz operator $T_{V}^{(l)}$ with symbol $V$ is defined by
\[
T_{V}^{(l)}g=\mathbb{P}(gV),\quad g\in H^{2}_{n\times n}.
\]


Let $\mathbb{E}_{kj} \in\mathbb{C}^{n\times n}$ denote the $n\times n$ matrix
whose $kj$-entry is $1$ and the other entries are all $0$, then $\{\mathbb{E}%
_{kj}z^{m}:1\leq k,j\leq n,m\geq0\}$ forms an orthonormal basis of the Hilbert
space $H^{2}_{n\times n}$. Set
\[
\Delta_{n\times n}:=\sum_{k,j=1}^{n} (\mathbb{E}_{kj}\otimes\mathbb{E}_{kj}),
\]
where $(\mathbb{E}_{kj}\otimes\mathbb{E}_{kj})$ denotes the rank-one operator.
Here we identify $\mathbb{E}_{kj}$ with the constant function $\mathbb{E}%
_{kj}z^{0}\in H^{2}_{n\times n}$. Then for all $f\in H^{2}_{n\times n}$,
\[
\Delta_{n\times n}(f)=\sum_{k,j=1}^{n}\langle f,\mathbb{E}_{kj}\rangle
_{H^{2}_{n\times n}}\mathbb{E}_{kj}.
\]

Let $S_{n\times n}$ denote the shift operator on $H^{2}_{n\times n}.$ That
is,
\[
S_{n\times n}f(z)=zf(z),f\in H^{2}_{n\times n}.
\]
The adjoint operator $S_{n\times n}^{\ast}$ is given by
\[
S_{n\times n}^{\ast}f(z)=\mathbb{P}(\bar{z}f(z))=\sum_{l=0}^{\infty}%
f_{l+1}z^{l},f(z)=\sum_{l=0}^{\infty}f_{l}z^{l}\in H^{2}_{n\times n}.
\]
Thus it is easy to see that shift operator $S_{n\times n}$, its adjoint
$S_{n\times n}^{\ast}$ \ and $\Delta_{n\times n}$ satisfy the identities%
\begin{equation}
S_{n\times n}^{\ast}S_{n\times n}=I,\quad\Delta_{n\times n}=I-S_{n\times
n}S_{n\times n}^{\ast}. \label{shift operaty}%
\end{equation}
where $(\mathbb{E}_{kj}\otimes\mathbb{E}_{kj})$ denotes the rank-one operator.
For simplicity, we will also use the shorthand notation
\[
S=S_{n\times n},\qquad\Delta=\Delta_{n\times n}.
\]

The corresponding shift operators on other vector- or matrix-valued
Hardy spaces are defined in the same manner; in particular, we use
$S_{N\times1}$ for the shift on the $N$-dimensional vector-valued Hardy
space.

For later use, we introduce the right and left shifted Hankel operators
\[
K^{(r)}_{U}:=H^{(r)}_{S^{\ast}U}, \quad K^{(l)}_{U}:=H^{(l)}_{S^{\ast}U}
\]
with domains
\[
\operatorname{Dom}K_{U}^{(r)} = \operatorname{Dom}H_{S^{\ast}U}^{(r)},
\qquad\operatorname{Dom}K_{U}^{(l)} = \operatorname{Dom}H_{S^{\ast}U}^{(l)}.
\]

\subsection{Unbounded Hankel operators and self-adjoint products}

\label{auxiliary}

It is well known that, for bounded Hankel operators, the adjoint of a right
Hankel operator is again a right Hankel operator, and similarly for left
Hankel operators. However, for unbounded operators, this fact is nontrivial.
The following theorem establishes this property together with the closedness
of the operators and the density of polynomials in their domains.

Define
\begin{align*}
\operatorname{Graph} \, (\Gamma^{(r)}_{U}):=\{(f,\Gamma^{(r)}_{U}(f)):
f\in\operatorname{Dom} \Gamma^{(r)}_{U}\}.
\end{align*}
The graph norm
\[
\Vert f\Vert_{\Gamma^{(r)}_{U}}:= \left(  \Vert f\Vert^{2}_{H^{2}_{n \times
n}}+\Vert\Gamma^{(r)}_{U}(f)\Vert^{2}_{H^{2}_{n \times n}}\right)  ^{1/2}.
\]

\begin{theorem}
[Proof deferred to Section \ref{Theorem-Domain dense}]\label{Domain dense} Let
$U \in H^{2}_{n \times n}$, and consider the right Hankel operator
$\Gamma^{(r)}_{U}$ defined by (\ref{def of right Hankel}) on the domain
(\ref{Dom of right Hankel}). Then: \newline(i) $\Gamma^{(r)}_{U}$ is closed
(i.e. $\operatorname{Dom} \Gamma^{(r)}_{U}$ is closed with respect to the
graph norm of $\Gamma^{(r)}_{U}$); \newline(ii) the set of all polynomials is
dense in $\operatorname{Dom} \Gamma^{(r)}_{U}$ with respect to the graph norm;
\newline(iii) the adjoint of $\Gamma^{(r)}_{U}$ is given by $(\Gamma^{(r)}%
_{U})^{\ast}=\Gamma^{(r)}_{U^{\#}}$, where $\Gamma^{(r)}_{U^{\#}}$ is defined
as above with $U$ replaced by $U^{\#}$. That is, $(\Gamma^{(r)}_{U})^{\ast}
\subset\Gamma^{(r)}_{U^{\#}}$ and $\Gamma^{(r)}_{U^{\#}} \subset(\Gamma
^{(r)}_{U})^{\ast}$. More precisely, $(\Gamma^{(r)}_{U})^{\ast} \subset
\Gamma^{(r)}_{U^{\#}}$ means that $\operatorname{Dom}(\Gamma^{(r)}_{U})^{\ast}
\subset\text{Dom} \,\Gamma^{(r)}_{U^{\#}}$ and $(\Gamma^{(r)}_{U})^{\ast}h
=\Gamma^{(r)}_{U^{\#}}h$ for all $h\in\text{Dom}\,(\Gamma^{(r)}_{U})^{\ast}$.
\end{theorem}

Analogous results hold for the left Hankel operator.

\begin{theorem}
\label{Domain dense left}Let $U \in H^{2}_{n \times n} $, and consider the
left Hankel operator $\Gamma^{(l)}_{U}$ defined by (\ref{def of left Hankel})
on the domain (\ref{Dom of left Hankel}). Then: \newline(i) $\Gamma^{(l)}_{U}$
is closed (i.e. $\operatorname{Dom} \Gamma^{(l)}_{U}$ is closed with respect
to the graph norm of $\Gamma^{(l)}_{U}$); \newline(ii) the set of all
polynomials is dense in $\operatorname{Dom} \Gamma^{(l)}_{U}$ with respect to
the graph norm; \newline(iii) the adjoint of $\Gamma^{(l)}_{U}$ is given by
$(\Gamma^{(l)}_{U})^{\ast}=\Gamma^{(l)}_{U^{\#}}$, where $\Gamma^{(l)}%
_{U^{\#}}$ is defined as above with $U$ replaced by $U^{\#}$.
\end{theorem}

Theorems \ref{Domain dense} and \ref{Domain dense left} can be reformulated in
terms of anti-linear Hankel operators as follows.

\begin{theorem}
\label{Domain dense left and right Hankel}Let $U \in H^{2}_{n \times n}$, and
let $H^{(r)}_{U}$ be the right Hankel operator, defined by (\ref{def-right-H})
on the domain (\ref{eq:dom-left-anti-hankel}); $H^{(l)}_{U}$ be the left
Hankel operator, defined by (\ref{def-left-H}) on the domain
(\ref{eq:dom-right-anti-hankel}). Then: \newline(i) $H ^{(r)}_{U}$ and $H
^{(l)}_{U}$ are closed (i.e. $\operatorname{Dom} H^{(r)}_{U}$ and
$\operatorname{Dom} H^{(l)}_{U}$ are closed with respect to their graph norms,
respectively); \newline(ii) the set of all polynomials is dense in
$\operatorname{Dom} H^{(r)}_{U}$ and in $\operatorname{Dom} H^{(l)}_{U}$ with
respect to their graph norms, respectively; \newline(iii) assume that there
exist $f \in H^{2}_{n \times n}$ and $h \in H^{2}_{n \times n}$ such that
\begin{align*}
\langle H ^{(l)}_{U}g,f\rangle_{H^{2}_{n \times n}}=\langle h,g\rangle
_{H^{2}_{n \times n}},\quad\forall g \in\operatorname{Dom}H ^{(l)}_{U}.
\end{align*}
Then $f \in\operatorname{Dom}H ^{(r)}_{U}$ and $H ^{(r)}_{U} f= h$;
\newline(iv) assume that there exist $g\in H^{2}_{n\times n}$ and $h_{0}\in
H^{2}_{n \times n}$ such that
\[
\langle H_{U}^{(r)}f,g\rangle_{H^{2}_{n\times n}}=\langle h_{0},f\rangle
_{H^{2}_{n \times n}},\quad\forall f\in\operatorname{Dom}H_{U}^{(r)}.
\]
Then $g\in\operatorname{Dom}H_{U}^{(l)}$ and $H_{U}^{(l)}g=h_{0}$.
\end{theorem}

\begin{proof}
Parts (i) and (ii) follow immediately from Theorems \ref{Domain dense} and
\ref{Domain dense left}. The proof of part (iv) is omitted, as it is similar
to that of part (iii). Therefore, only part (iii) requires proof.

Rewriting the hypothesis of (iii) in terms of the linear Hankel operators, we
have
\begin{align}
\label{equal left}\langle\Gamma^{(l)}_{U}g^{\#} ,f\rangle_{H^{2}_{n \times n}%
}=\langle h,g\rangle_{H^{2}_{n \times n}},\quad\forall g^{\#} \in
\operatorname{Dom}\Gamma^{(l)}_{U}.
\end{align}
Moreover,
\begin{align*}
\langle h,g\rangle_{H^{2}_{n \times n}}  &  = \int_{\mathbb{T}}
\operatorname{tr}( h(z) g(z)^{\ast}) d\lambda(z)\\
&  = \int_{\mathbb{T}} \operatorname{tr}(g(z)^{\ast} h(z)) d\lambda(z)\\
&  =\int_{\mathbb{T}} \operatorname{tr}( g(\overline{z})^{\ast}h(\overline{z})
) d\lambda(z)\\
&  =\langle g^{\#} ,h^{\#}\rangle_{H^{2}_{n \times n}},
\end{align*}
so that (\ref{equal left}) becomes
\begin{align*}
\langle\Gamma^{(l)}_{U}g^{\#} ,f\rangle_{H^{2}_{n \times n}}=\langle g^{\#}
,h^{\#}\rangle_{H^{2}_{n \times n}},\quad\forall g^{\#} \in\operatorname{Dom}%
\Gamma^{(l)}_{U}.
\end{align*}
This means that $f \in\operatorname{Dom}(\Gamma^{(l)}_{U})^{\ast}$ and
$(\Gamma^{(l)}_{U})^{\ast} f=h^{\#}$. By Theorem \ref{Domain dense left}(iii),
it follows that $f \in\operatorname{Dom}\Gamma^{(l)}_{U^{\#}}$ and
$\Gamma^{(l)}_{U^{\#}} f=h^{\#}$. Applying the conjugation operator $\#$ to
this identity yields
\[
f^{\#} \in\operatorname{Dom}\Gamma^{(r)}_{U}, \quad\Gamma^{(r)}_{U}f^{\#}=h.
\]
Equivalently,
\[
f \in\operatorname{Dom}H^{(r)}_{U}, \quad H^{(r)}_{U}f=h,
\]
which proves the claim.

The proof of Theorem \ref{Domain dense left and right Hankel} is complete.
\end{proof}

To prove Lemma \ref{self-adjoint}, we recall the following classical result of
von Neumann.

\begin{lemma}
[von Neumann]\textup{\cite[Theorem X.25]{Reed802}\label{von-Neumann}}  Let
$T:\text{Dom }(T)\subseteq G\rightarrow H$ be a closed densely defined
operator between Hilbert spaces $G,H$. Then the operator $T^{\ast}T$ is
self-adjoint on its domain%
\[
\text{Dom }(T^{\ast}T)=\{x\in\text{Dom }(T):Tx\in\text{Dom }(T^{\ast})\},
\]
and is also densely defined in $G$.
\end{lemma}

We can now establish the self-adjointness of the product of the
anti-linear Hankel operators $H_{U}^{(r)}H_{U}^{(l)}$.

\begin{lemma}
\label{self-adjoint}  For $U\in H^{2}_{n\times n}$, the operator
\[
H_{U}^{(r)}H_{U}^{(l)}
\]
is a densely defined self-adjoint operator with domain
\[
\operatorname{Dom}(H_{U}^{(r)}H_{U}^{(l)}) = \{g\in\operatorname{Dom}
H_{U}^{(l)}: H_{U}^{(l)}g\in\operatorname{Dom}H_{U}^{(r)} \}.
\]

\end{lemma}

\begin{proof}
By Theorem \ref{Domain dense}(i)-(iii), the linear right Hankel operator
$\Gamma_{U}^{(r)}$ is densely defined and closed, and the adjoint of
$\Gamma^{(r)}_{U}$ is the Hankel operator $\Gamma^{(r)}_{U^{\#}}$, that is,
\[
(\Gamma^{(r)}_{U})^{\ast}=\Gamma^{(r)}_{U^{\#}}
\]
as operators with the same domain. Moreover, by the definitions of
$\mathcal{C}$ and $U^{\#}=JU^{\ast}$, we have
\[
\mathcal{C}\Gamma_{U}^{(l)}\mathcal{C }= \Gamma_{U^{\#}}^{(r)}.
\]
Consequently,
\[
H^{(r)}_{U}H ^{(l)}_{U}= \Gamma^{(r)}_{U}\mathcal{C}\Gamma^{(l)}%
_{U}\mathcal{C}= \Gamma^{(r)}_{U}\Gamma^{(r)}_{U^{\#}} =\Gamma^{(r)}%
_{U}(\Gamma^{(r)}_{U})^{\ast}.
\]
Therefore, it suffices to study the operator $\Gamma_{U}^{(r)}(\Gamma
_{U}^{(r)})^{\ast}$.

By von Neumann's theorem (Lemma \ref{von-Neumann}), $\Gamma_{U}^{(r)}%
(\Gamma_{U}^{(r)})^{\ast}$ is a densely defined self-adjoint operator
with
\[
\operatorname{Dom} \Gamma_{U}^{(r)}(\Gamma_{U}^{(r)})^{\ast} = \left\{
g\in\operatorname{Dom}(\Gamma_{U}^{(r)})^{\ast}: (\Gamma_{U}^{(r)})^{\ast}g
\in\operatorname{Dom}\Gamma_{U}^{(r)} \right\} .
\]
By the operator identity above,
\[
\operatorname{Dom}(H_{U}^{(r)}H_{U}^{(l)}) = \left\{  g\in\operatorname{Dom}%
H_{U}^{(l)}: H_{U}^{(l)}g\in\operatorname{Dom}H_{U}^{(r)} \right\} .
\]
Therefore, $H_{U}^{(r)}H_{U}^{(l)}$ is also a densely defined self-adjoint
operator.
\end{proof}

To establish the strong resolvent convergence of $H^{(r)}_{U}H^{(l)}_{U}$,
which will be given in Theorem \ref{strong con-H}, we first recall
the notion of strong resolvent convergence.

\begin{definition}
\textup{\cite[Theorem VIII.19]{Reed80} } Let $\{A_{k}\}_{k\in\mathbb{N}}$ be a
sequence of self-adjoint operators on the Hilbert space $X$, and let $A$ be a
self-adjoint operator on $X$. The sequence $\{A_{k}\}$ is said to converge to
$A$ in the strong resolvent sense if for all $\lambda\in\mathbb{C}$ with
$\operatorname{Im}\lambda\neq0$, the resolvent sequence $\{(A_{k}%
-\lambda\,I)^{-1}\}$ converges to $(A-\lambda\,I)^{-1}$ in the strong operator
topology on $B(X)$; that is, for every $f\in X$,
\[
\left\Vert (A_{k}-\lambda\,I)^{-1}f-(A-\lambda\,I)^{-1}f\right\Vert
_{X}\rightarrow0\quad\text{as }k\rightarrow\infty.
\]

\end{definition}

The following standard consequence of strong resolvent convergence will be
used below.

\begin{lemma}
\textup{\cite[Theorem VIII.20]{Reed80}} \label{Reed80Thm} Let $\{A_{k}%
\}_{k\in\mathbb{N}}$ and $A$ be self-adjoint operators on the Hilbert space
$X$. If $A_{k}\rightarrow A$ in the strong resolvent sense and $\varphi$ is a
bounded continuous function, then
\[
\varphi(A_{k})f\rightarrow\varphi(A)f
\]
for all $f \in X$.
\end{lemma}

We now apply this framework to anti-linear Hankel operators and establish the
following strong resolvent convergence result.

\begin{theorem}
[Proof deferred to Section \ref{proof of strong convergence}]%
\label{strong con-H} Let $\{U_{k}\}_{k=1}^{\infty}$ be a sequence of elements
of $H^{2}_{n\times n}$ with $\Vert U_{k}-U\Vert_{H^{2}_{n\times n}}%
\rightarrow0$ as $k\rightarrow\infty$. Then
\[
H_{U_{k}}^{(r)}H_{U_{k}}^{(l)}\rightarrow H_{U}^{(r)}H_{U}^{(l)}%
\]
in the strong resolvent sense.
\end{theorem}

\subsection{Solutions of cubic matrix Szeg\H{o} equation}

\label{sun result}

In this subsection, we recall two fundamental results on the cubic matrix
Szeg\H{o} equation (\ref{cubic matrix}) from \cite{Sun}: the Lax pair
formulation (Theorem \ref{Sun lax pair}) and the explicit formula for the flow
map (Theorem \ref{Sun explicit formula}). The latter is stated for initial
data in $\mathcal{H}_{+}^{s}( \mathbb{T}, \mathbb{C}^{n\times n}),\;s>1/2 $.
In Section \ref{Auxiliary results}, this explicit formula will be extended to
the larger space $H^{2}_{n\times n}$, where it will play a central role in the
proof of our main result.

For the reader's convenience, we provide detailed proofs of these two
theorems in Appendix \ref{Alternative proof}. Although these results are
covered by the more general framework of \cite{Sun}, our more specialized
setting involves simpler notation, allowing for a more direct proof.
The proof of the Lax pair identities is obtained by adapting the
arguments used for the cubic scalar Szeg\H{o} equation (\ref{cubic scalar}) in
\cite{Gerard15} to the matrix-valued setting, while the proof of Sun's
explicit formula follows the argument in \cite{Sun}, rewritten in our notation
and setting.

\begin{remark}
The results in \cite{Sun} are established for arbitrary dimensions
$M,N,d\in\mathbb{N}_{+}$. Throughout this paper, we work with the square case
$M=d=N=n$, which is sufficient for our purposes.
\end{remark}

In the following theorem, let $U=U(t,z)\in C^{\infty}(\mathbb{R}%
\times\mathbb{T}, \mathcal{H}_{+}^{s}( \mathbb{T}, \mathbb{C}^{n\times n}))$
with $s> 1/2$. The anti-linear Hankel operators introduced in \cite{Sun} are
bounded on $H^{2}_{n\times n}$, and hence have no domain issues. Moreover,
they coincide with the anti-linear Hankel operators defined above.

\begin{theorem}
\textup{\cite[Theorem 1.3]{Sun}}\label{Sun lax pair} Let $U=U(t,z)\in
C^{\infty}(\mathbb{R} \times\mathbb{T},\mathcal{H}^{s}_{+}(\mathbb{T}%
,\mathbb{C}^{n \times n})) $, $s>1/2 $, be a solution of (\ref{cubic matrix}).
Then
\begin{align}
\frac{d}{dt} (H^{(r)}_{U}H^{(l)}_{U})=i[H^{(r)}_{U}H^{(l)}_{U},T^{(r)}%
_{UU^{\ast}}];\label{Lax-pair}\\
\frac{d}{dt}(H^{(l)}_{U}H^{(r)}_{U})=i[H^{(l)}_{U}H^{(r)}_{U},T^{(l)}%
_{U^{\ast}U}];\nonumber\\
\frac{d}{dt}(K^{(r)}_{U}K^{(l)}_{U})=i[K^{(r)}_{U}K^{(l)}_{U},T^{(r)}%
_{UU^{\ast}}];\nonumber\\
\frac{d}{dt}(K^{(l)}_{U}K^{(r)}_{U})=i[K^{(l)}_{U}K^{(r)}_{U},T^{(l)}%
_{U^{\ast}U}].\nonumber
\end{align}

\end{theorem}

\begin{theorem}
\textup{\cite[Theorem 1.6]{Sun}} \label{Sun explicit formula} Let $U = U(0,z)
\in\mathcal{H}^{s}_{+}(\mathbb{T},\mathbb{C}^{n \times n}),\, s>1/2 $, and
let  $U(t,z) \in C(\mathbb{R} , \allowbreak \mathcal{H}^{s} _{+}(\mathbb{T},\mathbb{C}^{n
\times n}))$ be the solution of equation  (\ref{cubic matrix}) with initial
data $U$. For $|z|<1,$ we have
\begin{align*}
\underline{U}(t,z)  &  =\Delta\left(  (I-ze^{-itH^{(r)} _{U}H^{(l)}_{U}%
}e^{itK^{(r)}_{U}K^{(l)}_{U}}S^{\ast} )^{-1}e^{-itH^{(r)}_{U}H^{(l)}_{U}%
}(U)\right) \\
&  =\Delta\left(  (I-ze^{-itH^{(l)}_{U}H^{(r)}_{U} }e^{itK^{(l)}_{U}%
K^{(r)}_{U}}S^{\ast})^{-1}e^{-itH^{(l)}_{U} H^{(r)}_{U}}(U)\right)  ,
\end{align*}
where $\underline{U}(t,z)$ denotes the Poisson integral of
$U(t,\cdot)$, which is holomorphic on the unit disk $\mathbb{D}$.
\end{theorem}

\subsection{Extension of the flow map to \texorpdfstring{$H^2_{n\times n}$}{}}

\label{Extension flow}

By Lemma \ref{self-adjoint}, applied to $U$ and $S^{\ast}U$, the operators
\[
H_{U}^{(r)}H_{U}^{(l)} \quad\text{and}\quad K_{U}^{(r)}K_{U}^{(l)}
\]
are densely defined self-adjoint operators. Hence, by Stone's theorem
\cite[Theorem X.5.6]{Conway},
\[
\{e^{-itH_{U}^{(r)}H_{U}^{(l)}}\}_{t\in\mathbb{R}} \quad\text{and}%
\quad\{e^{itK_{U}^{(r)}K_{U}^{(l)}}\}_{t\in\mathbb{R}}
\]
are unitary groups on $H^{2}_{n\times n}$ . Similarly, by the same argument,
together with Theorem \ref{Domain dense left}, the operators
\[
H_{U}^{(l)}H_{U}^{(r)} \quad\text{and}\quad K_{U}^{(l)}K_{U}^{(r)}
\]
are also densely defined self-adjoint operators. Moreover, since $S^{\ast}$ is
a contraction and the exponential operators are unitary, both
\[
e^{-itH_{U}^{(r)}H_{U}^{(l)}}e^{itK_{U}^{(r)}K_{U}^{(l)}}S^{\ast}
\quad\text{and}\quad e^{-itH_{U}^{(l)}H_{U}^{(r)}}e^{itK_{U}^{(l)}K_{U}^{(r)}%
}S^{\ast}
\]
are contractions. Therefore, $(\mathcal{W}(t)U)(z)$ is well-defined and
holomorphic for $|z|<1$.

This enables us to extend the formula in Theorem \ref{Sun formula}, initially
defined for $U\in\mathcal{H}_{+}^{s}(\mathbb{T},\mathbb{C}^{n\times n}), \quad
s>1/2$ to all $U\in H^{2}_{n\times n}.$ This extension constitutes one of the
key parts in the proof of our main result.

For the cubic matrix Szeg\H{o} equation (\ref{cubic matrix}) with
initial data $U(0,z)=U\in H^{2}(\mathbb{T},\mathbb{C}^{n\times n})$, we denote
the corresponding extended solution at time $t$ by $\mathcal{W}(t)U$. We
state this formally as the following theorem.

Furthermore, we can regard $(W(t)U)(z)$ as an element of
$H^{2}(\mathbb{D},C^{n\times n})$, i.e. as a matrix-valued analytic function
on the unit disk rather than merely a boundary function on the unit circle.

\begin{theorem}
\label{Sun formula} The flow map $\mathcal{W}(t)$ of the cubic matrix
Szeg\H{o} equation (\ref{cubic matrix}) can be extended to the whole
$H_{n\times n}^{2}.$ More precisely, for initial data $U\in H_{n\times n}^{2}%
$,
\begin{align}
(\mathcal{W}(t)U)(z) &  =\Delta\left(  (I-ze^{-itH_{U}^{(l)}H_{U}^{(r)}%
}e^{itK_{U}^{(l)}K_{U}^{(r)}}S^{\ast})^{-1}e^{-itH_{U}^{(l)}H_{U}^{(r)}%
}(U)\right)  \label{flow}\\
&  =\Delta\left(  (I-ze^{-itH_{U}^{(r)}H_{U}^{(l)}}e^{itK_{U}^{(r)}K_{U}%
^{(l)}}S^{\ast})^{-1}e^{-itH_{U}^{(r)}H_{U}^{(l)}}(U)\right)  \nonumber
\end{align}
where $z\in\mathbb{D}$.
\end{theorem}

We remark that formula (\ref{flow}) in Theorem \ref{Sun formula}
serves as the starting point for the proof of Theorem \ref{main result}. We
next introduce a pair of isometries that play a central role in the proof.

\subsection{Wold decomposition and the pure part of the isometry
\texorpdfstring{$\Omega^{(l)}$}{}}

\label{pure part sec}

An isometry $\Omega$ on a Hilbert space is said to be pure if $\Omega$ is not
unitary on any of its invariant subspaces. The dimension of
$\mathrm{Ran}(I-\Omega\Omega^{\ast})$ is called the multiplicity of
$\Omega.$ The following result is the Wold decomposition of an
isometry (see, e.g.,  \cite{Halmos61}, \cite[Theorem I.1.1]{Nagy10}) in the
case the pure part is of finite multiplicity.

\begin{theorem}[Wold decomposition]
\label{Wold decomposition}  Let $\Omega$ be an arbitrary
isometry on a Hilbert space $X$ with $\dim\mathrm{Ran} \allowbreak
(I-\Omega\Omega^{\ast})=j.$ Then $X$ can be represented as an orthogonal sum
$X=X_{\mathrm{u}}\oplus X_{\mathrm{pure}}$, s.t.
\begin{equation}
\label{Wold decomposition u part}\Omega=%
\begin{pmatrix}
\Omega_{\mathrm{u}} & 0\\
0 & \Omega_{\mathrm{pure}}%
\end{pmatrix}
\quad\text{on} \quad X_{\mathrm{u}}\oplus X_{\mathrm{pure}},
\end{equation}
where $\Omega_{\mathrm{u}}$ is unitary and $\Omega_{\mathrm{pure}}$ is a pure
isometry. Moreover, let $\{q_{k}:1\leq k\leq j\}$ be an orthonormal basis of
$\mathrm{Ran}(I-\Omega\Omega^{\ast})$. Then $X_{\mathrm{pure}}$ coincides with
the closed linear span of $\{\Omega^{l}q_{k}:1\leq k\leq j,l\geq0\}$ and
$\Omega_{\mathrm{pure}}$ is unitarily equivalent to the shift operator
$S_{j\times1}$ on $H^{2}(\mathbb{T},\mathbb{C}^{j\times1})$. In
particular, $\{\Omega^{m}q_{k}:1\leq k\leq j,m\geq0\}$ forms an orthonormal
basis of $X_{\mathrm{pure}}$.
\end{theorem}

From now on, for $U\in H^{2}_{n\times n}$ and $t\in\mathbb{R}$. Set
\begin{align*}
\Omega^{(r)}  &  :=e^{itH_{U}^{(r)}H_{U}^{(l)}}S e^{-itK_{U}^{(r)}K_{U}^{(l)}%
}:H^{2}_{n\times n}\to H^{2}_{n\times n},\\
Q^{(r)}  &  :=\{q^{(r)}_{kj}:1\leq k, j\leq n\},\text{ where}\quad
q^{(r)}_{kj}=e^{itH_{U}^{(r)}H_{U}^{(l)}}\mathbb{E}_{kj}.
\end{align*}

Similarly, let
\begin{align*}
\Omega^{(l)}  &  := e^{it H^{(l)}_{U}H^{(r)}_{U}}Se^{-it K^{(l)}_{U}%
K^{(r)}_{U}} :H^{2}_{n\times n}\to H^{2}_{n\times n},\\
Q^{(l)} &  :=\{q^{(l)}_{kj}:1\leq k, j\leq n\},\text{ where}\quad q^{(l)}%
_{kj}=e^{itH_{U}^{(l)}H_{U}^{(r)}}\mathbb{E}_{kj} .
\end{align*}
It is clear that $\Omega^{(r)}$ and $\Omega^{(l)}$ are isometries on
$H^{2}_{n\times n}$, since they are the product of unitary operators and the
shift isometry $S$.

Directly computing, thus,%
\begin{align}
I-\Omega^{(l)}\left( \Omega^{(l)}\right) ^{\ast}  &  =I-e^{itH_{U}^{(l)}%
H_{U}^{(r)}}S S^{\ast}e^{-itH_{U}^{(l)}H_{U}^{(r)}}= I-e^{itH_{U}^{(l)}%
H_{U}^{(r)}}\left(  I -\Delta\right)  e^{-itH_{U}^{(l)}H_{U}^{(r)}}\nonumber\\
&  =e^{itH_{U}^{(l)}H_{U}^{(r)}}\Delta e^{-itH_{U}^{(l)}H_{U}^{(r)}%
}\nonumber\\
&  =\sum_{k,j=1}^{n}(e^{itH_{U}^{(l)}H_{U}^{(r)}}\mathbb{E}_{kj} \otimes
e^{itH_{U}^{(l)}H_{U}^{(r)}} \mathbb{E}_{kj} )\nonumber\\
&  =\sum_{k,j}q^{(l)}_{kj}\otimes q^{(l)}_{kj} \label{finite rank}%
\end{align}
and $\Vert q^{(l)}_{kj}\Vert_{H^{2}_{n\times n}}=1$ (compare with
(\ref{shift operaty})).

We apply the Wold decomposition to the isometry $\Omega^{(l)}$ and write
$H^{2}_{n\times n}=X^{(l)}_{\mathrm{u}}\oplus X^{(l)}_{\mathrm{pure}}$, where
$\Omega^{(l)} =\Omega^{(l)}_{\mathrm{u}}\oplus\Omega^{(l)}_{\mathrm{pure}} .$
It follows from (\ref{finite rank}) that%
\begin{align*}
\mathrm{Ran}\left( I-\Omega^{(l)}\left( \Omega^{(l)}\right) ^{\ast}\right)
&  =\operatorname{span} Q^{(l)},\\
X^{(l)}_{\mathrm{pure}}  &  =%
{\textstyle\bigoplus_{n=0}^{\infty}}
\left( \Omega^{(l)}\right) ^{m}\mathrm{Ran}\left( I-\Omega^{(l)}\left(
\Omega^{(l)}\right) ^{\ast}\right) =\overline{\left\{ \left( \Omega
^{(l)}\right) ^{m}Q^{(l)}:m\ge0\right\}  }.
\end{align*}
Hence $\Omega^{(l)}_{\mathrm{pure}}$ is unitarily equivalent to $S$ on
$H^{2}_{n\times n}.$

Similarly, for $\Omega^{(r)}$, we also have
\begin{align}
\label{rpure}X^{(r)}_{\mathrm{pure}}  &  =%
{\textstyle\bigoplus_{n=0}^{\infty}}
\left( \Omega^{(r)}\right) ^{m}\mathrm{Ran}\left( I-\Omega^{(r)}\left(
\Omega^{(r)}\right) ^{\ast}\right) =\overline{\left\{ \left( \Omega
^{(r)}\right) ^{m}Q^{(r)}:m\ge0\right\}  }.
\end{align}

The following theorem, which shows that $U$ belongs to the pure part of the
isometry $\Omega^{(l)}$, will be used in the proof of Theorem
\ref{main result}.

\begin{theorem}
[Proof deferred to Section \ref{start}]\label{complicated} Assume that $U\in
H^{2}_{n\times n}$, then $\Omega^{(l)}$ is an isometry and
\[
U\in%
{\textstyle\bigoplus_{n=0}^{\infty}}
\left( \Omega^{(l)}\right) ^{m}\mathrm{Ran}\left( I-\Omega^{(l)}\left(
\Omega^{(l)}\right) ^{\ast}\right) =\overline{\left\{ \left( \Omega
^{(l)}\right) ^{m}Q^{(l)}:m\ge0\right\}  },
\]
where $\overline{\left\{ \left( \Omega^{(l)}\right) ^{m}Q^{(l)}:m\ge0\right\}
}$ denotes the closed linear span of $\left\{ \left( \Omega^{(l)}\right)
^{m}Q^{(l)}:m\ge0\right\}  $.
\end{theorem}

\section{Proof of Theorem \ref{main result}}

\label{Proof of main} This section proves Theorem \ref{main result}. Part (i)
contains the core argument, while parts (ii)-(iv) build on the
auxiliary results stated below. However, the proof of Theorem
\ref{main result} (i) presents the main line of argument, which relies on four
auxiliary theorems in Subsection \ref{Auxiliary results}.

\label{proof of main result}

\subsection{Proof of Theorem \ref{main result}(i)}

\label{1sec}

Before proving Theorem \ref{main result}(i), we first establish one lemma on
contractions $\Omega^{(l)}$ and $\left( \Omega^{(l)}\right) ^{\ast}$.

\begin{lemma}
\label{unitary equivalent} Recall the definition of $\Omega^{(l)}$.
Then, for any $|z|<1$, we have
\begin{align*}
&  \left( I - z e^{-itH_{U}^{(l)}H_{U}^{(r)}} e^{itK_{U}^{(l)}K_{U}^{(r)}} S
^{\ast}\right) ^{-1} e^{-itH_{U}^{(l)}H_{U}^{(r)}}\\
=\,  &  e^{-itH_{U}^{(l)}H_{U}^{(r)}} \left( I - z \left( \Omega^{(l)}\right)
^{\ast}\right) ^{-1}.
\end{align*}

\end{lemma}

\begin{proof}
Set $\mathcal{U} := e^{-itH_{U}^{(l)}H_{U}^{(r)}}$, which is unitary. For
convenience, define
\[
A := e^{itK_{U}^{(l)}K_{U}^{(r)}} S ^{\ast} e^{-itH_{U}^{(l)}H_{U}^{(r)}},
\qquad B := e^{-itH_{U}^{(l)}H_{U}^{(r)}} e^{itK_{U}^{(l)}K_{U}^{(r)}} S
^{\ast} = \mathcal{U} A \mathcal{U}^{\ast}.
\]
Note that
\begin{align*}
A  & = \left( \Omega^{(l)}\right) ^{\ast} =\left( e^{it H^{(l)}_{U}H^{(r)}%
_{U}}Se^{-it K^{(l)}_{U}K^{(r)}_{U}}\right) ^{\ast}\\
& = e^{itK_{U}^{(l)}K_{U}^{(r)}} S ^{\ast} e^{-itH_{U}^{(l)}H_{U}^{(r)}} .
\end{align*}

The lemma then follows immediately from the standard identity for unitary
equivalence of resolvents:
\[
(I - z B)^{-1} \mathcal{U} = \mathcal{U} (I - z A)^{-1}.
\]

\end{proof}

We are now in position to prove Theorem \ref{main result}(i).

\textit{\textbf{Proof of Theorem \ref{main result}(i):}} By (\ref{flow})\ and
Lemma \ref{unitary equivalent}, for a fixed $z\in\mathbb{D}$, we have%
\begin{align*}
&  (\mathcal{W}(t)U)(z)=\Delta\left(  I-ze^{-itH_{U}^{(l)}H_{U}^{(r)}%
}e^{itK_{U}^{(l)}K_{U}^{(r)}}S^{\ast})^{-1}e^{-itH_{U}^{(l)}H_{U}^{(r)}%
}(U)\right) \\
&  =\Delta\left(  e^{-itH_{U}^{(l)}H_{U}^{(r)}}\left( I-z\left( \Omega
^{(l)}\right) ^{\ast}\right) ^{-1}U\right) \\
&  = \sum_{k,j=1}^{n} \left\langle e^{-itH_{U}^{(l)}H_{U}^{(r)}}\left(
I-z\left( \Omega^{(l)}\right) ^{\ast}\right) ^{-1}U, \mathbb{E}_{kj}
\right\rangle \mathbb{E}_{kj}\\
&  =\sum_{k,j}\left\langle \left( I-z\left( \Omega^{(l)}\right) ^{\ast}\right)
^{-1}U,q^{(l)}_{kj}\right\rangle \mathbb{E}_{kj}\\
&  =\sum_{k,j}\left\langle \left( \sum_{m=0}^{\infty}z^{m} \left( \left(
\Omega^{(l)}\right) ^{\ast}\right) ^{m}\right) U,q^{(l)}_{kj} \right\rangle
\mathbb{E}_{kj}\\
&  =\sum_{m=0}^{\infty}z^{m} \sum_{k,j} \left\langle U,\left( \Omega
^{(l)}\right) ^{m}q^{(l)}_{kj}\right\rangle _{H^{2}_{n\times n}}%
\mathbb{E}_{kj}.
\end{align*}

By Theorem \ref{complicated}, $U$ belongs to the closed subspace
$\overline{\left\{ \left( \Omega^{(l)}\right) ^{m}Q^{(l)}:m\ge0\right\}  }$.
Furthermore, according to the Wold decomposition (Theorem
\ref{Wold decomposition}), 
\[
\left\{ \left( \Omega^{(l)}\right) ^{m}Q^{(l)}:m\ge0\right\}
\]
is an orthonormal basis of this closed subspace. Consequently, by Parseval's
identity,
\[
\Vert U\Vert_{H^{2}_{n\times n}}^{2}  = \sum_{m=0}^{\infty}\sum_{k,j}
\left\vert \left\langle U,\left( \Omega^{(l)}\right) ^{m}q^{(l)}%
_{kj}\right\rangle _{H^{2}_{n\times n}}\right\vert ^{2}.
\]

On the other hand, since $\{\mathbb{E}_{kj}z^{m}:1\leq k,j\leq n,m\geq0\}$
forms an orthonormal basis for $H^{2}_{n\times n}$, we have
\begin{align*}
\Vert\mathcal{W}(t)U\Vert^{2}_{H^{2}_{n\times n}}  & =\left\Vert \sum
_{m=0}^{\infty} \sum_{k,j} \left\langle U,\left( \Omega^{(l)}\right)
^{m}q^{(l)}_{kj}\right\rangle _{H^{2}_{n\times n}} \mathbb{E}_{kj} z^{m}
\right\Vert ^{2}_{H^{2}_{n\times n}}\\
& = \sum_{m=0}^{\infty}\sum_{k,j} \left\vert \left\langle U,\left(
\Omega^{(l)}\right) ^{m}q^{(l)}_{kj}\right\rangle _{H^{2}_{n\times n}}
\right\vert ^{2} \Vert\mathbb{E}_{kj} z^{m} \Vert^{2}_{H^{2}_{n\times n}}\\
& = \Vert U\Vert_{H^{2}_{n\times n}}^{2} .
\end{align*}

\subsection{Proof of Theorem \ref{main result} (ii) and (iii)}

\label{23sec}

\textit{\textbf{Proof of Theorem \ref{main result}(ii):}} Let $\Vert
U_{k}-U\Vert_{H^{2}_{n\times n}}\rightarrow0$ in $H^{2}_{n\times n}$. By
Theorem \ref{strong con-H} we have the strong resolvent convergence
\[
H_{U_{k}}^{(r)}H_{U_{k}}^{(l)}\rightarrow H_{U}^{(r)}H_{U}^{(l)}.
\]
It then follows from Lemma \ref{Reed80Thm} that, for every bounded continuous
function $\varphi$ on $\mathbb{R}$,
\[
\varphi(H_{U_{k}}^{(r)}H_{U_{k}}^{(l)})\rightarrow\varphi(H_{U}^{(r)}%
H_{U}^{(l)})
\]
in the strong operator topology. In particular, taking
\[
\varphi_{t}(x)=\exp(-itx),\quad x\in\mathbb{R},
\]
we obtain
\[
\exp(-itH_{U_{k}}^{(r)}H_{U_{k}}^{(l)})\rightarrow\exp(-itH_{U}^{(r)}%
H_{U}^{(l)})
\]
in the strong operator topology as $k\rightarrow\infty$.

Similarly, we have
\[
\exp(-itK_{U_{k}}^{(r)}K_{U_{k}}^{(l)})\rightarrow\exp(-itK_{U}^{(r)}%
K_{U}^{(l)})
\]
in the strong operator topology as $k\rightarrow\infty$. This yields pointwise
convergence $(\mathcal{W}(t)U_{k})(z)\rightarrow(\mathcal{W}(t)U)(z) $ for all
$|z|<1$, which further implies that
\[
\mathcal{W}(t) U_{k}\rightarrow\mathcal{W}(t)U
\]
weakly in $H^{2}_{n\times n}$. Moreover, Theorem \ref{main result}(i) yields
the conservation of the $H^{2}_{n\times n}$-norm:
\[
\Vert\mathcal{W}(t)U_{k}\Vert_{H^{2}_{n\times n}}=\Vert U_{k}\Vert
_{H^{2}_{n\times n}}\rightarrow\Vert U\Vert_{H^{2}_{n\times n}}=\Vert
\mathcal{W}(t)U\Vert_{H^{2}_{n\times n}},
\]
and so we conclude that
\[
\Vert\mathcal{W}(t)U_{k}-\mathcal{W}(t)U\Vert_{H^{2}_{n\times n}}\rightarrow0
\]
as $k\rightarrow\infty$.

Next, we will prove Theorem \ref{main result}(iii) as follows.

\textit{\textbf{Proof of Theorem \ref{main result}(iii):}} Let $t_{k}%
\rightarrow t$ as $k\rightarrow\infty$. Following the argument from the proof
of part (ii), the strong convergence of the associated operators implies that
$(\mathcal{W}(t_{k})U)(z)\rightarrow(\mathcal{W}(t)U)(z)$ for any $|z|<1.$
Hence, $\mathcal{W}(t_{k})U\rightarrow\mathcal{W}(t)U$ weakly in
$H^{2}_{n\times n}$. Moreover, by Theorem \ref{main result}(i),
\[
\Vert\mathcal{W}(t_{k})U\Vert_{H^{2}_{n\times n}}=\Vert\mathcal{W}%
(t)U\Vert_{H^{2}_{n\times n}},
\]
we conclude that
\[
\Vert\mathcal{W}(t_{k})U-\mathcal{W}(t)U\Vert_{H^{2}_{n\times n}}\rightarrow0
\]
as $k\rightarrow\infty$.

\subsection{Proof of Theorem \ref{main result}(iv)}

\label{4sec}

In the study of the matrix Szeg\H{o} equation (\ref{cubic matrix}),
conservation laws play an important role. In particular, the following lemmas
provide the conserved quantities that will be used in the proof of Theorem
\ref{main result}(iv). These results extend the conservation laws
established in the scalar case \cite{Gerard10} to the matrix-valued setting.

\begin{lemma}
\label{conservation law} Assume that $U_{0} \in\mathcal{H}_{+}^{s}%
(\mathbb{T},\mathbb{C}^{n\times n})$ for $s > 1/2$, and let $U(t) =
\mathcal{W}(t)U_{0} \in C(\mathbb{R}, \mathcal{H}_{+}^{s}(\mathbb{T}%
,\mathbb{C}^{n\times n}))$ be the corresponding solution of
(\ref{cubic matrix}). For every $\alpha\in\mathbb{N}^{+}$, set
\[
J_{\alpha}(U(t)) = \sum_{k,j=1}^{n} \left\langle \left( H_{U(t)}^{(r)}%
H_{U(t)}^{(l)}\right) ^{\alpha} \mathbb{E}_{kj}, \mathbb{E}_{kj} \right\rangle
_{H^{2}_{n\times n}} .
\]
Then
\[
\frac{d}{dt}J_{\alpha}(U(t))=0, \qquad\forall\alpha\in\mathbb{N}^{+}.
\]

\end{lemma}

\begin{proof}
For simplicity, we write $U$ for $U(t)$ throughout the proof. Let
$\mathcal{L}:=H^{(r)}_{U }H^{(l)}_{U }$, by (\ref{Lax-pair}),
\begin{align*}
\frac{d}{dt} \mathcal{L}=i\left[  \mathcal{L},T^{(r)}_{U U ^{\ast}}\right]  .
\end{align*}
By induction on $\alpha$, we have
\[
\frac{d}{dt} \mathcal{L}^{\alpha}= i [\mathcal{L}^{\alpha}, T^{(r)}_{UU^{\ast
}}].
\]

Using the definition of $J_{\alpha}(U)$, we have:
\begin{align*}
\frac{d}{dt} J_{\alpha}(U)  &  = \sum_{k,j=1}^{n} \left\langle \frac{d}{dt}
\mathcal{L}^{\alpha}\mathbb{E}_{kj}, \mathbb{E}_{kj} \right\rangle \\
&  = i \sum_{k,j} \left\langle \left(  \mathcal{L}^{\alpha}T^{(r)}_{UU^{\ast}}
- T^{(r)}_{UU^{\ast}} \mathcal{L}^{\alpha}\right)  \mathbb{E}_{kj},
\mathbb{E}_{kj} \right\rangle \\
&  = i \sum_{k,j} \left(  \left\langle \mathcal{L}^{\alpha}T^{(r)}_{UU^{\ast}}
\mathbb{E}_{kj}, \mathbb{E}_{kj} \right\rangle - \left\langle T^{(r)}%
_{UU^{\ast}} \mathcal{L}^{\alpha}\mathbb{E}_{kj}, \mathbb{E}_{kj}
\right\rangle \right)  .
\end{align*}

We next identify the action of $T^{(r)}_{UU^{\ast}}$ on the constant matrix
$\mathbb{E}_{kj}$. Indeed, $\mathcal{L} \mathbb{E}_{kj}=H^{(r)}_{U}H^{(l)}%
_{U}\mathbb{E}_{kj}=\mathbb{P}\left(  U \left(  \left( \mathbb{E}_{kj}\right)
^{\ast}U \right)  ^{\ast}\right)  = \mathbb{P}\left(  U U^{\ast}%
\mathbb{E}_{kj}\right)  =T^{(r)}_{UU^{\ast}}\mathbb{E}_{kj}$. Consequently,
\[
\left\langle \mathcal{L}^{\alpha}T^{(r)}_{UU^{\ast}}\mathbb{E}_{kj},
\mathbb{E}_{kj} \right\rangle = \left\langle \mathcal{L}^{\alpha+1}%
\mathbb{E}_{kj}, \mathbb{E}_{kj} \right\rangle .
\]
On the other hand, since $\mathcal{L}$ and $T^{(r)}_{UU^{\ast}}$ are
self-adjoint operators on $H^{2}_{n \times n}$,
\begin{align*}
\left\langle T^{(r)}_{UU^{\ast}}\mathcal{L}^{\alpha}\mathbb{E}_{kj},
\mathbb{E}_{kj} \right\rangle  & = \left\langle \mathcal{L}^{\alpha}%
\mathbb{E}_{kj}, T^{(r)}_{UU^{\ast}}\mathbb{E}_{kj} \right\rangle \\
& = \left\langle \mathcal{L}^{\alpha}\mathbb{E}_{kj}, \mathcal{L}%
\mathbb{E}_{kj} \right\rangle \\
& = \left\langle \mathcal{L}^{\alpha+1}\mathbb{E}_{kj}, \mathbb{E}_{kj}
\right\rangle .
\end{align*}
Therefore
\[
\frac{d}{dt}J_{\alpha}(U)=0.
\]

\end{proof}

The following lemma provides another conserved quantity of the flow, needed
for Theorem \ref{main result}(iv).

\begin{lemma}
\label{conservation-Q} Let $U_{0} \in\mathcal{H}_{+}^{s}(\mathbb{T}%
,\mathbb{C}^{n\times n})$ for $s > 1/2$, and let $U(t) = \mathcal{W}(t)U_{0}
\in C(\mathbb{R}, \mathcal{H}_{+}^{s}(\mathbb{T},\mathbb{C}^{n\times n}))$ be
the corresponding solution of (\ref{cubic matrix}). Set
\[
\mathcal{Q}(U(t)) = \int_{\mathbb{T}} U(t,z) U(t,z)^{\ast} d\lambda(z).
\]
Then $\mathcal{Q}(U(t))$ is conserved along the flow of (\ref{cubic matrix}),
that is,
\[
\frac{d}{dt}\mathcal{Q}(U(t)) = 0.
\]

\end{lemma}

\begin{proof}
For simplicity, we write $U$ for $U(t)$. If we write $\Phi\in L^{2}_{n\times
n}$ as $\Phi=\Phi_{+}+\Phi_{-}$, where $\Phi_{+}=\mathbb{P}(\Phi)\in
H^{2}_{n\times n},$ $\Phi_{-}=(I-\mathbb{P})(\Phi)\in\overline{zH^{2}_{n\times
n}},$ and since $U^{\ast}\in\overline{H^{2}_{n\times n}}$ then
\[
\Phi_{-}U^{\ast}\in\overline{zH^{1}_{n\times n}},
\]
has no 0-th Fourier coefficient. Hence its integral over $\mathbb{T}$
vanishes:
\[
\int_{ \mathbb{T}}\Phi_{-}U^{\ast}\,d\lambda(z)=0.
\]
This yields
\begin{align}
\label{pro1}\int_{ \mathbb{T}}\mathbb{P}(\Phi)U^{\ast}\,d\lambda(z) = \int_{
\mathbb{T}}\Phi\,U^{\ast}\,d\lambda(z).
\end{align}
Similarly, $U(\Phi_{-})^{\ast}\in zH^{1}_{n\times n} $ has no 0-th Fourier
coefficient, so
\[
\int_{ \mathbb{T}}U(\Phi_{-})^{\ast}d\lambda(z)=0.
\]
Then
\begin{align}
\label{pro2}\int_{ \mathbb{T}} U(\mathbb{P}(\Phi))^{\ast}\,d\lambda(z) =
\int_{ \mathbb{T}}U \Phi^{\ast}\,d\lambda(z).
\end{align}
Now set $\Phi= UU^{\ast}U$. By the cubic matrix Szeg\H{o} equation
(\ref{cubic matrix}), we have
\begin{align*}
\frac{d\mathcal{Q}}{dt}  & = \int_{\mathbb{T}} \left( \partial_{t} U\right)
U^{\ast} d\lambda(z)+\int_{\mathbb{T}} U \left( \partial_{t}U^{\ast} \right)
d\lambda(z)\\
& = -i\int_{\mathbb{T}} \mathbb{P}(UU^{\ast}U) U^{\ast} d\lambda
(z)+i\int_{\mathbb{T}} U (\mathbb{P}(UU^{\ast}U))^{\ast} d\lambda(z)\\
& = -i\int_{\mathbb{T}} (UU^{\ast}U) U^{\ast} d\lambda(z)+i\int_{\mathbb{T}} U
(UU^{\ast}U)^{\ast} d\lambda(z) \quad\text{by (\ref{pro1}), (\ref{pro2})}\\
& = -i\int_{\mathbb{T}} UU^{\ast}U U^{\ast} d\lambda(z)+i\int_{\mathbb{T}}
UU^{\ast}U U^{\ast} d\lambda(z)\\
& = 0.
\end{align*}
We therefore conclude that $\frac{d}{dt}\mathcal{Q}(U(t))=0$, which means
$\mathcal{Q}(U(t))$ is a conservation law of (\ref{cubic matrix}).
\end{proof}

\begin{remark}
The conserved quantity $\mathcal{Q}(U(t))$ is a matrix-valued analogue of the
$L^{2}$-mass conservation law for the cubic scalar Szeg\H{o} equation
(\ref{cubic scalar}) (see \cite{Gerard10}). Indeed, when $n=1$, it reduces to
the standard mass:
\[
\mathcal{Q}(u(t))=\int_{\mathbb{T}}|u(t,z)|^{2}\,d\lambda(z)=\Vert
u(t)\Vert_{L^{2}}^{2}.
\]
Furthermore, in the general matrix-valued setting, taking the trace of
$\mathcal{Q}(U(t))$ directly implies the conservation of the $L_{n\times
n}^{2}$-norm of the solution:
\[
\operatorname{tr}(\mathcal{Q}(U(t)))=\int_{\mathbb{T}}\operatorname{tr}\left(
U(t,z)U(t,z)^{\ast}\right)  \,d\lambda(z)=\Vert U(t)\Vert_{L_{n\times n}^{2}%
}^{2}.
\]
Similarly, the quantity
\[
\widehat{\mathcal{Q}}(U(t))=\int_{\mathbb{T}}U(t,z)^{\ast}U(t,z)d\lambda(z)
\]
is also conserved along the flow of
(\ref{cubic matrix}).
\end{remark}

The following lemma provides an identity for self-adjoint matrix-valued
functions that will be used later.

\begin{lemma}
\label{Fourier-identity} Let $A\in L^{2}( \mathbb{T}, \mathbb{C}^{n\times n})$
be self-adjoint, that is, $A(z)^{\ast}=A(z)$ for almost every $z\in\mathbb{T}%
$. Then
\begin{align*}
2\sum_{k,j=1}^{n} \|\mathbb{P}(A\mathbb{E}_{kj})\|_{H^{2}_{n\times n}}^{2} = n
\left\Vert A \right\Vert _{L^{2}_{n\times n}}^{2} + n \left\Vert
\int_{\mathbb{T}} A d\lambda(z) \right\Vert _{\mathbb{C}^{n\times n}}^{2}.
\end{align*}

\end{lemma}

\begin{proof}
By the Fourier representation of the $H^{2}_{n\times n}$-norm,
\begin{align*}
\sum_{k,j} \|\mathbb{P}(A\mathbb{E}_{kj})\|_{H^{2}_{n\times n}}^{2}  & =
\sum_{k,j}\sum_{m=0}^{\infty} \left\Vert \widehat{A \mathbb{E}_{kj}}(m)
\right\Vert _{\mathbb{C}^{n\times n}}^{2}\\
& = \sum_{k,j}\sum_{m=0}^{\infty} \left\Vert \widehat{A}(m) \mathbb{E}%
_{kj}\right\Vert _{\mathbb{C}^{n\times n}}^{2}\\
& = \sum_{m=0}^{\infty} \operatorname{tr} \left(  \widehat{A}(m) \left(
\sum_{k,j} \mathbb{E}_{kj} (\mathbb{E}_{kj})^{\ast}\right)  \left( \widehat
{A}(m) \right) ^{\ast}\right) \\
& = \sum_{m=0}^{\infty} n \left\Vert \widehat{A}(m) \right\Vert _{\mathbb{C}%
^{n\times n}}^{2} .
\end{align*}
Here, the above equality follows from the identity
\[
\sum_{k,j} \mathbb{E}_{kj} (\mathbb{E}_{kj})^{\ast} = n I_{n},
\]
where $I_{n}$ denotes the $n\times n$ identity matrix.

Since $A$ is self-adjoint, its Fourier coefficients satisfy
\[
\widehat A(-m)=\widehat A(m)^{\ast},\qquad m\in\mathbb{Z}.
\]
Moreover,
\[
\|\widehat A(-m)\|_{ \mathbb{C}^{n\times n}} = \|\widehat A(m)\|_{
\mathbb{C}^{n\times n}},
\]
and hence
\begin{align*}
2 \sum_{m=0}^{\infty} \left\Vert \widehat{A}(m) \right\Vert _{\mathbb{C}%
^{n\times n}}^{2}  & = \sum_{m=-\infty}^{\infty} \left\Vert \widehat{A}(m)
\right\Vert _{\mathbb{C}^{n\times n}}^{2} + \left\Vert \widehat{A}(0)
\right\Vert _{\mathbb{C}^{n\times n}}^{2}\\
& = \left\Vert A \right\Vert _{L^{2}_{n\times n}}^{2} + \left\Vert
\int_{\mathbb{T}} A d\lambda(z) \right\Vert _{\mathbb{C}^{n\times n}}^{2}.
\end{align*}
Combining the above identities gives the desired result.
\end{proof}

Having established the above conservation laws and the auxiliary identity, we
now turn to the proof of part (iv) of Theorem \ref{main result}. Compared with
the scalar case in \cite{Gerard23}, the main new ingredient is the identity
(\ref{second deri}) below.

\textit{\textbf{Proof of Theorem \ref{main result}(iv):}} Assume first that
$UU^{\ast}\in L^{2}_{n\times n}$. We will show that the $L^{2}$-norm of
$UU^{\ast}$ is preserved by the flow. The first step is to recover this
quantity from the second derivative of the resolvent functional $J(x,U)$ at
$x=0_{+}$.

For every $U\in H^{2}_{n\times n}$ and $x\ge0$, define
\begin{align*}
J(x,U):=\sum_{k,j=1}^{n} \langle(I+xH^{(r)}_{U}H^{(l)}_{U})^{-1}%
\mathbb{E}_{kj},\mathbb{E}_{kj}\rangle_{H^{2}_{n\times n}}.
\end{align*}
For each fixed $U$, the function $x \mapsto J(x,U)$ is of class $C^{\infty}$
on $(0,\infty)$ with
\begin{align*}
&  \partial_{x}J(x,U)\\
&  =-\sum_{k,j} \langle(I+xH^{(r)}_{U}H^{(l)}_{U})^{-1}H^{(r)}_{U}H^{(l)}%
_{U}(I+xH^{(r)}_{U}H^{(l)}_{U})^{-1}\mathbb{E}_{kj},\mathbb{E}_{kj}\rangle\\
&  = -\sum_{k,j} \langle H^{(l)}_{U}(I+xH^{(r)}_{U}H^{(l)}_{U})^{-1}%
\mathbb{E}_{kj},H^{(l)}_{U}(I+xH^{(r)}_{U}H^{(l)}_{U})^{-1}\mathbb{E}%
_{kj}\rangle\\
&  = -\sum_{k,j}\Vert H^{(l)}_{U}(I+xH^{(r)}_{U}H^{(l)}_{U})^{-1}%
\mathbb{E}_{kj}\Vert^{2}.
\end{align*}
Furthermore, for $x>0$,
we have
\begin{align*}
&  \partial^{2} _{x}J(x,U) \\
&  = 2\sum_{k,j} \langle(H^{(r)}_{U}H^{(l)}_{U})^{2}(I+xH^{(r)}_{U}H^{(l)}%
_{U})^{-3}\mathbb{E}_{kj},\mathbb{E}_{kj}\rangle\\
&  = 2\sum_{k,j}\Vert H^{(r)}_{U}H^{(l)}_{U} (I+xH^{(r)}_{U}H^{(l)}%
_{U})^{-3/2}\mathbb{E}_{kj}\Vert^{2}\\
&  =2\sum_{k,j}\Vert H^{(r)}_{U}(I+xH^{(l)}_{U}H^{(r)}_{U})^{-3/2} H^{(l)}_{U}
\mathbb{E}_{kj} \Vert^{2}\\
& =2\sum_{k,j}\Vert H^{(r)}_{U}(I+xH^{(l)}_{U}H^{(r)}_{U})^{-3/2} \left(
(\mathbb{E}_{kj})^{\ast}U \right) \Vert^{2} .
\end{align*}

Observe that $(\mathbb{E}_{kj})^{\ast}U \in\operatorname{Dom} H^{(r)}_{U}$ for
all $k,j$ if and only if
\begin{align*}
\sup_{\substack{x > 0}} \sum_{k,j}\Vert H^{(r)}_{U}(I+xH^{(l)}_{U}H^{(r)}%
_{U})^{-3/2} \left( (\mathbb{E}_{kj})^{\ast}U \right) \Vert^{2}_{H^{2}%
_{n\times n}} < \infty,
\end{align*}
Equivalently, this condition holds if and only if the second derivative
$\partial^{2} _{x}J(x,U)$ remains bounded as $x \to0_{+}$, and in this case
the identity
\begin{align*}
\partial^{2} _{x}J(0_{+},U)=2\sum_{k,j}\Vert H^{(r)}_{U}\left( (\mathbb{E}%
_{kj})^{\ast}U \right) \Vert^{2}%
\end{align*}
holds.

Next step, we will show that the identity (\ref{second deri}) holds.

By the definition of $\operatorname{Dom} H^{(r)}_{U}$,
\[
(\mathbb{E}_{kj})^{\ast}U \in\operatorname{Dom} H^{(r)}_{U} \quad\text{if and
only if} \quad\mathbb{P}\left(  U \cdot\left(  (\mathbb{E}_{kj})^{\ast}U
\right)  ^{\ast}\right)  \in H^{2}_{n\times n}.
\]
Consequently, by Lemma \ref{Fourier-identity},
\begin{align*}
2\sum_{k,j }\left\Vert H^{(r)}_{U} \left(  (\mathbb{E}_{kj})^{\ast}U \right)
\right\Vert _{H^{2}_{n\times n}}^{2}  &  = 2\sum_{k,j}\left\Vert
\mathbb{P}\left(  U \left(  (\mathbb{E}_{kj})^{\ast}U \right) ^{\ast} \right)
\right\Vert _{H^{2}_{n\times n}}^{2}\\
& = 2\sum_{k,j} \|\mathbb{P}(UU^{\ast} \mathbb{E}_{kj})\|_{H^{2}_{n\times n}%
}^{2}\\
& = n \Vert U U^{\ast} \Vert_{L^{2}_{n\times n}}^{2} +n \left\Vert
\int_{\mathbb{T}} U U^{\ast} d\lambda(z) \right\Vert _{\mathbb{C}^{n \times
n}}^{2} .
\end{align*}
And in this case, $(\mathbb{E}_{kj})^{\ast}U \in\text{Dom} H^{(r)}_{U}$ for
all $k,j$ if and only if $UU^{\ast} \in L^{2}_{n\times n}$.

Thus we conclude that $UU^{\ast}\in L^{2}_{n\times n}$ if and only if the
second derivative $\partial_{x}^{2}J(x,U)$ remains bounded as $x\rightarrow
0_{+}$, and in this case
\begin{equation}
\label{second deri}\partial^{2}_{x}J(0_{+},U) = n \Vert U U^{\ast}
\Vert_{L^{2}_{n\times n}}^{2} +n \left\Vert \int_{\mathbb{T}} U U^{\ast}
d\lambda(z) \right\Vert _{\mathbb{C}^{n \times n}}^{2}.
\end{equation}

Furthermore, for initial data $U\in\mathcal{H}_{+}^{s}(\mathbb{T}%
,\mathbb{C}^{n\times n})$, $s>1/2$, the corresponding solution
\[
\mathcal{W}(t)U\in C(\mathbb{R}, \mathcal{H}_{+}^{s}(\mathbb{T},\mathbb{C}%
^{n\times n}))
\]
satisfies, by the Lax pair structure, that $J(x,\mathcal{W}(t)U)$ is a
conservation law of the cubic matrix Szeg\H{o} equation (\ref{cubic matrix});
see Lemma \ref{conservation law}. Hence,
\[
J(x,\mathcal{W}(t)U)=J(x,U)
\]
for all sufficiently small $x>0$. Moreover, Lemma \ref{conservation-Q} gives
the conservation law
\[
\int_{ \mathbb{T}} (\mathcal{W}(t)U)(\mathcal{W}(t)U)^{\ast}d\lambda(z) =
\int_{ \mathbb{T}}UU^{\ast}d\lambda(z).
\]

We next extend this conservation property to general $H^{2}_{n\times n}$ data.
If $U_{k}\rightarrow U$ in $H^{2}_{n\times n}$, then Theorem
\ref{strong con-H} gives
\[
H_{U_{k}}^{(r)}H_{U_{k}}^{(l)} \rightarrow H_{U}^{(r)}H_{U}^{(l)}
\]
in the strong resolvent sense. Consequently, for each fixed $x > 0$, the map
$U\mapsto J(x,U)$ from $H^{2}_{n\times n}$ to $\mathbb{C}$ is continuous.
Combining this continuity with Theorem \ref{main result}(ii), we obtain,
\[
J(x,\mathcal{W}(t)U)=J(x,U), \quad U\in H^{2}_{n\times n}.
\]

Now assume that $UU^{\ast}\in L^{2}_{n\times n}$. By identity
(\ref{second deri}), the above two conservation laws imply that
\[
(\mathcal{W}(t)U)(\mathcal{W}(t)U)^{\ast}\in L^{2}_{n\times n}
\]
and moreover,
\begin{align*}
\Vert(\mathcal{W}(t)U)(\mathcal{W}(t)U)^{\ast}\Vert_{L^{2}_{n\times n}}^{2}  &
=\partial^{2}_{x}J(0_{+},\mathcal{W}(t)U) -n \left\Vert \int_{\mathbb{T}}
\left( \mathcal{W}(t)U\right) \left( \mathcal{W}(t)U\right) ^{\ast}
d\lambda(z) \right\Vert _{\mathbb{C}^{n \times n}}^{2}\\
&  =\partial^{2}_{x}J(0_{+}, U) -n \left\Vert \int_{\mathbb{T}} UU^{\ast}
d\lambda(z) \right\Vert _{\mathbb{C}^{n \times n}}^{2}\\
&  =\Vert UU^{\ast}\Vert_{L^{2}_{n\times n}}^{2}.
\end{align*}
This completes the proof of Theorem \ref{main result}(iv).

\section{Proof of the basic properties of the right Hankel operator}

\label{Polynomials are dense} This section is devoted to proving Theorem
\ref{Domain dense}. The argument relies on the commutation relation for the
right Hankel operator $\Gamma_{U}^{(r)}$, which implies that its graph is
invariant under $S\oplus S ^{\ast}$. This allows us to apply the
description of the invariant subspaces of $S_{E} \oplus S_{F}^{\ast}$ on the
vector-valued Hardy space $H^{2}_{E} \oplus H^{2}_{F}$ established by Gu and
Luo \cite{GuLuo}.

\subsection{The commutation relation and the graph of
\texorpdfstring{$\Gamma^{(r)}_{U}$}{}}

\label{commu sec}

The following lemma establishes a commutation relation that reflects the
invariance of the underlying Hardy spaces. This relation is crucial for
proving Theorem \ref{Domain dense}.

\begin{lemma}
\label{commutation lemma} Let $U \in H^{2}_{n \times n}$, and let
$\Gamma^{(r)}_{U}$ be the right Hankel operator, defined by
(\ref{def of right Hankel}) on the domain (\ref{Dom of right Hankel}). Then
$\Gamma^{(r)}_{U}$ satisfies the following commutation relation
\begin{align}
\label{commutation relation}\Gamma^{(r)}_{U}S =S ^{\ast}\Gamma^{(r)}_{U}.
\end{align}

\end{lemma}

The above formula exactly meaning: for any $f \in\operatorname{Dom}
\Gamma^{(r)}_{U}$, we have $Sf \in\operatorname{Dom} \Gamma^{(r)}_{U}$ and
$\Gamma^{(r)}_{U}Sf=S^{\ast}\Gamma^{(r)}_{U}f$. This commutation relation can
be verified directly on the domain (for more details, see V. Peller \cite[Page
71]{Peller}).

The commutation relation also holds for the left Hankel operator.

\begin{lemma}
Let $U \in H^{2}_{n \times n}$, and let $\Gamma^{(l)}_{U}$ be the left Hankel
operator, defined by (\ref{def of left Hankel}) on the domain
(\ref{Dom of left Hankel}). Then $\Gamma^{(l)}_{U}$ satisfies the following
commutation relation
\[
\Gamma^{(l)}_{U}S =S ^{\ast}\Gamma^{(l)}_{U}.
\]

\end{lemma}


We also establish the following two commutation relations.

By Lemma \ref{commutation lemma}, we have, for any $f \in\operatorname{Dom}
H^{(r)}_{U}$,
\begin{align*}
S ^{\ast}H^{(r)}_{U}f  &  = S^{\ast}\,\Gamma^{(r)}_{U}f^{\#} \overset
{\text{by} (\ref{commutation relation})}{=} \Gamma^{(r)}_{U}(Sf^{\#})\\
&  =\mathbb{P}(U\overline{z}(f)^{\ast})\\
&  = \Gamma^{(r)}_{U}(S f)^{\#} = H^{(r)}_{U}(S f).
\end{align*}
Thus the right shifted Hankel operator $K^{(r)}_{U} := H^{(r)}_{S^{\ast}U}$
can be alternatively written as
\begin{align}
\label{K-right}K^{(r)}_{U} = S^{\ast}H^{(r)}_{U} = H^{(r)}_{U}S
\end{align}
on $\operatorname{Dom} H^{(r)}_{U}$, and the left shifted Hankel operator
$K^{(l)}_{U} := H^{(l)}_{S^{\ast}U}$ can be written as
\begin{align}
\label{K-left}K^{(l)}_{U} = S^{\ast}H^{(l)}_{U} = H^{(l)}_{U}S
\end{align}
on $\operatorname{Dom} H^{(l)}_{U}$. Moreover,
\[
\operatorname{Dom} K^{(r)}_{U} = \operatorname{Dom} H^{(r)}_{U},
\qquad\text{Dom}\,K^{(l)}_{U} = \operatorname{Dom} H^{(l)}_{U}.
\]

Next, our proof of Theorem \ref{Domain dense} is based on the following observation.

Clearly, using Lemma \ref{commutation lemma}, $\text{Dom} \,\Gamma^{(r)}_{U}$
is invariant for $S$. In addition, using Lemma \ref{commutation lemma},
\begin{align*}
(S \oplus S ^{\ast})(f,\Gamma^{(r)}_{U}(f))  &  =(S f,S ^{\ast}\Gamma
^{(r)}_{U}(f))\\
&  \overset{\text{by} (\ref{commutation relation})}{=}(S f,\Gamma^{(r)}_{U}S
(f)) \in\operatorname{Graph}\, (\Gamma^{(r)}_{U}),
\end{align*}
that is
\begin{align*}
(S\oplus S ^{\ast})\operatorname{Graph}\, (\Gamma^{(r)}_{U})\subset
\operatorname{Graph}\, (\Gamma^{(r)}_{U}).
\end{align*}
Then, we also see that $\operatorname{Graph}\, (\Gamma^{(r)}_{U})$ is
invariant under $S \oplus S ^{\ast}$. Therefore, we need to study the
subspaces of $H^{2}_{n \times n}\oplus H^{2}_{n \times n}$ that are invariant
under $S \oplus S ^{\ast}$.

\subsection{Invariant subspaces of \texorpdfstring{$S_E \oplus S_F^\ast$}{}}

Recently, Gu and Luo \cite{GuLuo} characterized the invariant
subspaces of $S_{E}\oplus S_{F}^{\ast}$ on the vector-valued Hardy space
$H^{2}_{E}\oplus H^{2}_{F}$, where $S_{E}$ and $S_{F}$ denote the shift
operators on $H_{E}^{2}$ and $H_{F}^{2}$, respectively. In order to state
their main results, we need to introduce the following notation.

Let
\begin{equation}
\Phi(z)=\left[
\begin{array}
[c]{cc}%
A(z) & B(z)\\
C(z) & D(z)
\end{array}
\right]  , \label{defphi}%
\end{equation}
where $A\in H_{B(E)}^{\infty}$, $B\in H_{B(F,E)}^{\infty}$, $C\in
L_{B(E,F)}^{\infty}$, $D\in L_{B(F)}^{\infty}$. Let $\Gamma_{E\oplus
F}^{\infty}$ denote the set of all symbols $\Phi$ above. Let $W_{\Phi}$ be a
bounded operator defined by
\[
W_{\Phi}=\left[
\begin{array}
[c]{cc}%
T_{A} & T_{B}\\
H_{C} & H_{D}%
\end{array}
\right]  ,
\]
where the corresponding definitions of Toeplitz and Hankel operators
can be found in Appendix \ref{notions change}.

Recall $\Theta\in H_{B\left(  F_{1}, F\right)  }^{\infty}$ is left inner if
$\Theta(z)^{\ast} \Theta(z)=I_{F_{1}}$ a.e. $z \in\mathbb{T}$ for some $F_{1}
\subset F, F_{1} \neq\{0\}$. In addition, we say $\Theta$ is right extremal if
$\Theta$ can not be factored as $\Theta_{1} \eta$ where $\Theta_{1}$ is left
inner and $\eta$ is a non-constant square inner function (i.e. $\eta\in
H_{B\left(  F_{1}\right)  }^{\infty}$ is not constant valued, and
$\eta(z)^{\ast} \eta(z)=\eta(z) \eta(z)^{\ast}=I_{F_{1}}$ a.e. $z
\in\mathbb{T}$ ). Let $K_{\Theta}=H_{F}^{2}\ominus\Theta H_{F_{1}}^{2}$ denote
the model space. $\Phi(z)$ is partial isometry-valued means that
$\Phi(z)^{\ast} \Phi(z)=I_{E_{0}}$ a.e. $z \in\mathbb{T}$ and $E_{0} \subset E
\oplus F$.

\begin{theorem}
\textup{\cite[Theorem 4.15]{GuLuo}} \label{GuLuo} A closed subspace $M$ of
$H_{E}^{2} \oplus H_{F}^{2}$ is an invariant subspace of $S_{E} \oplus
S_{F}^{\ast}$, if and only if
\[
M=\text{Span}\left\{  \operatorname{Ran}\left(  W_{\Phi}\right)  , K_{\Theta
}\right\}  ,
\]
where $\text{Span}\left\{  \operatorname{Ran}\left(  W_{\Phi}\right)
, K_{\Theta}\right\} $ denotes the closed linear span of $\operatorname{Ran}%
(W_{\Phi})$ and $K_{\Theta}$, $\Phi\in\Gamma_{E \oplus F}^{\infty}$ is zero
or partial isometry-valued, and $\Theta\in H_{B\left(  F_{1}, F\right)
}^{\infty}$ is zero or left inner, right extremal.
\end{theorem}

Thus, we obtain the description of the invariant subspaces of $S_{n
\times1}\oplus S_{n \times1}^{\ast}$ on the matrix-valued Hardy space
$H^{2}_{n \times1}\oplus H^{2}_{n \times1}$, which will be used as a key tool
in the proof of Theorem \ref{Domain dense}. Note that the notation in
Corollary \ref{GuLuo Corol} slightly differs from that in \cite{GuLuo}. These
changes are presented in Appendix \ref{notions change} for the interested
readers.

\begin{corollary}
\label{GuLuo Corol} A closed subspace $M$ of $H^{2}_{n \times1}\oplus H^{2}_{n
\times1}$ is an invariant subspace of $S_{n \times1}\oplus S_{n \times1}%
^{\ast}$, if and only if $M=\text{Span}\left\{  \operatorname{Ran}\left(
W_{\Phi}\right)  , K_{\Theta}\right\}  $, where $\Phi\in\Gamma_{\mathbb{C}^{n
\times1} \oplus\mathbb{C}^{n\times1}}^{\infty}$ is partial isometry-valued,
and $\Theta\in H_{B\left(  F_{1}, \mathbb{C}^{n \times1}\right)  }^{\infty}$
is zero or left inner, right extremal.
\end{corollary}

\subsection{Proof of Theorem \ref{Domain dense}}

\label{Theorem-Domain dense}

Before proving the theorem, we first reduce the matrix-valued setting to the
standard vector-valued one.

Let $N:=n^{2}.$ We identify
\[
H^{2}_{n\times n}\cong H^{2}_{N\times1}
\]
via the vectorization map
\[
\operatorname{Vec}: \mathbb{C}^{n\times n}\to\mathbb{C}^{N\times1},
\]
which stacks the columns of a matrix into a single column vector.

Next, define the symbol
\[
w=
\begin{pmatrix}
U & 0 & \cdots & 0\\
0 & U & \cdots & 0\\
\vdots & \vdots & \ddots & \vdots\\
0 & 0 & \cdots & U
\end{pmatrix}
\in H^{2}_{n^{2}\times n^{2}}.
\]
Then the right Hankel operator $\Gamma_{U}^{(r)}$ is unitarily equivalent to
the vector-valued Hankel operator $\Gamma_{w}$. Consequently, recall
(\ref{Dom of right Hankel}),
\[
\operatorname{Dom}\Gamma_{U}^{(r)} \cong\operatorname{Dom}\Gamma_{w},
\]
where
\[
\operatorname{Dom}\Gamma_{w} = \left\{  h\in H^{2}_{N\times1}: \mathbb{P}%
(wJh)\in H^{2}_{N\times1} \right\} .
\]

In this subsection, we will present the proof of Theorem \ref{Domain dense},
following the idea of \cite{Gerard23}. The proof combines the commutation
relations established above with the description of invariant
subspaces of $S_{N\times1}\oplus S_{N\times1}^{\ast}$ on $H^{2}_{N\times
1}\oplus H^{2}_{N\times1}$, provided by Corollary \ref{GuLuo Corol}, applied
with dimension $N$. Therefore, by the above unitary equivalence, it
suffices to establish Theorem \ref{Domain dense} for the vector-valued Hankel
operator $\Gamma_{w}$.



\textbf{\textit{Proof of Theorem \ref{Domain dense}(i):}} Suppose that
$\{f^{(k)}\}$ is a sequence in $\operatorname{Dom}\Gamma_{w}$ which is Cauchy
with respect to the graph norm. The sequence $\{f^{(k)}\}$ is Cauchy in the
usual norm of $H^{2}_{N \times1} $ and therefore
\[
\Vert f^{(k)}-f \Vert_{H^{2}_{N \times1}} \to0
\]
for some $f \in H^{2}_{N \times1}$. Moreover, there exists $h\in
H^{2}_{N\times1}$ such that
\[
\|\Gamma_{w}f^{(k)}-h\|_{H^{2}_{N\times1}}\to0 .
\]
It remains to show that $f \in\operatorname{Dom} \Gamma_{w}$ and $\Gamma
_{w}f=h$.

Denote by $\mathbb{P}_{M}$ the projection from $L^{1}_{N\times1}$ onto the
subspace of analytic polynomials of degree at most $M$. Clearly, $w Jf^{(k)},$
$wJf\in L^{1}_{N\times1}.$ Convergence $\Vert\Gamma_{w}f^{(k)}-h\Vert
_{H^{2}_{N \times1}}\rightarrow0$ implies that, for every fixed $M$,
\begin{equation}
\Vert\mathbb{P}_{M}(\Gamma_{w}f^{(k)}-h)\Vert_{H^{2}_{N \times1}}\rightarrow0.
\label{poly}%
\end{equation}
On the other hand, the convergence $\Vert f^{(k)}-f\Vert_{H^{2}_{N \times1}%
}\rightarrow0$ implies
\[
\Vert w Jf^{(k)}-wJf\Vert_{L^{1}_{N\times1}}\rightarrow0.
\]
Consequently, the Fourier coefficients of $w Jf^{(k)}$ converge to those of
$wJf$. Therefore, for every fixed $M$,
\[
\Vert\mathbb{P}_{M}(\Gamma_{w}f^{(k)})-\mathbb{P}_{M}(w Jf)\Vert_{H^{2}_{N
\times1}}\rightarrow0,\quad k\rightarrow\infty.
\]
Combining this with (\ref{poly}), we obtain
\begin{equation}
\label{poly1}\mathbb{P}_{M}(wJf)=\mathbb{P}_{M}(h)
\end{equation}
for every $M$. Hence all nonnegative Fourier coefficients of $wJf$ coincide
with those of $h$, which implies that
\[
\mathbb{P}(wJf)=h\in H^{2}_{N \times1}.
\]
Thus, $f\in\operatorname{Dom}\Gamma_{w}$, and (\ref{poly1}) can be rewritten
as
\[
\mathbb{P}_{M}(\Gamma_{w}f)=\mathbb{P}_{M}(h)
\]
for every $M$. Since $M$ is arbitrary, the Fourier coefficients of $\Gamma
_{w}f$ and $h$ coincide for all nonnegative indices, and hence $\Gamma_{w}f=h
$.  The proof of Theorem \ref{Domain dense}(i) is complete.

\textbf{\textit{Proof of Theorem \ref{Domain dense}(ii):}} Let $Y$ denote the
graph of $\Gamma_{w}:$%
\[
Y=\operatorname{Graph}\,(\Gamma_{w})=\left\{  (f,\Gamma_{w}f):f\in
\operatorname{Dom}\Gamma_{w}\right\}  \subset H^{2}_{N \times1}\oplus H^{2}_{N
\times1}.
\]
Let $Y_{0}$ be the nonclosed subspace of $Y:$%
\[
Y_{0}=\left\{  (p,\Gamma_{w}p):p\in\mathcal{P}_{N\times1}\right\}  .
\]

We have proved that $Y$ is closed in $H^{2}_{N \times1} \oplus H^{2}_{N
\times1}$ in the above Theorem \ref{Domain dense}(i), and now we prove that
the closure of $Y_{0}$ coincides with $Y$. We will do this by proving that
$Y_{0}^{\perp}\subset Y^{\perp}$, where the orthogonal complements are taken
in $H^{2}_{N \times1}\oplus H^{2}_{N \times1}$.

Since $Y_{0}$ is invariant for $S_{N \times1}\oplus S_{N \times1}^{\ast}$, its
orthogonal complement $Y_{0}^{\perp}$ is invariant for $S_{N \times1}^{\ast
}\oplus S_{N \times1}$. Define $\mathbf{V}$ to be the unitary operator by
\begin{align*}
\mathbf{V}  &  \mathbf{:}H^{2}_{N \times1}\oplus H^{2}_{N \times1}\rightarrow
H^{2}_{N \times1}\oplus H^{2}_{N \times1},\\
\mathbf{V} & (f\oplus h) =h\oplus f,\quad f\in H^{2}_{N \times1},\quad h\in
H^{2}_{N \times1}.
\end{align*}
Then
\[
S_{N \times1}^{\ast}\oplus S_{N \times1}=\mathbf{V}^{\ast}(S_{N \times1}\oplus
S_{N \times1}^{\ast})\mathbf{V}.
\]
Thus, $\mathbf{V}Y_{0}^{\perp}$ is invariant for $S_{N \times1}\oplus S_{N
\times1}^{\ast}$, if and only if $Y_{0}^{\perp}$ is invariant for $S_{N
\times1}^{\ast}\oplus S_{N \times1}$.

Applying Corollary \ref{GuLuo Corol} with the multiplicity parameter
chosen as $N=n^{2}$, we obtain a description of the invariant subspace
$\mathbf{V}Y_{0}^{\perp}$. Let $\Phi,\Theta$ be as in Corollary
\ref{GuLuo Corol}, write%
\[
\mathbf{V}Y_{0}^{\perp}=\text{Span}\left\{  \operatorname{Ran}\left(  W_{\Phi
}\right)  ,K_{\Theta}\right\}  .
\]
Set%
\begin{equation}%
\begin{bmatrix}
\alpha\\
\beta
\end{bmatrix}
:=W_{\Phi}%
\begin{bmatrix}
h_{1}\\
h_{2}%
\end{bmatrix}
=%
\begin{bmatrix}
Ah_{1}+Bh_{2}\\
\Gamma_{\overline{z}C}h_{1}+\Gamma_{\overline{z}D}h_{2}%
\end{bmatrix}
\label{defwv}%
\end{equation}
where $h_{1} ,h_{2}\in H^{2}_{N \times1}.$ Then%
\begin{equation}
\mathbf{V}Y_{0}^{\perp}=\left\{
\begin{bmatrix}
\alpha\\
\beta+k_{\Theta}%
\end{bmatrix}
:h_{1} ,h_{2}\in H^{2}_{N \times1},k_{\Theta}\in K_{\Theta}\right\}  ^{\_}
\label{vperp}%
\end{equation}

Note that $\mathbf{V}Y_{0}^{\perp}=(\mathbf{V}Y_{0})^{\perp}$ and
\[
\mathbf{V}Y_{0}=\left\{  (\Gamma_{w}p,p):p\in\mathcal{P}_{N \times1}\right\}
.
\]
The column vector of the form $\left[ \alpha,\ \beta+k_{\Theta}\right]
^{\mathsf{T}}$ belongs to $\mathbf{V}Y_{0}^{\perp}$ if and only if all
polynomials $p\in\mathcal{P}_{N \times1}$,
\[
\langle\alpha,\Gamma_{w}p\rangle_{H^{2}_{N \times1}}+\langle\beta+k_{\Theta
},p\rangle_{H^{2}_{N \times1}}=0.
\]
By (\ref{defwv}), equivalently
\begin{equation}
\langle Ah_{1}+Bh_{2},\mathbb{P}\left(  w Jp\right)  \rangle_{H^{2}_{N
\times1}}+ \langle\Gamma_{\overline{z}C}h_{1}+\Gamma_{\overline{z}D}%
h_{2}+k_{\Theta},p\rangle_{H^{2}_{N \times1}} =0. \label{orthogonality1}%
\end{equation}
Rewrite the first term in (\ref{orthogonality1}) as
\begin{align}
&  \langle Ah_{1}+Bh_{2},\mathbb{P}\left(  w Jp\right)  \rangle_{H^{2}_{N
\times1}}\label{ptof}\\
=  &  \langle Ah_{1}+Bh_{2},w Jp\rangle_{L^{2}_{N \times1}}\nonumber\\
=  &  \int_{\mathbb{T}}\mathrm{tr}\left( (Ah_{1}+Bh_{2})(z)\, \left(
w(z)\,p(\overline{z})\right) ^{\ast} \right) \,d\lambda(z)\quad\nonumber\\
=  &  \int_{\mathbb{T}}\mathrm{tr}\left( ((Ah_{1}+Bh_{2})p^{\#}%
)(z)\,w(z)^{\ast}\right) \,d\lambda(z)\nonumber\\
=  &  \langle(Ah_{1}+Bh_{2})p^{\#},w\rangle_{H^{2}_{N\times N}} .\nonumber
\end{align}
Thus, (\ref{orthogonality1}) can be rewritten as:
\begin{equation}
\langle(Ah_{1}+Bh_{2})p^{\#},w\rangle_{H^{2}_{N\times N} }+\langle
\Gamma_{\overline{z}C}h_{1}+\Gamma_{\overline{z}D}h_{2}+k_{\Theta}%
,p\rangle_{H^{2}_{N \times1}}=0 \label{orthogonality2}%
\end{equation}
for all polynomials $p\in\mathcal{P}_{N \times1}.$

To prove $Y_{0}^{\perp}\subset Y^{\perp}$, we need to show $Y\perp
Y_{0}^{\perp}$. Equivalently, we will show $\mathbf{V}Y\perp\mathbf{V}%
Y_{0}^{\perp}$. Recall that
\[
\mathbf{V}Y=\left\{  (\Gamma_{w}f,f):f\in\operatorname{Dom}(\Gamma
_{w})\right\}  .
\]
Take any $f\in\text{Dom}\,\Gamma_{w}$, we will check $(\Gamma_{w}%
f,f)\perp\mathbf{V}Y_{0}^{\perp}$. Since
\[
\left\{
\begin{bmatrix}
\alpha\\
\beta+k_{\Theta}%
\end{bmatrix}
:h_{1},h_{2}\in\mathcal{P}_{N\times1},k_{\Theta}\in K_{\Theta}\right\}
\]
is a dense subset of $\mathbf{V}Y_{0}^{\perp}$ because of $A,B\in H_{N\times
N}^{\infty}$, $C,D\in L_{N\times N}^{\infty}$ by Corollary \ref{GuLuo Corol}.
Thus, it suffices to check that
\[
\langle\alpha,\Gamma_{w}f\rangle_{H_{N\times1}^{2}}+\langle\beta+k_{\Theta
},f\rangle_{H_{N\times1}^{2}}=0,
\]
for all polynomials $h_{1},h_{2}\in\mathcal{P}_{N\times1}$, $k_{\Theta}\in
K_{\Theta}.$ Equivalently,
\begin{equation}
\langle Ah_{1}+Bh_{2},\mathbb{P}\left(  wJf\right)  \rangle_{H_{N\times1}^{2}%
}+\langle\Gamma_{\overline{z}C}h_{1}+\Gamma_{\overline{z}D}h_{2}+k_{\Theta
},f\rangle_{H_{N\times1}^{2}}=0.\label{check}%
\end{equation}
By the same computation as in (\ref{ptof}) with $p$ replaced by $f,$
(\ref{check}) can be rewritten as
\begin{equation}
\langle(Ah_{1}+Bh_{2})f^{\#},w\rangle_{H_{N\times N}^{2}}+\langle
\Gamma_{\overline{z}C}h_{1}+\Gamma_{\overline{z}D}h_{2}+k_{\Theta}%
,f\rangle_{H_{N\times1}^{2}}=0\label{check1}%
\end{equation}
for all $h_{1},h_{2}\in\mathcal{P}_{N\times1}$, $k_{\Theta}\in K_{\Theta}.$

Now let $\{p_{l}\}$ be a sequence of polynomials such that $\Vert p_{l}%
-f\Vert\rightarrow0$ as $l\rightarrow\infty$. For each $p=p_{l}$, each $h_{1},
h_{2}\in\mathcal{P}_{N\times1}$, and $k_{\Theta}\in K_{\Theta},$ orthogonality
condition (\ref{orthogonality2})
\[
\langle(Ah_{1}+Bh_{2})p_{l}^{\#},w\rangle_{H^{2}_{N\times N}}+\langle
\Gamma_{\overline{z}C} h_{1}+\Gamma_{\overline{z}D} h_{2}+k_{\Theta}%
,p_{l}\rangle_{H^{2}_{N \times1}}=0
\]
holds true. We pass to the limit $l\rightarrow\infty$ in the above equation.
Since $Ah_{1}+Bh_{2}\in H^{\infty}_{N \times1},$
\[
\Vert(Ah_{1}+Bh_{2})p_{l}^{\#}-(Ah_{1}+Bh_{2})f^{\#}\Vert\rightarrow0.
\]
This limit yields (\ref{check1}) which proves $\mathbf{V}Y\perp\mathbf{V}%
Y_{0}^{\perp}$. The proof of Theorem \ref{Domain dense}(ii) is complete.

\textbf{ \textit{Proof of Theorem \ref{Domain dense}(iii):}} Suppose that
$f\in\operatorname{Dom} \Gamma_{w} ^{\ast}$, then there exists $h\in
H^{2}_{N \times1}$ such that
\begin{equation}
\langle\Gamma_{w}h_{1},f\rangle_{H^{2}_{N \times1}}=\langle h_{1}%
,h\rangle_{H^{2}_{N \times1}},\quad\forall h_{1}\in\operatorname{Dom}%
\Gamma_{w} \label{inner product}%
\end{equation}
In this case, $h= \Gamma_{w} ^{\ast}f$. In particular, the above relation
holds for any polynomial $h_{1}\in\mathcal{P}_{N \times1}$. For polynomial
$h_{1}$, we have
\begin{align}
\langle\Gamma_{w}h_{1},f\rangle_{H^{2}_{N \times1}}  & = \langle\mathbb{P}(w
Jh_{1}),f\rangle\label{inner2}\\
& = \langle w Jh_{1},f\rangle_{L^{2}_{N \times1}}\nonumber\\
& = \int_{ \mathbb{T}} \operatorname{tr}\left(  w(z)h_{1}(\overline
z)f(z)^{\ast} \right) \,d\lambda(z)\nonumber\\
& = \int_{ \mathbb{T}} \operatorname{tr}\left(  h_{1}(\overline z)f(z)^{\ast
}w(z) \right) \,d\lambda(z)\nonumber\\
& = \int_{ \mathbb{T}} \operatorname{tr}\left(  h_{1}(z)f(\overline z)^{\ast
}w(\overline z) \right) \,d\lambda(z)\nonumber\\
& = \langle h_{1},w^{\#}Jf\rangle_{L^{2}_{N \times1}} .\nonumber
\end{align}
where $h_{\ast}=w ^{\#}Jf\in L^{1}_{N \times1}$. Since this identity agrees
with (\ref{inner product}) for all polynomials $h_{1}\in\mathcal{P}_{N
\times1}$, it follows that $f\in\operatorname{Dom}\Gamma_{w^{\#}} $, and
\[
\Gamma_{w^{\#}}f=\mathbb{P}h_{\ast}=h= \Gamma_{w} ^{\ast}f.
\]
This proves $\Gamma_{w} ^{\ast}\subset\Gamma_{w^{\#}}$.

Conversely, let $f\in\operatorname{Dom}\Gamma_{w^{\#}}$, then $h:=\mathbb{P}%
(w^{\#}Jf)\in H^{2}_{N \times1}$. For every polynomial $h_{1}\in\mathcal{P}_{N
\times1}$, by the same computation as in (\ref{inner2}), we have
\[
\langle\Gamma_{w}h_{1},f\rangle_{H^{2}_{N \times1}}=\langle h_{1}%
,\mathbb{P}\left(  w^{\#}Jf\right)  \rangle_{H^{2}_{N \times1}}.
\]
By Theorem \ref{Domain dense}-(ii), polynomials are dense in
$\operatorname{Dom}\Gamma_{w}$ with respect to the graph norm, this identity
extends to all $h_{1}\in\operatorname{Dom}\Gamma_{w}$. Therefore, $f\in
\Gamma_{w} ^{\ast}$ and $\mathbb{P}\left(  w^{\#}Jf\right)  =h= \Gamma_{w}
^{\ast}f$ for all $f\in\operatorname{Dom}\Gamma_{w^{\#}}$. This proves
$\Gamma_{w^{\#}}\subset\Gamma_{w} ^{\ast}.$ In summary, $\Gamma_{w} ^{\ast
}=\Gamma_{w^{\#}}$ and the proof of Theorem \ref{Domain dense}(iii) is complete.

\section{Wold decomposition and spectral properties of
\texorpdfstring{$\Omega^{(r)}$ and $\Omega^{(l)}$}{}}

\label{proof of complicated}

In this section, we prove Theorem \ref{complicated}, which shows that $U$
belongs to the pure part generated by the wandering subspace of the isometry
$\Omega^{(l)}$. The proof is organized into several steps.

In Subsection \ref{start}, we introduce the auxiliary anti-linear operators,
and state Theorem \ref{Wold reduce}, which shows that $\mathcal{H}^{(l)}$
maps the pure part of $\Omega^{(r)}$  into the pure part of $\Omega^{(l)}$ .
This property provides the main structural ingredient for the proof of Theorem
\ref{complicated}.

In Subsection \ref{ac section}, we establish the absence of absolutely
continuous spectral components in the unitary parts of $\Omega^{(r)}$ and
$\Omega^{(l)}$. 

Finally, in Subsection \ref{Lemma sub}, we prove Theorem
\ref{Wold reduce}. The argument relies on the result from Subsection
\ref{ac section}, stated in Lemma \ref{no a.c.part}, as well as on the
symmetry and commutation relations established in Lemma \ref{symmetryall}.

\subsection{The pure part of \texorpdfstring{$\Omega^{(l)}$}{} and Hankel
operator invariance}

\label{start}

Here we introduce two auxiliary anti-linear operators whose properties mirror
those of $H_{U}^{(r)}$ and $H_{U}^{(l)}$. We set
\begin{align*}
\mathcal{H}^{(r)}  &  =e^{itH_{U}^{(r)}H_{U}^{(l)}}H_{U}^{(r)},\quad
\operatorname{Dom}\mathcal{H}^{(r)}=\operatorname{Dom}H_{U}^{(r)}.\\
\mathcal{H}^{(l)}  &  =e^{itH_{U}^{(l)}H_{U}^{(r)}}H_{U}^{(l)},\quad
\operatorname{Dom}\mathcal{H}^{(l)}=\operatorname{Dom}H_{U}^{(l)}.
\end{align*}

Since $H_{U}^{(r)}$ and $H_{U}^{(l)}$ are closed operators, $\mathcal{H}%
^{(r)}$ and $\mathcal{H}^{(l)}$ are also closed operators.

The following theorem shows that $\mathcal{H}^{(l)}$  maps the pure part of
the Wold decomposition of $\Omega^{(r)}$  into the pure part of the Wold
decomposition of $\Omega^{(l)}$, and will be used in the proof of Theorem
\ref{complicated}, which is postponed to Section \ref{Lemma sub}.

\begin{theorem}
[Proof deferred to Section \ref{Lemma sub}]\label{Wold reduce}

The operator $\mathcal{H}^{(l)}$ maps the pure part of $\Omega^{(r)}$ into the
pure part of $\Omega^{(l)}$ on its domain. Specifically, let $P^{(r)}%
_{\mathrm{pure}}$ be the orthogonal projection in $X:=H^{2}_{n\times n}$ onto
the subspace $X^{(r)}_{\mathrm{pure}}$ in the Wold decomposition
(\ref{Wold decomposition u part}) of $\Omega^{(r)}$, and let $P^{(l)}%
_{\mathrm{pure}}$ be the orthogonal projection in $X$ onto the subspace
$X^{(l)}_{\mathrm{pure}}$ in the Wold decomposition of $\Omega^{(l)}$. Then we
have
\begin{align}
\label{eqdom} & P^{(r)}_{\mathrm{pure}}\operatorname{Dom}\mathcal{H}^{(l)}
\subset\operatorname{Dom}\mathcal{H}^{(l)} \text{ and }\\
&  \mathcal{H}^{(l)}P^{(r)}_{\mathrm{pure}}(\operatorname{Dom}\mathcal{H}%
^{(l)})\subset X^{(l)}_{\mathrm{pure}}.\label{eqpure}%
\end{align}

\end{theorem}

\textbf{\textit{Proof of Theorem \ref{complicated}:}} Throughout this
argument, we assume that
\[
U\in H^{2}_{n\times n}.
\]

For
\[
Q^{(r)} =\{q^{(r)}_{kj}:1\leq k, j\leq n\},\text{ where}\quad q^{(r)}%
_{kj}=e^{itH_{U}^{(r)}H_{U}^{(l)}}\mathbb{E}_{kj}.
\]
Since $H_{U}^{(r)}H_{U}^{(l)}$ is self-adjoint, $e^{itH_{U}^{(r)}H_{U}^{(l)}}$
is unitary on $H^{2}_{n\times n}$. Hence,
\[
q^{(r)}_{kj} = e^{itH_{U}^{(r)}H_{U}^{(l)}}\mathbb{E}_{kj} \in H^{2}_{n\times
n}.
\]

We also have
\begin{align*}
\mathcal{H}^{(l)} q^{(r)}_{kj}  &  =(H^{(l)}_{U}e^{-i t H^{(r)}_{U}H^{(l)}%
_{U}}) (e^{it H^{(r)}_{U}H^{(l)}_{U}}\mathbb{E}_{kj})\\
&  =H^{(l)}_{U}\mathbb{E}_{kj}\\
&  =\mathbb{P}((\mathbb{E}_{kj})^{\ast}U) =(\mathbb{E}_{kj})^{\ast}U \in
H^{2}_{n\times n}%
\end{align*}
for all $k,j$. Hence $q^{(r)}_{kj} \in\operatorname{Dom} \mathcal{H}^{(l)}$.
Recall that (\ref{rpure}),
\[
Q^{(r)}\subset X^{(r)}_{\mathrm{pure}} ,
\]
we obtain
\[
q^{(r)}_{kj}\in X^{(r)}_{\mathrm{pure}}\cap\operatorname{Dom}\mathcal{H}^{(l)}
.
\]
Together with (\ref{eqpure}) in Theorem \ref{Wold reduce}, we have for all
$k,j$,
\[
(\mathbb{E}_{kj})^{\ast}U =\mathcal{H}^{(l)}q^{(r)}_{kj} \in X^{(l)}%
_{\mathrm{pure}}.
\]

We now choose $j=k$. Then, for each $k=1,\dots,n$ , we have
\[
(\mathbb{E}_{kk})^{\ast}U = \mathcal{H}^{(l)}q_{kk}^{(r)} \in X^{(l)}%
_{\mathrm{pure}} .
\]
Moreover,
\[
U=\sum_{k=1}^{n}(\mathbb{E}_{kk})^{\ast}U .
\]
Since $X^{(l)}_{\mathrm{pure}} =\overline{\left\{ \left( \Omega^{(l)}\right)
^{m}Q^{(l)}:m\ge0\right\}  }$ is a closed linear subspace, it follows that
\[
U\in\overline{\left\{ \left( \Omega^{(l)}\right) ^{m}Q^{(l)}:m\ge0\right\}
}.
\]

Therefore, Theorem \ref{complicated} has been proved.

\subsection{Absence of absolutely continuous spectrum}

\label{ac section}

In this subsection, we show that the unitary parts of $\Omega^{(r)}$
and $\Omega^{(l)}$ have no absolutely continuous spectrum in Lemma
\ref{no a.c.part}. We begin by establishing a trace class perturbation framework.

\begin{lemma}
\label{trace class}Fix $U \in H^{2}_{n \times n}$ and $t \in\mathbb{R}$.
Recall
\begin{align*}
\Omega^{(r)} = e^{it H^{(r)}_{U}H^{(l)}_{U}} S e^{-it K^{(r)}_{U}K^{(l)}_{U}%
}:H^{2}_{n\times n}\to H^{2}_{n\times n},
\end{align*}
and define
\[
\Omega_{0}=e^{i t H^{(r)}_{U}H^{(l)}_{U}} S e^{-i t H^{(r)}_{U}H^{(l)}_{U}}.
\]
Then $\Omega^{(r)}-\Omega_{0}$ is trace class.
\end{lemma}

\begin{proof}
By (\ref{K-right}) and (\ref{K-left}), we have%
\begin{align}
\label{difference rank}H_{U}^{(r)}H_{U}^{(l)}-K_{U}^{(r)}K_{U}^{(l)}  &
=H_{U}^{(r)}H_{U}^{(l)}-H_{U}^{(r)}S S^{\ast}H_{U}^{(l)}\\
&  =H_{U}^{(r)}\left(  I -SS ^{\ast}\right)  H_{U}^{(l)}\nonumber\\
&  =H_{U}^{(r)} \Delta H_{U}^{(l)} .\nonumber
\end{align}

It follows from (\ref{difference rank}) that $H^{(r)}_{U}H^{(l)}%
_{U}-K^{(r)}_{U}K^{(l)}_{U}$ is a finite-rank operator. Since $H^{(r)}%
_{U}H^{(l)}_{U}$ and $K^{(r)}_{U}K^{(l)}_{U}$ are self-adjoint operators, by
the Duhamel formula \cite{Rodn04}, we have
\begin{align*}
&  e^{-itH^{(r)}_{U}H^{(l)}_{U}}-e^{-itK^{(r)}_{U}K^{(l)}_{U}}\\
=  &  -i\int_{0}^{t}e^{-i(t-s)K^{(r)}_{U}K^{(l)}_{U}} (H^{(r)}_{U}H^{(l)}%
_{U}-K^{(r)}_{U}K^{(l)}_{U})e^{-isH^{(r)}_{U}H^{(l)}_{U}}ds\\
=  &  - i\int_{0}^{t}e^{-i(t-s)K^{(r)}_{U}K^{(l)}_{U}} (H_{U}^{(r)} \Delta
H_{U}^{(l)} )e^{-isH^{(r)}_{U}H^{(l)}_{U}}ds,
\end{align*}
Multiply both sides of the above equation by $e^{itK^{(r)}_{U}K^{(l)}_{U}}$,
we have
\begin{align*}
&  e^{itK^{(r)}_{U}K^{(l)}_{U}}e^{-itH^{(r)}_{U}H^{(l)}_{U}}-I\\
=  &  -i\int_{0}^{t}e^{isK^{(r)}_{U}K^{(l)}_{U}} (H_{U}^{(r)} \Delta
H_{U}^{(l)} )e^{-isH^{(r)}_{U}H^{(l)}_{U}}ds
\end{align*}
is a trace class operator. It follows that the operator
\begin{align*}
\Omega^{(r)} -\Omega_{0}=e^{i t H^{(r)}_{U}H^{(l)}_{U}} S e^{-i t K^{(r)}%
_{U}K^{(l)}_{U}}\left(  I-e^{i t K^{(r)}_{U}K^{(l)}_{U}} e^{-i t H^{(r)}%
_{U}H^{(l)}_{U}}\right)
\end{align*}
is trace class.
\end{proof}

We recall the following Kato-Rosenblum's theorem, which will be used later.

\begin{theorem}[Kato-Rosenblum]
\label{Kato} \textup{\cite[Theorem X4.3]{Kato66}}  Let $A$ and
$B$ be self-adjoint operators in a Hilbert space $X$ such that the difference
$A-B$ is trace class. Then the absolutely continuous parts of $A$ and $B$ are
unitarily equivalent.
\end{theorem}

Before proving Lemma \ref{no a.c.part}, we firstly write $A^{\mathrm{ac}}
\simeq B^{\mathrm{ac}}$ to denote that the absolutely continuous parts of $A$
and $B$ are unitarily equivalent.

\begin{lemma}
\label{no a.c.part}Let
\begin{align*}
\Omega^{(r)}  &  = e^{it H^{(r)}_{U}H^{(l)}_{U}} S e^{-it K^{(r)}_{U}%
K^{(l)}_{U}}:H^{2}_{n\times n}\to H^{2}_{n\times n}\\
\Omega^{(l)}  &  = e^{it H^{(l)}_{U}H^{(r)}_{U}} S e^{-it K^{(l)}_{U}%
K^{(r)}_{U}}:H^{2}_{n \times n}\to H^{2}_{n \times n}.
\end{align*}
Then the unitary parts of both $\Omega^{(r)}$ and $\Omega^{(l)}$ have no
absolutely continuous components.
\end{lemma}

\begin{proof}
Writing $\mathrm{Re}\, \Omega^{(r)} =\left(  \Omega^{(r)} +\left( \Omega^{(r)}
\right) ^{\ast}\right)  / 2$, Lemma \ref{trace class} implies that the
difference
\begin{align*}
\mathrm{Re}\,\Omega^{(r)} -\mathrm{Re}\,\Omega_{0}%
\end{align*}
is a trace class operator. By Theorem \ref{Kato}, we have
\begin{align*}
(\mathrm{Re}\,\Omega^{(r)} )^{\mathrm{ac}} \simeq(\mathrm{Re}\,S)^{\mathrm{ac}%
} .
\end{align*}

Since
\begin{align*}
\mathrm{Re}\,S=1/2(S+S^{\ast})= 1/2(T_{z}+T_{\bar{z}})=T_{\frac{z+\bar{z}}{2}%
},
\end{align*}
thus $\sigma( \mathrm{Re}\, S)=[-1,1] $, and that $\mathrm{Re}\, S $ has
purely absolutely continuous spectrum \cite{Rose,Sobolev}, where
$\sigma(\mathrm{Re}\,S)$ denotes the spectrum of $\mathrm{Re}\,S$.

Let $N=n^{2}$. We identify $H^{2}_{n\times n}\cong H^{2}_{N\times1}$, via the
vectorization map
\[
\text{Vec}:\mathbb{C}^{n\times n} \to\mathbb{C}^{N\times1},
\]
which stacks the columns of a matrix into a single column vector. More
precisely, for $F = (f_{kj})_{1 \le k ,j \le n}\in H^{2}_{n\times n}$ we
define
\[
\text{Vec}(F) =
\begin{pmatrix}
f_{11}, \cdots, f_{n1}, f_{12}, \cdots, f_{n2}, \cdots, f_{1n}, \cdots, f_{nn}%
\end{pmatrix}
^{\mathsf{T}} \in H^{2}_{N \times1}.
\]
Under this identification, the shift operator
\begin{align*}
S_{n\times n}F(z)=zF(z),\, F\in H^{2}_{n\times n},
\end{align*}
is unitarily equivalent to the standard vector-valued shift of multiplicity
$N$:%
\[
S^{(N)}:H^{2}_{N\times1}\to H^{2}_{N\times1}, \quad(S^{(N)}\text{Vec}(F))(z) =
z\text{Vec}(F)(z).
\]

The multiplicity of a shift operator is the dimension of its wandering
subspace, i.e.
\[
\dim\ker\left(  (S^{(N)})^{\ast} \right) .
\]
For the scalar shift on $H^{2}(\mathbb{T}) $, this dimension is $1 $, since
the wandering subspace is spanned by the constant function $1 $. For
$H^{2}_{N\times1} $, we have
\[
\ker\left( (S^{(N)})^{\ast}\right)  = \{F\in H^{2}_{N\times1}:F(z)\equiv
c,\ c\in\mathbb{C}^{N\times1}\},
\]
thus the multiplicity of $S^{(N)}$ is $N$.

Therefore, from Toeplitz spectral theory, $\mathrm{Re}\,S_{n\times n}$
has purely absolutely continuous spectrum, spectrum $[-1,1]$, and spectral
multiplicity $n^{2}$.

Thus, we have
\begin{align}
\label{ac-part}(\mathrm{Re}\,\Omega^{(r)})^{\mathrm{ac}} \simeq\mathrm{Re}\,S.
\end{align}

In addition, by the Wold decomposition (\ref{Wold decomposition u part}), we
obtain
\[
(\mathrm{Re}\,\Omega^{(r)})^{\mathrm{ac}} \simeq(\mathrm{Re}\,\Omega
^{(r)}_{\mathrm{u}})^{\mathrm{ac}} \oplus(\mathrm{Re}\,\Omega^{(r)}%
_{\mathrm{pure}})^{\mathrm{ac}},
\]
where $\Omega^{(r)}_{\mathrm{pure}} $ is unitarily equivalent to $S$. Since
the real part of unitarily equivalent operators is also unitarily equivalent,
we have $(\mathrm{Re}\,\Omega^{(r)}_{\mathrm{pure}})^{\mathrm{ac}}
\simeq\mathrm{Re}\,S$. Therefore,
\[
(\mathrm{Re}\,\Omega^{(r)})^{\mathrm{ac}} \simeq(\mathrm{Re}\,\Omega
^{(r)}_{\mathrm{u}})^{\mathrm{ac}} \oplus\mathrm{Re}\,S.
\]
Comparing with (\ref{ac-part}), we have
\begin{align}
\label{ReS2-decomp}\mathrm{Re}\,S \simeq(\mathrm{Re}\,\Omega^{(r)}%
_{\mathrm{u}})^{\mathrm{ac}} \oplus\mathrm{Re}\,S.
\end{align}

Let $m_{(\cdot)}$ denote the multiplicity function of the absolutely
continuous spectrum. Then (\ref{ReS2-decomp}) implies
\[
m_{(\mathrm{Re}\,S)}(x) = m_{(\mathrm{Re}\,\Omega_{\mathrm{u}})^{\mathrm{ac}}%
}(x) + m_{(\mathrm{Re}\,S)}(x), \quad\text{a.e. } x \in\mathbb{R}.
\]
Since multiplicity functions are non-negative, it follows that
\[
(\mathrm{Re}\,\Omega^{(r)}_{\mathrm{u}})^{\mathrm{ac}} = 0.
\]

Let $E_{\Omega^{(r)}_{\mathrm{u}}}(\cdot)$ and $E_{\mathrm{Re}\,\Omega
^{(r)}_{\mathrm{u}}}(\cdot)$ denote the spectral measures of the unitary
operator $\Omega^{(r)} _{\mathrm{u}}$ and the self-adjoint operator
$\mathrm{Re}\,\Omega^{(r)}_{\mathrm{u}}$, respectively. By the functional
calculus for the unitary operator $\Omega^{(r)}_{\mathrm{u}}$, these
are related by
\[
E_{\mathrm{Re}\,\Omega^{(r)}_{\mathrm{u}}}(B) = E_{\Omega^{(r)}_{\mathrm{u}}%
}\left( \{ z \in\mathbb{T} : \mathrm{Re}\,z \in B \}\right) , \quad B
\subset\mathbb{R} \text{ Borel}.
\]
Hence the absolutely continuous part of the spectral measure of
$\operatorname{Re}\Omega^{(r)}_{\mathrm{u}}$ is the push-forward of that of
$\Omega^{(r)}_{\mathrm{u}}$. Since $(\mathrm{Re}\,\Omega^{(r)}_{\mathrm{u}%
})^{\mathrm{ac}} = 0$, it follows that $\Omega^{(r)}_{\mathrm{u}}$ has no
absolutely continuous part.

Finally, by defining
\[
\Omega_{0}^{\prime}=e^{it H^{(l)}_{U}H^{(r)}_{U}} S e^{-it H^{(l)}_{U}%
H^{(r)}_{U}},
\]
one can show that the unitary part of $\Omega^{(l)}$ has no absolutely
continuous component by an analogous argument. This completes the proof of
Lemma \ref{no a.c.part}.
\end{proof}

\subsection{The Wold decomposition and the action of
\texorpdfstring{$\mathcal{H}^{(l)}$}{}}

\label{Lemma sub}

In this subsection, we give the proof of Theorem \ref{Wold reduce} used above.
The proof relies on the fact that the unitary parts of $\Omega^{(r)}$ and
$\Omega^{(l)}$ have no absolutely continuous components, as proved in the
preceding section in Lemma \ref{no a.c.part}. We also establish in Lemma
\ref{symmetryall} the symmetry relations for the anti-linear left and right
Hankel operators and their associated operators $\mathcal{H}^{(r)}$ and
$\mathcal{H}^{(l)}$, as well as the commutation relation needed in the proof.

\begin{lemma}
\label{symmetryall}  Let $U \in H^{2}_{n \times n}$, then  \newline(i) for all
$f \in\operatorname{Dom}H ^{(r)}_{U}$, $g \in\operatorname{Dom}H ^{(l)}_{U}$,
we have
\begin{align}
\label{symmetry-rl}\langle H ^{(r)}_{U}f,g\rangle_{H^{2}_{n \times n}}=\langle
H ^{(l)}_{U}g,f\rangle_{H^{2}_{n \times n}};
\end{align}
\newline(ii) let $\mathcal{H}^{(r)} =e^{itH_{U}^{(r)}H_{U}^{(l)}}H_{U}^{(r)}$
and $\mathcal{H}^{(l)} =e^{itH_{U}^{(l)}H_{U}^{(r)}}H_{U}^{(l)}$. We have
\begin{align}
\label{symmetry2}\langle\mathcal{H} ^{(r)}f,g\rangle_{H^{2}_{n \times n}}=
\langle\mathcal{H}^{(l)}g, f\rangle_{H^{2}_{n \times n}}%
\end{align}
for all $f \in\operatorname{Dom}\mathcal{H} ^{(r)}$, $g \in\text{Dom}%
\,\mathcal{H} ^{(l)}$; \newline(iii)
\begin{align}
\label{commutation-left}\mathcal{H}^{(l)} \Omega^{(r)}=\left( \Omega
^{(l)}\right) ^{\ast}\mathcal{H}^{(l)}%
\end{align}
holds on the domain of $\mathcal{H}^{(l)} $.
\end{lemma}

\begin{proof}
(i) If $f,g$ are polynomials, this is a direct calculation. Now, use Theorem
\ref{Domain dense left and right Hankel} to extend to all $f \in
\operatorname{Dom}H ^{(r)}_{U}$, $g \in\operatorname{Dom}H ^{(l)}_{U}$.

(ii) Using (\ref{symmetry-rl}), we find
\begin{align*}
\langle\mathcal{H} ^{(r)}f,g\rangle=  &  \langle e^{i t H^{(r)}_{U}H^{(l)}%
_{U}} H^{(r)}_{U}f,g\rangle\\
=  &  \langle H^{(r)}_{U}f,e^{-i t H^{(r)}_{U}H^{(l)}_{U}}g\rangle\\
=  &  \langle H ^{(l)}_{U}e^{-i t H^{(r)}_{U}H^{(l)}_{U}}g,f\rangle\\
=  &  \langle\mathcal{H}^{(l)}g,f\rangle
\end{align*}
for all $f \in\operatorname{Dom}\mathcal{H} ^{(r)}$, $g\in\text{Dom}%
\,\mathcal{H} ^{(l)}$.

(iii) By the relation
\begin{align*}
K^{(l)}_{U}= S ^{\ast} H^{(l)}_{U}=H^{(l)}_{U}S,
\end{align*}
we find
\begin{align*}
\mathcal{H}^{(l)}\Omega^{(r)}  &  =(H^{(l)}_{U}e^{-i t H^{(r)}_{U}H^{(l)}_{U}%
})(e^{it H^{(r)}_{U}H^{(l)}_{U}}Se^{-it K^{(r)}_{U}K^{(l)}_{U}})\\
&  =H^{(l)}_{U}Se^{-it K^{(r)}_{U}K^{(l)}_{U}} =K^{(l)}_{U}e^{-it K^{(r)}%
_{U}K^{(l)}_{U}}=e^{itK^{(l)}_{U}K^{(r)}_{U}}K^{(l)}_{U},
\end{align*}
\begin{align*}
\left( \Omega^{(l)}\right) ^{\ast}\mathcal{H}^{(l)}  &  =(e^{itK^{(l)}%
_{U}K^{(r)}_{U}}S^{\ast}e^{-itH^{(l)}_{U}H^{(r)}_{U}})(e^{i t H^{(l)}%
_{U}H^{(r)}_{U}} H^{(l)}_{U})\\
&  =e^{itK^{(l)}_{U}K^{(r)}_{U}}S^{\ast}H^{(l)}_{U} =e^{itK^{(l)}_{U}%
K^{(r)}_{U}}K^{(l)}_{U}.
\end{align*}
Then, we have (\ref{commutation-left}) holds on the domain of $\mathcal{H}%
^{(l)} $.
\end{proof}

We now have the tools that we need to prove the crucial Theorem
\ref{Wold reduce}.

\textbf{\textit{Proof of Theorem \ref{Wold reduce}:}} First, we discuss the
spectral measures associated with the isometry $\Omega^{(r)}$ (and
similarly for $\Omega^{(l)}$). Recall that $X:=H^{2}_{n\times n}$. Let
$P^{(r)}_{\mathrm{u}}$ (resp. $P^{(r)}_{\mathrm{pure}}$ ) denote the
orthogonal projection from $X$ onto $X^{(r)}_{\mathrm{u}}$ (resp.
$X^{(r)}_{\mathrm{pure}}$) in the Wold decomposition
(\ref{Wold decomposition u part}) of $\Omega^{(r)}$. For $f ,g\in X$ and $m
\geq0$, we have
\begin{align}
\label{inner-Wold}\left\langle \left( \Omega^{(r)}\right) ^{m}f,
g\right\rangle  &  =\left\langle \left( \Omega^{(r)}\right) ^{m}
(P^{(r)}_{\mathrm{u}} f + P^{(r)}_{\mathrm{pure}} f), \, P^{(r)}_{\mathrm{u}}
g + P^{(r)}_{\mathrm{pure}} g\right\rangle \nonumber\\
&  = \left\langle \left( \Omega_{\mathrm{u}}^{(r)}\right) ^{m} P^{(r)}%
_{\mathrm{u}} f, P^{(r)}_{\mathrm{u}} g\right\rangle +\left\langle \left(
\Omega_{\mathrm{pure}}^{(r)}\right) ^{m} P^{(r)}_{\mathrm{pure}} f,
P^{(r)}_{\mathrm{pure}} g\right\rangle .
\end{align}
Indeed, for $m \ge0$,
\[
\left\langle \left( \Omega_{\mathrm{u}}^{(r)}\right) ^{m} P^{(r)}_{\mathrm{u}}
f,\, P^{(r)}_{\mathrm{pure}} g \right\rangle = 0, \qquad\left\langle \left(
\Omega_{\mathrm{pure}}^{(r)}\right) ^{m} P^{(r)}_{\mathrm{pure}} f,\,
P^{(r)}_{\mathrm{u}} g \right\rangle = 0,
\]
because $\left( \Omega_{\mathrm{u}}^{(r)}\right) ^{m} P^{(r)}_{\mathrm{u}} f
\in X^{(r)}_{\mathrm{u}}$, $P^{(r)}_{\mathrm{pure}} g \in X^{(r)}%
_{\mathrm{pure}}$, $\left( \Omega_{\mathrm{pure}}^{(r)}\right) ^{m}
P^{(r)}_{\mathrm{pure}} f \in X^{(r)}_{\mathrm{pure}}$, $P^{(r)}_{\mathrm{u}}
g \in X^{(r)}_{\mathrm{u}}$, and $X^{(r)}_{\mathrm{u}} \perp X^{(r)}%
_{\mathrm{pure}}$.

Then we write
\begin{align}
\label{write-u}\left\langle \left( \Omega_{\mathrm{u}}^{(r)}\right) ^{m}
P^{(r)}_{\mathrm{u}} f, P^{(r)}_{\mathrm{u}} g\right\rangle =\int_{\mathbb{T}}
z^{m} d \nu_{f, g}^{(r)}(z),
\end{align}
where $\nu_{f, g}^{(r)}$ denotes the complex spectral measure of $\Omega
^{(r)}_{\mathrm{u}}$, associated with the pair of vectors $P_{\mathrm{u}} f,
P_{\mathrm{u}} g$. By Lemma \ref{no a.c.part}, $\Omega^{(r)}_{\mathrm{u}}$ has
no absolutely continuous component. Hence, the measure $\nu_{f, g}^{(r)}$ is
singular with respect to the Lebesgue measure.

In order to handle the second term on the right hand side of (\ref{inner-Wold}%
), we use the Wold decomposition theorem. By Theorem \ref{Wold decomposition},
$\Omega^{(r)}_{\mathrm{pure}}$ is unitarily equivalent to the shift operator
$S$ on $H^{2}_{n \times n}$, i.e. there exists a unitary map $\mathbf{V}%
:X^{(r)}_{\mathrm{pure}} \to H^{2}_{n \times n}$ such that
\[
\mathbf{V} \Omega^{(r)}_{\mathrm{pure}} \mathbf{V}^{\ast}= S.
\]
Equivalently, for any $m \ge0$,
\[
\mathbf{V} \left( \Omega_{\mathrm{pure}}^{(r)}\right) ^{m} = S^{m}
\mathbf{V}.
\]
Set
\[
f_{\ast}:= \mathbf{V}(P^{(r)}_{\mathrm{pure}} f), \quad g_{\ast}:=
\mathbf{V}(P^{(r)}_{\mathrm{pure}} g) \in H^{2}_{n \times n}.
\]
Then
\[
\mathbf{V} \left( \Omega_{\mathrm{pure}}^{(r)}\right) ^{m} P^{(r)}%
_{\mathrm{pure}} f = S^{m} f_{\ast},\quad\left\langle \left( \Omega
_{\mathrm{pure}}^{(r)}\right) ^{m} P^{(r)}_{\mathrm{pure}} f, P^{(r)}%
_{\mathrm{pure}} g \right\rangle = \langle S^{m} f_{\ast}, g_{\ast}%
\rangle_{H^{2}_{n \times n}}.
\]
Recalling the inner product,
\[
\langle F, G \rangle_{H^{2}_{n \times n}} = \int_{\mathbb{T}}
\operatorname{tr}(F(z) G(z)^{\ast}) \, d\lambda(z),
\]
we obtain
\begin{align}
\label{write-pure}\left\langle \left( \Omega_{\mathrm{pure}}^{(r)}\right) ^{m}
P^{(r)}_{\mathrm{pure}} f, P^{(r)}_{\mathrm{pure}} g \right\rangle =
\int_{\mathbb{T}} z^{m} \operatorname{tr}(f_{\ast}(z) g_{\ast}(z)^{\ast}) \,
d\lambda(z) = \int_{\mathbb{T}} z^{m} \, d\mu^{(r)}_{f,g}(z),
\end{align}
where $\mu^{(r)}_{f, g}$ is the absolutely continuous measure on $\mathbb{T}$
defined by the above identity. Combining equations (\ref{write-u}) and
(\ref{write-pure}), we obtain
\begin{align}
\label{inner-sum}\left\langle \left( \Omega^{(r)} \right) ^{m} f,
g\right\rangle =\int_{\mathbb{T}} z^{m} d \nu_{f, g}^{(r)}(z)+\int
_{\mathbb{T}} z^{m} d \mu_{f, g}^{(r)}(z), \quad m \geq0 ,
\end{align}
where $\nu_{f, g}^{(r)}$ is singular and $\mu_{f, g}^{(r)}$ is absolutely
continuous with respect to the Lebesgue measure.

Similarly, by Lemma \ref{no a.c.part}, $\Omega^{(l)}_{\mathrm{u}}$ has no
absolutely continuous component. Hence, the measure $\nu_{f, g}^{(l)}$ is
singular with respect to the Lebesgue measure. The same argument yields
\begin{align}
\label{inner-sum1}\left\langle \left( \Omega^{(l)} \right) ^{m} f,
g\right\rangle =\int_{\mathbb{T}} z^{m} d \nu_{f, g}^{(l)}(z)+\int
_{\mathbb{T}} z^{m} d \mu_{f, g}^{(l)}(z), \quad m \geq0 ,
\end{align}
where $\nu_{f, g}^{(l)}$ is singular and $\mu_{f, g}^{(l)}$ is absolutely
continuous with respect to the Lebesgue measure.

Now, taking $f \in\operatorname{Dom}\mathcal{H}^{(r)}$ and $g
\in\operatorname{Dom}\mathcal{H}^{(l)}$, we aim to prove the two statements
(\ref{eqdom}) and (\ref{eqpure}).

We find that
\begin{align}
\label{Omegadom}\Omega^{(r)}\operatorname{Dom}\mathcal{H} ^{(l)}%
\subset\operatorname{Dom}\mathcal{H} ^{(l)}.
\end{align}
By $g \in\operatorname{Dom}\mathcal{H} ^{(l)}$, that is $g\in H^{2}_{n \times
n}$, and $\mathcal{H} ^{(l)}g\in H^{2}_{n \times n}$. Then $\Omega^{(r)} g\in
H^{2}_{n\times n}$, and by equation (\ref{commutation-left}), we have
$\mathcal{H}^{(l)} \Omega^{(r)} g=\left( \Omega^{(l)}\right) ^{\ast
}\mathcal{H}^{(l)}g\in H^{2}_{n \times n}$, that is $\Omega^{(r)}
g\in\operatorname{Dom}\mathcal{H} ^{(l)}$. Thus (\ref{Omegadom}) holds.

Iterating the commutation relation (\ref{commutation-left}) in Lemma
\ref{symmetryall}, we obtain
\[
\mathcal{H}^{(l)} \left(  \Omega^{(r)}\right) ^{m}=\left( \left( \Omega
^{(l)}\right) ^{\ast}\right) ^{m}\mathcal{H}^{(l)}
\]
for all $m \geq0$. Applying the inner product to this identity, we obtain
\begin{align}
\label{left hand}\left\langle f,\mathcal{H}^{(l)} \left(  \Omega^{(r)}\right)
^{m} g\right\rangle =\left\langle f,\left( \left( \Omega^{(l)}\right) ^{\ast
}\right) ^{m}\mathcal{H}^{(l)} g\right\rangle =\left\langle \left(
\Omega^{(l)}\right) ^{m}f,\mathcal{H}^{(l)}g \right\rangle .
\end{align}

Taking complex conjugates on (\ref{symmetry2}) yields
\begin{align}
\label{symmetry (a)}\langle g,\mathcal{H} ^{(r)} f\rangle= \langle
f,\mathcal{H}^{(l)}g \rangle.
\end{align}
By the symmetry relation (\ref{symmetry (a)}), the left hand side of
(\ref{left hand}) can be written as
\[
\left\langle f, \mathcal{H}^{(l)} \left( \Omega^{(r)}\right) ^{m}
g\right\rangle =\left\langle \left(  \Omega^{(r)}\right) ^{m} g,\mathcal{H}%
^{(r)} f \right\rangle .
\]
Combining this with (\ref{left hand}), we obtain
\begin{align}
\label{inner-equal}\left\langle \left(  \Omega^{(r)}\right) ^{m} g,
\mathcal{H}^{(r)} f\right\rangle =\left\langle \left( \Omega^{(l)}\right) ^{m}
f , \mathcal{H}^{(l)} g \right\rangle .
\end{align}
By (\ref{inner-sum}), (\ref{inner-sum1}) and (\ref{inner-equal}), this yields
\begin{align*}
\int_{\mathbb{T}} z^{m} d \nu^{(r)}_{g, \mathcal{H}^{(r)} f}(z)+\int
_{\mathbb{T}} z^{m} d \mu^{(r)}_{g, \mathcal{H}^{(r)} f}(z)=\int_{\mathbb{T}}
z^{m} d \nu^{(l)}_{f, \mathcal{H}^{(l)} g}(z)+\int_{\mathbb{T}} z^{m} d
\mu^{(l)}_{f, \mathcal{H}^{(l)} g}(z)
\end{align*}
for all $m \geq0$.

Then
\begin{align*}
\int_{\mathbb{T}} z^{m} d (\nu^{(r)}_{g, \mathcal{H}^{(r)} f}-\nu^{(l)}_{f,
\mathcal{H}^{(l)} g})(z)=\int_{\mathbb{T}} z^{m} d (\mu^{(l)}_{f,
\mathcal{H}^{(l)} g}-\mu^{(r)}_{g, \mathcal{H}^{(r)} f})(z).
\end{align*}

By the F. and M. Riesz theorem, the singular parts of the two measures
coincide, namely,
\[
\nu^{(r)}_{g, \mathcal{H}^{(r)} f}=\nu^{(l)}_{f, \mathcal{H}^{(l)} g}.
\]
In particular,
\[
\nu^{(r)}_{g, \mathcal{H}^{(r)} f}(\mathbb{T})=\nu^{(l)}_{f, \mathcal{H}^{(l)}
g}(\mathbb{T});
\]
Using the definition (\ref{write-u}) of the measure $\nu^{(r)}_{f,g}$ (an
analogous definition applies to $\nu_{f,g}^{(l)}$), we obtain
\begin{align*}
\langle P^{(r)}_{\mathrm{u}} g, P^{(r)}_{\mathrm{u}} \mathcal{H}^{(r)}%
f\rangle=\langle P^{(l)}_{\mathrm{u}} f, P^{(l)}_{\mathrm{u}} \mathcal{H}%
^{(l)}g\rangle.
\end{align*}
Equivalently, since $P^{(r)}_{\mathrm{u}}$ and $P^{(l)}_{\mathrm{u}}$ are
orthogonal projections, taking complex conjugates gives
\begin{align}
\label{another}\langle\mathcal{H}^{(r)}f, P^{(r)}_{\mathrm{u}} g\rangle
=\langle P^{(l)}_{\mathrm{u}} \mathcal{H}^{(l)}g, f \rangle.
\end{align}

\textbf{On the other hand, by Theorem \ref{Domain dense left and right Hankel}
(iv) and } $\operatorname{Dom}\mathcal{H} ^{(r)}=\operatorname{Dom}
H^{(r)}_{U}$, there exists $h_{0}\in H^{2}_{n\times n}$ such that
\[
\langle\mathcal{H}^{(r)}f,P^{(r)}_{\mathrm{u}}g\rangle= \langle h_{0}%
,f\rangle,\quad\forall f \in\operatorname{Dom}\mathcal{H} ^{(r)},
\]
then $P^{(r)}_{\mathrm{u}} g\in\operatorname{Dom}\mathcal{H}^{(l)}$ and
$h_{0}=\mathcal{H}^{(l)}P^{(r)}_{\mathrm{u}} g$.

Combining the above equation with (\ref{another}), we can obtain
\[
\mathcal{H}^{(l)}P^{(r)}_{\mathrm{u}} g=P^{(l)}_{\mathrm{u}} \mathcal{H}%
^{(l)}g.
\]

Since
\begin{align*}
\mathcal{H}^{(l)}P^{(r)}_{\mathrm{pure}}g  &  = \mathcal{H}^{(l)}%
(g-P^{(r)}_{\mathrm{u}}g)\\
&  = \mathcal{H}^{(l)}g-\mathcal{H}^{(l)}P^{(r)}_{\mathrm{u}}g
\end{align*}
and $g \in\operatorname{Dom}\mathcal{H}^{(l)}$, we also have $P^{(r)}%
_{\mathrm{pure}}g\in\operatorname{Dom}\mathcal{H}^{(l)}$, and furthermore
$\mathcal{H}^{(l)}P^{(r)}_{\mathrm{pure}} g=P^{(l)}_{\mathrm{pure}}
\mathcal{H}^{(l)}g $.

\section{Proof of the strong resolvent convergence of
\texorpdfstring{$H^{(r)}_{U}H^{(l)}_{U}$}{}}

\label{proof of strong convergence}

This section is devoted to the proof of Theorem \ref{strong con-H}. The main
idea is to first establish the strong resolvent convergence of $\Gamma
^{(r)}_{U}(\Gamma^{(r)}_{U})^{\ast}$, from which the corresponding result for
$H^{(r)}_{U}H^{(l)}_{U}$ follows immediately.

Let $U\in H^{2}_{n \times n}$ and let $\Gamma^{(r)}_{U}$ and $\Gamma^{(l)}%
_{U}$ be defined as in Theorems \ref{Domain dense} and \ref{Domain dense left}%
, respectively. Consider the operator $(\Gamma^{(r)}_{U})^{\ast}\Gamma
^{(r)}_{U}$ with the domain
\begin{align*}
\operatorname{Dom} (\Gamma^{(r)}_{U})^{\ast}\Gamma^{(r)}_{U}=\left\{  f \in\,
\operatorname{Dom} \Gamma^{(r)}_{U}: \Gamma^{(r)}_{U}f \in\operatorname{Dom}%
(\Gamma^{(r)}_{U})^{\ast}\right\}  .
\end{align*}
Since $\Gamma_{U}^{(r)}$ is closed, $(\Gamma_{U}^{(r)})^{\ast}\Gamma_{U}%
^{(r)}$ is a self-adjoint operator. Moreover, it is the self-adjoint
operator associated with the closed sesquilinear form
\begin{align*}
\langle\Gamma^{(r)}_{U}f,\Gamma^{(r)}_{U}h \rangle,
\end{align*}
where $f , h \in\operatorname{Dom} \Gamma^{(r)}_{U}$, see e.g. \cite[Section
VIII.6]{Reed80}.

A simple calculation with the first resolvent identity can obtain the
following result.

\begin{lemma}
\textup{\cite[Problem VIII.20]{Reed80}} \label{the first resolvent identity}
Let $\{A_{k}\}_{k \in\mathbb{N}} $ and $A $ be self-adjoint operators on the
Hilbert space $X$ and suppose that for all $f,g \in X$ and all $\lambda
\in\mathbb{C} $ with $\operatorname{Im}\,\lambda\neq0$,
\[
\langle(A_{k} - \lambda\, I )^{-1}f ,g\rangle\to\langle(A - \lambda\, I
)^{-1}f ,g\rangle.
\]
Then $A_{k}\to A $ in the strong resolvent sense.
\end{lemma}

The proof of the following theorem is similar to that of \cite[Theorem
2.2]{Gerard23}, using Lemma \ref{the first resolvent identity}.

\begin{theorem}
\label{strong con} Let $\{U_{k}\}_{k=1}^{\infty}$ be a sequence of elements of
$H^{2}_{n \times n}$ with $\Vert U_{k}-U\Vert_{H^{2}_{n \times n}}\to0$ as $k
\to\infty$. Then
\begin{align*}
(\Gamma^{(r)}_{U_{k}})^{\ast}\Gamma^{(r)}_{U_{k}} \to(\Gamma^{(r)}_{U})^{\ast
}\Gamma^{(r)}_{U}%
\end{align*}
and
\begin{align*}
\Gamma^{(r)}_{U_{k}} (\Gamma^{(r)}_{U_{k}})^{\ast} \to\Gamma^{(r)}_{U}%
(\Gamma^{(r)}_{U})^{\ast}%
\end{align*}
in the strong resolvent sense.
\end{theorem}

\begin{proof}
We only prove the convergence of $(\Gamma^{(r)}_{U_{k}})^{\ast}\Gamma
^{(r)}_{U_{k}}$; the corresponding convergence of $\Gamma^{(r)}_{U_{k}}%
(\Gamma^{(r)}_{U_{k}})^{\ast}$ then follows by replacing $U_{k}$ with
$U_{k}^{\#}$.  Fix $\lambda$ with $\text{Im} \lambda\neq0$ and $f \in H^{2}_{n
\times n}$; denote
\[
\varphi_{k}
=((\Gamma^{(r)}_{U_{k}})^{\ast}\Gamma^{(r)}_{U_{k}}-\lambda\, I )^{-1}f,
\quad\varphi= ((\Gamma^{(r)}_{U})^{\ast}\Gamma^{(r)}_{U}-\lambda\, I )^{-1}f.
\]
Our goal is to show that $\|\varphi_{k}-\varphi\| \to0$ as $k \to\infty$. By
Lemma \ref{the first resolvent identity}, it is enough to prove that
$\varphi_{k} \to\varphi$ weakly.

Since $\operatorname{Im}\lambda\neq0$, we have
\[
\|\varphi_{k}\| \leq C_{1}(\lambda)\|f\|.
\]
Next,
\begin{align*}
\|\Gamma^{(r)}_{U_{k}}\varphi_{k}\|^{2}-\lambda\|\varphi_{k}\|^{2}  & =
\langle(\Gamma^{(r)}_{U_{k}})^{\ast}\Gamma^{(r)}_{U_{k}}\varphi_{k}%
,\varphi_{k}\rangle- \lambda\|\varphi_{k}\|^{2}\\
& = \langle f,\varphi_{k}\rangle= \langle f,( (\Gamma^{(r)}_{U_{k}})^{\ast
}\Gamma^{(r)}_{U_{k}}-\lambda\, I )^{-1}f\rangle,
\end{align*}
and so
\begin{align*}
\|\Gamma^{(r)}_{U_{k}}\varphi_{k}\|^{2}  &  \leq\left|  \langle f,(\Gamma
^{(r)}_{U_{k}})^{\ast}\Gamma^{(r)}_{U_{k}}\varphi_{k},\varphi_{k}%
\rangle-\lambda\, I )^{-1}f\rangle\right|  + |\lambda|\|\varphi_{k}\|^{2}\\
&  \leq C_{2}(\lambda)\|f\|^{2}.
\end{align*}

Hence, we can select a subsequence such that $\varphi_{k} \to\tilde{\varphi}$
weakly and $\Gamma^{(r)}_{U_{k}} \varphi_{k} \to\psi$ weakly for some
$\tilde{\varphi}\in H^{2}_{n \times n}$, $\psi\in H^{2}_{n \times n}$.

We first show that $\tilde{\varphi} \in\text{Dom } \Gamma^{(r)}_{U}$
and $\Gamma^{(r)}_{U} \tilde{\varphi}=\psi$. For any polynomial $p_{1}%
\in\mathcal{P}_{n \times n}$, we have
\[
\langle\varphi_{k},(\Gamma^{(r)}_{U_{k}})^{\ast}p_{1}\rangle= \langle
\Gamma^{(r)}_{U_{k}}\varphi_{k},p_{1}\rangle\to\langle\psi,p_{1}\rangle.
\]
Since $p_{1} \in H^{\infty}_{n \times n}$, we also have
\[
\|(\Gamma^{(r)}_{U_{k}})^{\ast}p_{1} -(\Gamma^{(r)}_{U})^{\ast}p_{1}\| \to0,
\]
and so
\[
\langle\varphi_{k},(\Gamma^{(r)}_{U_{k}})^{\ast}p_{1}\rangle\to\langle
\tilde{\varphi},(\Gamma^{(r)}_{U})^{\ast}p_{1},\rangle.
\]
Therefore,
\[
\langle\tilde{\varphi},(\Gamma^{(r)}_{U})^{\ast}p_{1}\rangle= \langle
\psi,p_{1}\rangle
\]
for any polynomial $p_{1}$. By Theorem \ref{Domain dense}, polynomials are
dense in $\text{Dom } (\Gamma^{(r)}_{U})^{\ast}$ with respect to the graph
norm and therefore the last identity extends by continuity to all $p_{1}
\in\text{Dom } (\Gamma^{(r)}_{U})^{\ast}$. And by
\[
\langle\psi,p_{1}\rangle=\langle\tilde{\varphi},(\Gamma^{(r)}_{U})^{\ast}%
p_{1}\rangle= \langle\Gamma^{(r)}_{U}\tilde{\varphi},p_{1}\rangle,
\]
we have $\tilde{\varphi}\in\text{Dom } \Gamma^{(r)}_{U}$ and
$\Gamma^{(r)}_{U}\tilde{\varphi}=\psi$.

By the definition of $\varphi_{k}$, for any polynomial $p_{2}\in\mathcal{P}_{n
\times n}$ we have
\begin{align}
\label{identity used}\langle\Gamma^{(r)}_{U_{k}}\varphi_{k},\Gamma
^{(r)}_{U_{k}}p_{2}\rangle- \lambda\langle\varphi_{k},p_{2}\rangle= \langle
f,p_{2}\rangle.
\end{align}
From what has already been established, we have $\varphi_{k} \to\tilde
{\varphi}$ and $\Gamma^{(r)}_{U_{k}}\varphi_{k} \to\Gamma^{(r)}_{U}%
\tilde{\varphi}$ weakly. Since $p_{2} \in H^{\infty}_{n \times n}$, we also
have
\[
\|\Gamma^{(r)}_{U_{k}}p_{2} -\Gamma^{(r)}_{U}p_{2}\| \to0.
\]
Passing to the limit (\ref{identity used}) yields
\[
\langle\Gamma^{(r)}_{U}\tilde{\varphi},\Gamma^{(r)}_{U}p_{2}\rangle-
\lambda\langle\tilde{\varphi},p_{2}\rangle= \langle f,p_{2}\rangle.
\]
for all polynomials $p_{2}$. Applying Theorem \ref{Domain dense} once again,
we extend this identity from polynomials $p_{2}$ to all elements $p_{2}%
\in\text{Dom } \Gamma^{(r)}_{U}$. It follows that $\Gamma^{(r)}_{U}%
\tilde{\varphi}\in\text{Dom } (\Gamma^{(r)}_{U})^{\ast}$ and
\[
(\Gamma^{(r)}_{U})^{\ast}\Gamma^{(r)}_{U}\tilde{\varphi} = \lambda
\tilde{\varphi} + f.
\]
Equivalently,
\[
\tilde{\varphi} = ( (\Gamma^{(r)}_{U})^{\ast}\Gamma^{(r)}_{U}-\lambda\, I
)^{-1}f.
\]
By the definition of $\varphi$, the right-hand side is precisely $\varphi$.
Hence, we have proved that $\varphi_{k} \to\tilde{\varphi}=\varphi$
weakly along a subsequence. Therefore, $\varphi_{k}\rightharpoonup\varphi$
weakly. Indeed, otherwise one could find a subsequence that does not converge
weakly to $\varphi$, while the same argument yields a further subsequence
converging weakly to $\varphi$, a contradiction. This completes the proof of
Theorem \ref{strong con}.
\end{proof}

By Theorem \ref{strong con}, together with the identity
\[
H^{(r)}_{U}H ^{(l)}_{U}= \Gamma^{(r)}_{U}\mathcal{C}\Gamma^{(l)}%
_{U}\mathcal{C}= \Gamma^{(r)}_{U}\Gamma^{(r)}_{U^{\#}} =\Gamma^{(r)}%
_{U}(\Gamma^{(r)}_{U})^{\ast},
\]
we immediately obtain Theorem \ref{strong con-H}.

\appendix

\section{Appendix}

\label{Appendix}

\subsection{Some notions in the Corollary \ref{GuLuo Corol}}

\label{notions change} Here, Corollary \ref{GuLuo Corol} has some changes to
the notation that differ from those in \cite{GuLuo}. These changes
are presented in the Appendix \ref{notions change} for interested readers.

In this appendix, we introduce a change of notation and a definition of the
Hankel operator that differs from the one in \cite{GuLuo}; however, the two
definitions can be easily converted into each other. The details are given below.

Let $\Psi\in L^{\infty}_{B(E,F)}$. The Toeplitz operator $T_{\Psi}$ from
$H^{2}_{E}$ into $H^{2}_{F}$ is defined by
\begin{align*}
T_{\Psi}(f)=\mathbb{P}(\Psi f),f \in H^{2}_{E}.
\end{align*}
Let $J_{1}$ be defined on $L^{2}_{E}$ by
\begin{align*}
J_{1} f(z)=\bar{z}f(\bar{z}), f(z) \in L^{2}_{E}.
\end{align*}
The $J_{1}$ maps $\overline{z H^{2}_{E}}$ onto $H^{2}_{E}$, and $J_{1}$ maps
$H^{2}_{E}$ onto $\overline{z H^{2}_{E}}$.

The Hankel operator $H_{\Psi}$ from $H^{2}_{E}$ into $H^{2}_{F}$ is defined
by
\begin{align*}
H_{\Psi}(f)=\mathbb{P}J_{1}(\Psi f),f \in H^{2}_{E}.
\end{align*}

Then if we let $\Psi\in L^{\infty}_{n \times n}$; the Toeplitz operator
$T_{\Psi}$ from $H^{2}_{n \times1}$ into $H^{2}_{n \times1}$ is defined by
\begin{align*}
T_{\Psi}(f)=\mathbb{P}(\Psi f),f \in H^{2}_{n \times1}.
\end{align*}
Let $J_{1}$ be defined on $L^{2}_{n \times1}$ by
\begin{align*}
J_{1} f(z)=\bar{z}f(\bar{z}), f(z) \in L^{2}_{n \times1}.
\end{align*}

The Hankel operator $H_{\Psi}$ from $H^{2}_{n \times1}$ into $H^{2}_{n
\times1}$ is defined by
\begin{align*}
H_{\Psi}(f)=\mathbb{P}J_{1}(\Psi f),f \in H^{2}_{n \times1}.
\end{align*}

For convenience, we rewrite the notation in Corollary
\ref{GuLuo Corol} in the form used in the proof of Theorem \ref{Domain dense}.
Let
\begin{align*}
\Phi(z)=\left[
\begin{array}
[c]{cc}%
A(z) & B(z)\\
C(\overline{z}) & D(\overline{z})
\end{array}
\right]  ,
\end{align*}
where $A ,B\in H^{\infty}_{n \times n}$, $C(\overline{z}), D(\overline{z}) \in
L^{\infty}_{n \times n}$.
\begin{align*}
W_{\Phi}=\left[
\begin{array}
[c]{cc}%
T_{A} & T_{B}\\
H_{C(\overline{z})} & H_{D(\overline{z})}%
\end{array}
\right]  ,
\end{align*}
For $\left[
\begin{array}
[c]{c}%
h_{1}\\
h_{2}%
\end{array}
\right]  \in H^{2}_{n \times1}\oplus H^{2}_{n \times1}$, where $h_{1} ,h_{2}
\in H^{2}_{n \times1}$.

In Corollary \ref{GuLuo Corol}, if
\[%
\begin{bmatrix}
\alpha\\
\beta
\end{bmatrix}
\in\operatorname{Ran}(W_{\Phi}),
\]
then there exist $h_{1},h_{2}\in H^{2}_{n\times1}$ such that
\[%
\begin{bmatrix}
\alpha\\
\beta
\end{bmatrix}
= W_{\Phi}
\begin{bmatrix}
h_{1}\\
h_{2}%
\end{bmatrix}
=
\begin{bmatrix}
T_{A}h_{1}+T_{B}h_{2}\\
H_{C(\overline{z})}h_{1}+H_{D(\overline{z})}h_{2}%
\end{bmatrix}
.
\]
Hence,
\begin{equation}
\label{Range1}\alpha= T_{A}h_{1}+T_{B}h_{2},
\end{equation}
and
\begin{equation}
\label{Range2}\beta= H_{C(\overline{z})}h_{1} + H_{D(\overline{z})}h_{2}.
\end{equation}

Let us rewrite (\ref{Range1}) and (\ref{Range2}) in another form, we have
\begin{align*}
H_{C(\overline{z})}h_{1}(z)=\mathbb{P}(\overline{z}C(z) h_{1}(\overline{z})
)=\mathbb{P}((\overline{z}C) Jh_{1} ) =\Gamma_{\overline{z}C}h_{1}.
\end{align*}
Then
\begin{align}
\label{Range4}\alpha &  =Ah_{1} + Bh_{2}\\
\beta &  =\Gamma_{\overline{z}C }h_{1} +\Gamma_{\overline{z}D }h_{2}.
\end{align}
The condition on $\Phi$ in Corollary \ref{GuLuo Corol} will not play a role in
our proof of Theorem \ref{Domain dense}(ii) below; all that matters to us is
that $A, B \in H^{\infty}_{n \times n}$, $\overline{z}C , \overline{z}D \in
L^{\infty}_{n \times n}$.

\subsection{Proof of Theorem \ref{Sun lax pair} and Theorem
\ref{Sun explicit formula}}

\label{Alternative proof}

In \cite[Theorem 1.3]{Sun}, Sun established the two-Lax pair
structure of the cubic matrix Szeg\H{o} equation (\ref{cubic matrix}), which
we refer to as Theorem \ref{Sun lax pair}. Although all the results in this
appendix are contained in the more general framework of \cite{Sun}, the more
specialized setting considered here allows for simpler notation and more
direct proofs. Therefore, for the reader's convenience and for completeness,
we present self-contained proofs of all the lemmas in this appendix by
adapting the arguments used for the cubic scalar Szeg\H{o} equation in
\cite{Gerard15} to the matrix-valued setting.

To facilitate the proof of Theorem \ref{Sun lax pair}, we present the
following relations between Hankel and Toeplitz operators.

\begin{lemma}
\textup{\cite[Lemma 3.11]{Sun}} \label{Hankel-Toep relation} Let
$U\in\mathcal{H}_{+}^{s}(\mathbb{T},\mathbb{C}^{n\times n})$ for $s > 1/2$.
Then
\begin{align*}
H^{(r)}_{\mathbb{P}(UU^{\ast}U)}  &  =T^{(r)}_{UU^{\ast}} H^{(r)}_{U}%
+H^{(r)}_{U}T^{(l)}_{U^{\ast}U}-H^{(r)}_{U}H^{(l)}_{U}H^{(r)}_{U}%
:H^{2}_{n\times n}\to H^{2}_{n\times n},\\
H^{(l)}_{\mathbb{P}(UU^{\ast}U)}  &  =T^{(l)}_{U^{\ast}U} H^{(l)}_{U}%
+H^{(l)}_{U}T^{(r)}_{UU^{\ast}}-H^{(l)}_{U}H^{(r)}_{U}H^{(l)}_{U}%
:H^{2}_{n\times n}\to H^{2}_{n\times n}.
\end{align*}

\end{lemma}

\begin{proof}
We only give the proof of the first equality. Given $f\in H^{2}_{n\times n}$,
\begin{align*}
H^{(r)}_{\mathbb{P}(UU^{\ast}U)}f  &  = \mathbb{P}( \mathbb{P}(UU^{\ast
}U)f^{\ast})= \mathbb{P}( UU^{\ast}Uf^{\ast})+\mathbb{P}( (I-\mathbb{P}%
)(UU^{\ast}U)f^{\ast})\\
&  =\mathbb{P}( UU^{\ast}Uf^{\ast})\\
&  =\mathbb{P}( UU^{\ast}\mathbb{P}(Uf^{\ast}))+ \mathbb{P}( UU^{\ast
}((I-\mathbb{P})(Uf^{\ast}))\\
&  =T^{(r)}_{UU^{\ast}} H^{(r)}_{U}f+H^{(r)}_{U}h ,
\end{align*}
where $h =\overline{U^{\ast}(I-\mathbb{P})(Uf^{\ast})}=\overline{U^{\ast
}Uf^{\ast}}-\overline{U^{\ast}\mathbb{P}(Uf^{\ast})}$.

Since $h\in H^{2}_{n\times n}$,
\begin{align*}
h  &  = \mathbb{P}(h) = \mathbb{P}\left(  (U^{\ast} Uf^{\ast})^{\ast}\right)
- \mathbb{P}(\overline{U^{\ast} \mathbb{P}(U f^{\ast})})\\
&  = \mathbb{P}(f U^{\ast} U) - \mathbb{P}(\overline{\mathbb{P}(U f^{\ast})}
U)\\
&  = \mathbb{P}(f U^{\ast} U) - H_{U}^{(l)} H_{U}^{(r)} f\\
&  = T_{U^{\ast} U}^{(l)} f - H_{U}^{(l)} H_{U}^{(r)} f.
\end{align*}
This completes the proof.
\end{proof}

Using Lemma \ref{Hankel-Toep relation} and the property of shift operators, we
have the following result.

\begin{lemma}
\textup{\cite[Lemma 3.11]{Sun}} \label{K-relation} Let $U\in\mathcal{H}%
_{+}^{s}(\mathbb{T},\mathbb{C}^{n\times n})$ for $s > 1/2$. Then
\begin{align*}
K^{(r)}_{\mathbb{P}(UU^{\ast}U)}  &  =T^{(r)}_{UU^{\ast}} K^{(r)}_{U}%
+K^{(r)}_{U}T^{(l)}_{U^{\ast}U}-K^{(r)}_{U}K^{(l)}_{U}K^{(r)}_{U}%
:H^{2}_{n\times n}\to H^{2}_{n\times n},\\
K^{(l)}_{\mathbb{P}(UU^{\ast}U)}  &  =T^{(l)}_{U^{\ast}U} K^{(l)}_{U}%
+K^{(l)}_{U}T^{(r)}_{UU^{\ast}}-K^{(l)}_{U}K^{(r)}_{U}K^{(l)}_{U}%
:H^{2}_{n\times n}\to H^{2}_{n\times n}.
\end{align*}

\end{lemma}

\begin{proof}
We observe that
\begin{align}
K_{\mathbb{P}(UU^{\ast}U)}^{(r)}  &  =H_{\mathbb{P}(UU^{\ast}U)}%
^{(r)}S\label{H-shift}\\
&  =T_{UU^{\ast}}^{(r)}H_{U}^{(r)}S +H_{U}^{(r)}T_{U^{\ast}U}^{(l)}S
-H_{U}^{(r)}H_{U}^{(l)}H_{U}^{(r)}S\nonumber\\
&  =T_{UU^{\ast}}^{(r)}K_{U}^{(r)} +H_{U}^{(r)} \left( T_{U^{\ast}U}^{(l)}S -
H_{U}^{(l)}K_{U}^{(r)} \right)  .\nonumber
\end{align}
Moreover, it is well known that the left Toeplitz operator $T_{U^{\ast}%
U}^{(l)}$ satisfy
\[
S ^{\ast}T_{U^{\ast}U}^{(l)}S =T_{U^{\ast}U}^{(l)}.
\]
By $S S ^{\ast}=I-\Delta$, then for given $f\in H^{2}_{n\times n}$,
\begin{align*}
[S, T^{(l)}_{U^{\ast}U}] f  & = S T_{U^{\ast}U}^{(l)}f-T_{U^{\ast}U}^{(l)}S
f\\
& = -\Delta T_{U^{\ast}U}^{(l)}S f\\
& =-\sum_{k,j=1}^{n}\langle\mathbb{P}(S f U^{\ast}U),\mathbb{E}_{kj}
\rangle_{H^{2}_{n\times n}}\mathbb{E}_{kj}\\
&  =-\sum_{k,j}\langle S f U^{\ast}U,\mathbb{E}_{kj} \rangle_{L^{2}_{n\times
n}}\mathbb{E}_{kj} .
\end{align*}
And
\begin{align*}
\sum_{k,j}\langle S f U^{\ast}U,\mathbb{E}_{kj}\rangle_{L^{2}_{n\times n}}  &
=\sum_{k,j}\langle U(\mathbb{E}_{kj} )^{\ast},U(Sf)^{\ast}\rangle\\
&  =\sum_{k,j}\langle H_{U}^{(r)}\mathbb{E}_{kj} ,H_{U}^{(r)}(Sf)\rangle\\
&  =\sum_{k,j}\langle H_{U}^{(r)}\mathbb{E}_{kj} ,K_{U}^{(r)}f\rangle\\
&  =\sum_{k,j}\langle H_{U}^{(l)}K_{U}^{(r)}f,\mathbb{E}_{kj} \rangle
_{H^{2}_{n\times n}},
\end{align*}
where the last equality follows from Theorem
\ref{Domain dense left and right Hankel}(iv). Consequently,
\begin{align*}
[S, T^{(l)}_{U^{\ast}U}] f =-\sum_{k,j}\langle H_{U}^{(l)}K_{U}^{(r)}%
f,\mathbb{E}_{kj} \rangle_{H^{2}_{n\times n}} \mathbb{E}_{kj} .
\end{align*}
That is
\begin{align*}
[S, T^{(l)}_{U^{\ast}U}]  & =-\Delta H_{U}^{(l)}K_{U}^{(r)}\\
&  =(S S ^{\ast}-I)H_{U}^{(l)}K_{U}^{(r)}\\
&  = SK_{U}^{(l)}K_{U}^{(r)} - H_{U}^{(l)}K_{U}^{(r)}.
\end{align*}
Coming back to (\ref{H-shift}), we obtain
\begin{align*}
K_{\mathbb{P}(UU^{\ast}U)}^{(r)}  &  =T_{UU^{\ast}}^{(r)}K_{U}^{(r)}
+H_{U}^{(r)} \left( T_{U^{\ast}U}^{(l)}S - H_{U}^{(l)}K_{U}^{(r)} \right) \\
&  =T_{UU^{\ast}}^{(r)}K_{U}^{(r)} +H_{U}^{(r)} \left( ST_{U^{\ast}U}^{(l)} -
S K_{U}^{(l)}K_{U}^{(r)} \right) \\
& =T^{(r)}_{UU^{\ast}} K^{(r)}_{U}+K^{(r)}_{U}T^{(l)}_{U^{\ast}U}-K^{(r)}%
_{U}K^{(l)}_{U}K^{(r)}_{U} .
\end{align*}

\end{proof}

\textit{\textbf{Proof of Theorem \ref{Sun lax pair}:}} Here, we only provide
the proof of formula (\ref{Lax-pair}) using Lemma \ref{Hankel-Toep relation}.
The remaining three formulas can be derived in a similar way by applying Lemma
\ref{Hankel-Toep relation} together with Lemma \ref{K-relation}.

Using equation (\ref{cubic matrix}) and Lemma \ref{Hankel-Toep relation},
\begin{align*}
\partial_{t}H_{U}^{(r)}  &  =H_{\partial_{t}U}^{(r)} =H_{-i\mathbb{P}%
(UU^{\ast}U)}^{(r)} =-iH_{\mathbb{P}(UU^{\ast}U)}^{(r)}\\
&  =-i(T^{(r)}_{UU^{\ast}} H^{(r)}_{U}+H^{(r)}_{U}T^{(l)}_{U^{\ast}U}%
-H^{(r)}_{U}H^{(l)}_{U}H^{(r)}_{U}),\\
\partial_{t}H_{U}^{(l)}  &  =-iH_{\mathbb{P}(UU^{\ast}U)}^{(l)}\\
&  =-i(T^{(l)}_{U^{\ast}U} H^{(l)}_{U}+H^{(l)}_{U}T^{(r)}_{UU^{\ast}}%
-H^{(l)}_{U}H^{(r)}_{U}H^{(l)}_{U}).
\end{align*}
Using the anti-linearity of $H_{U}^{(r)}$, then
\begin{align*}
\frac{d}{dt}(H_{U}^{(r)} H_{U}^{(l)})  &  = (\partial_{t}H_{U}^{(r)}%
)H_{U}^{(l)}+H_{U}^{(r)} (\partial_{t}H_{U}^{(l)})\\
&  =-i\left(  T^{(r)}_{UU^{\ast}} H^{(r)}_{U}+H^{(r)}_{U}T^{(l)}_{U^{\ast}%
U}-H^{(r)}_{U}H^{(l)}_{U}H^{(r)}_{U}\right)  H_{U}^{(l)}\\
&  \quad+H_{U}^{(r)}\left(  -i(T^{(l)}_{U^{\ast}U} H^{(l)}_{U}+H^{(l)}%
_{U}T^{(r)}_{UU^{\ast}}-H^{(l)}_{U}H^{(r)}_{U}H^{(l)}_{U})\right) \\
&  =-i\left(  T^{(r)}_{UU^{\ast}} H^{(r)}_{U}+H^{(r)}_{U}T^{(l)}_{U^{\ast}%
U}-H^{(r)}_{U}H^{(l)}_{U}H^{(r)}_{U}\right)  H_{U}^{(l)}\\
&  \quad+iH_{U}^{(r)}\left(  T^{(l)}_{U^{\ast}U} H^{(l)}_{U}+H^{(l)}%
_{U}T^{(r)}_{UU^{\ast}}-H^{(l)}_{U}H^{(r)}_{U}H^{(l)}_{U}\right) \\
&  =iH^{(r)}_{U}H^{(l)}_{U}T^{(r)}_{UU^{\ast}}-iT^{(r)}_{UU^{\ast}}H^{(r)}%
_{U}H^{(l)}_{U}\\
&  =i[H^{(r)}_{U}H^{(l)}_{U},T^{(r)}_{UU^{\ast}}].
\end{align*}
This completes the proof of formula (\ref{Lax-pair}) in Theorem
\ref{Sun lax pair}.

In \cite[Lemma 3.11]{Sun}, the following commutator formula is established.
Using Lemmas \ref{Hankel-Toep relation} and \ref{K-relation}, the
corresponding result in the particular case considered in \cite{Sun} follows immediately.

\begin{lemma}
\label{commutator-STH} Notations as above. Then
\[
[S ^{\ast} ,T_{UU^{\ast}}^{(r)} ]H^{(r)}_{U}=K^{(r)}_{U}(H^{(l)}_{U}
H^{(r)}_{U}-K^{(l)}_{U}K^{(r)}_{U}):H^{2}_{n\times n}\to H^{2}_{n\times n}.
\]

\end{lemma}

\begin{proof}
By Lemma \ref{Hankel-Toep relation} and $K^{(r)}_{U} =S^{\ast} H^{(r)}_{U}$,
\begin{align*}
K^{(r)}_{\mathbb{P}(UU^{\ast}U)}  &  = S^{\ast} H^{(r)}_{\mathbb{P}(UU^{\ast
}U)}\\
&  =S^{\ast}T^{(r)}_{UU^{\ast}} H^{(r)}_{U}+S^{\ast}H^{(r)}_{U}T^{(l)}%
_{U^{\ast}U}-S^{\ast}H^{(r)}_{U}H^{(l)}_{U}H^{(r)}_{U}\\
&  =S^{\ast}T^{(r)}_{UU^{\ast}} H^{(r)}_{U}+K^{(r)}_{U}T^{(l)}_{U^{\ast}%
U}-K^{(r)}_{U}H^{(l)}_{U}H^{(r)}_{U} .
\end{align*}
Combining the equality above with Lemma \ref{K-relation}, we obtain the
following result.
\begin{align*}
S^{\ast}T^{(r)}_{UU^{\ast}} H^{(r)}_{U}+K^{(r)}_{U}T^{(l)}_{U^{\ast}U}%
-K^{(r)}_{U}H^{(l)}_{U}H^{(r)}_{U}= T^{(r)}_{UU^{\ast}} K^{(r)}_{U}%
+K^{(r)}_{U}T^{(l)}_{U^{\ast}U}-K^{(r)}_{U}K^{(l)}_{U}K^{(r)}_{U},
\end{align*}
then
\begin{align*}
S^{\ast}T^{(r)}_{UU^{\ast}} H^{(r)}_{U}-K^{(r)}_{U}H^{(l)}_{U}H^{(r)}_{U}=
T^{(r)}_{UU^{\ast}} S^{\ast}H^{(r)}_{U}-K^{(r)}_{U}K^{(l)}_{U}K^{(r)}_{U},
\end{align*}
that is
\[
[S^{\ast} ,T_{UU^{\ast}}^{(r)} ]H^{(r)}_{U}=K^{(r)}_{U}(H^{(l)}_{U}
H^{(r)}_{U}-K^{(l)}_{U}K^{(r)}_{U}).
\]

\end{proof}

By applying formula (\ref{Lax-pair}) in Theorem \ref{Sun lax pair}, and
following the classical argument of Lax \cite{Lax68}, we obtain the following result.

\begin{lemma}
\textup{\cite[Corollary 4.1]{Sun}} \label{Sun-Corollary4.1} Let $U \in
C^{\infty}(\mathbb{R}; \mathcal{H}_{+}^{s}(\mathbb{T}, \mathbb{C}^{n \times
n})) $ , $s > 1/2$ be solution of equation (\ref{cubic matrix}), and let
$\mathbf{W} \in C^{1}(\mathbb{R}; B(H^{2}_{n \times n})) $ denote the unique
solution of the following equation:
\begin{align}
\label{wequation}\frac{\mathrm{d}}{\mathrm{d}t} \mathbf{W}(t) = -i T_{U(t)
U(t)^{\ast}}^{(r)} \mathbf{W}(t)
\end{align}
with initial data $\mathbf{W}(0) = I_{H^{2}_{n \times n}}$.

Then for any $t \in\mathbb{R} $, $\mathbf{W} \in C^{1}(\mathbb{R}; B(H^{2}_{n
\times n})) $ is unitary operator and the following identities of unitary
equivalences hold:
\begin{align}
H_{U(t)}^{(r)} H_{U(t)}^{(l)} & =\mathbf{W}(t)H_{U(0)}^{(r)} H_{U(0)}^{(l)}
\mathbf{W}(t)^{\ast},\label{unitary equivalences}\\
K_{U(t)}^{(r)} K_{U(t)}^{(l)} & =\mathbf{W}(t)K_{U(0)}^{(r)} K_{U(0)}^{(l)}
\mathbf{W}(t)^{\ast}.\nonumber
\end{align}

\end{lemma}

\begin{proof}
We only give the proof of the first identity.

Let $\mathcal{L}(t):= H^{(r)}_{U(t)}H^{(l)}_{U(t)}$, $t\in\mathbb{R}.$ In
particular, $\mathcal{L}(0)=H^{(r)}_{U(0)}H^{(l)}_{U(0)}$. Since
$\mathbf{W}(t)$ is unitary, we have $\mathbf{W}(t)^{\ast}= \mathbf{W}(t)^{-1}
$. To get%
\[
\mathcal{L}(t)=\mathbf{W}(t)\mathcal{L}(0) \mathbf{W}(t)^{\ast},
\]
we need to prove
\[
\mathbf{W}(t)^{-1} \mathcal{L}(t) \mathbf{W}(t)=\mathcal{L}(0).
\]
Indeed, differentiating $\mathbf{W}(t)^{-1}\mathcal{L}(t)\mathbf{W}(t)$ with
respect to $t$, we obtain
\begin{align*}
& \frac{d}{dt} \left(  \mathbf{W}(t)^{-1}\mathcal{L}(t)\mathbf{W}(t) \right)
\\
= &  \frac{d}{dt}\mathbf{W}(t)^{-1}\mathcal{L}(t)\mathbf{W}(t) +
\mathbf{W}(t)^{-1}\frac{d}{dt}\mathcal{L}(t)\mathbf{W}(t) + \mathbf{W}%
(t)^{-1}\mathcal{L}(t) \frac{d}{dt}\mathbf{W}(t).
\end{align*}

From (\ref{wequation}), we have
\[
\frac{d}{dt}\mathbf{W}(t)^{-1} =i\mathbf{W}(t)^{-1}T^{(r)}_{U(t)U(t)^{\ast}}
.
\]
Therefore,
\begin{align*}
& \frac{d}{dt} \left(  \mathbf{W}(t)^{-1}\mathcal{L}(t)\mathbf{W}(t) \right)
\\
=  &  i\mathbf{W}(t)^{-1} T^{(r)}_{U(t)U(t)^{\ast}} \mathcal{L}(t)\mathbf{W}%
(t) + \mathbf{W}(t)^{-1} \frac{d}{dt}\mathcal{L}(t) \mathbf{W}(t) -
i\mathbf{W}(t)^{-1} \mathcal{L}(t) T^{(r)}_{U(t)U(t)^{\ast}} \mathbf{W}(t)\\
=  &  \mathbf{W}(t)^{-1} \left(  T^{(r)}_{U(t)U(t)^{\ast}} \mathcal{L}(t)
+\frac{d}{dt}\mathcal{L}(t) - i \mathcal{L}(t) T^{(r)}_{U(t)U(t)^{\ast}}
\right)  \mathbf{W}(t).
\end{align*}

By the Lax pair identity (\ref{Lax-pair})
\[
\frac{d}{dt}\mathcal{L}(t) =i[\mathcal{L}(t),T^{(r)}_{U(t)U(t)^{\ast}}].
\]
Hence,
\[
\frac{d}{dt} \left(  \mathbf{W}(t)^{-1}\mathcal{L}(t)\mathbf{W}(t) \right)
=0.
\]
It follows that
\[
\mathbf{W}(t)^{-1}\mathcal{L}(t)\mathbf{W}(t) = \mathbf{W}(0)^{-1}%
\mathcal{L}(0)\mathbf{W}(0) = \mathcal{L}(0).
\]

\end{proof}

For the reader's convenience, we provide some additional details of the
argument in \cite{Sun}.

\begin{lemma}
\textup{\cite[Lemma 4.2]{Sun}}\label{Sun-Lemma4.2} As in the previous lemma,
let $U$ and $\mathbf{W}$ be as above. Then the following identity holds on
$H^{2}_{n \times n}$ for all $t\in\mathbb{R}$ :
\begin{equation}
\label{identity-1}\mathbf{W}(t)^{\ast} S^{\ast} \mathbf{W}(t) H_{U(0)}^{(r)}
H_{U(0)}^{(l)} = e^{itK_{U(0)}^{(r)} K_{U(0)}^{(l)}} S^{\ast} e^{-itH_{U(0)}%
^{(r)} H_{U(0)}^{(l)}} H_{U(0)}^{(r)} H_{U(0)}^{(l)} .
\end{equation}
Moreover, for any constant matrix $P \in\mathbb{C}^{n \times n} $ , we have
\begin{align}
\label{identity-2}\mathbf{W}(t)^{\ast}(P) = e^{itH_{U(0)}^{(r)} H_{U(0)}%
^{(l)}}(P) \in H^{2}_{n \times n}.
\end{align}
We also have
\begin{align}
\label{identity-3}\mathbf{W}(t)^{\ast}(U(t)) = U(0) \in\mathcal{H}_{+}%
^{s}(\mathbb{T}, \mathbb{C}^{n \times n}).
\end{align}

\end{lemma}

\begin{proof}
Set $\mathcal{Y}(t):=\mathbf{W}(t)^{\ast} S^{\ast} \mathbf{W}(t) H^{(r)}_{U(0)
}H^{(l)}_{U(0) } $. By (\ref{wequation}), we have
\begin{align}
\label{Cauchy ast}\frac{\mathrm{d}}{\mathrm{d}t} \mathbf{W}(t)^{\ast} = i
\mathbf{W}(t)^{\ast} T_{U(t) U(t)^{\ast}}^{(r)}.
\end{align}
As before, let $\mathcal{L}(t):= H^{(r)}_{U(t)}H^{(l)}_{U(t)} $. Then
\begin{align*}
\frac{d}{dt} \mathcal{Y}(t)=  &  \frac{d(\mathbf{W}(t)^{\ast} S^{\ast} )}{dt}
\mathbf{W}(t) \mathcal{L}(0)+\mathbf{W}(t)^{\ast} S^{\ast} \frac
{d\mathbf{W}(t)}{dt} \mathcal{L}(0)\\
=  &  i\mathbf{W}(t)^{\ast} T_{U(t) U(t)^{\ast}}^{(r)}S^{\ast}\mathbf{W}(t)
\mathcal{L}(0)+\mathbf{W}(t)^{\ast}S^{\ast}\left( -i T_{U(t) U(t)^{\ast}%
}^{(r)} \mathbf{W}(t)\right) \mathcal{L}(0)\\
= &  -i\mathbf{W}(t)^{\ast}\left[  S^{\ast} ,T_{U(t) U(t)^{\ast}}^{(r)}
\right]  \mathbf{W}(t)\mathcal{L}(0)\\
=  &  - i\mathbf{W}(t)^{\ast} \left[  S^{\ast} ,T_{U(t) U(t)^{\ast}}^{(r)}
\right]  \mathcal{L}(t)\mathbf{W}(t),
\end{align*}
where the last equality follows from  the unitary equivalences
(\ref{unitary equivalences}) in Lemma \ref{Sun-Corollary4.1}.

By Lemma \ref{commutator-STH}, we have
\begin{align*}
\left[  S^{\ast} ,T_{U(t) U(t)^{\ast}}^{(r)} \right]  \mathcal{L}(t) =  &
K^{(r)}_{U(t)}\left(  H^{(l)}_{U(t)} H^{(r)}_{U(t)}-K^{(l)}_{U(t)}%
K^{(r)}_{U(t)}\right)  H_{U(t)}^{(l)}\\
=  &  K^{(r)}_{U(t)}H^{(l)}_{U(t)} \mathcal{L}(t)- K^{(r)}_{U(t)}%
K^{(l)}_{U(t)}K^{(r)}_{U(t)}H_{U(t)}^{(l)}\\
=  &  S^{\ast} \mathcal{L}(t) \mathcal{L}(t)- K^{(r)}_{U(t)}K^{(l)}%
_{U(t)}K^{(r)}_{U(t)}H_{U(t)}^{(l)} .
\end{align*}
Again by Lemma \ref{Sun-Corollary4.1},
\begin{align*}
\frac{d}{dt} \mathcal{Y}(t)=  &  - i\mathbf{W}(t)^{\ast}\left( S^{\ast
}\mathcal{L}(t) \mathcal{L}(t)- K^{(r)}_{U(t)}K^{(l)}_{U(t)}K^{(r)}%
_{U(t)}H_{U(t)}^{(l)} \right)  \mathbf{W}(t)\\
=  &  - i\mathbf{W}(t)^{\ast}S^{\ast}\mathcal{L}(t) \mathcal{L}(t)\mathbf{W}%
(t)+ i\mathbf{W}(t)^{\ast}K^{(r)}_{U(t)}K^{(l)}_{U(t)}S^{\ast}H^{(r)}%
_{U(t)}H_{U(t)}^{(l)} \mathbf{W}(t)\\
\overset{\text{by} (\ref{unitary equivalences})}{=}  & -i\mathbf{W}(t)^{\ast
}S^{\ast}\mathcal{L}(t)\mathbf{W}(t)\mathcal{L}(0)+ i\mathbf{W}(t)^{\ast}
\mathbf{W}(t)K^{(r)}_{U(0)}K^{(l)}_{U(0)} \mathbf{W}(t)^{\ast}S^{\ast}
\mathbf{W}(t) \mathcal{L}(0)\\
=  &  - i \mathcal{Y}(t)\mathcal{L}(0)+ iK^{(r)}_{U(0)}K^{(l)}_{U(0)}
\mathcal{Y}(t).
\end{align*}

Now we have
\begin{align*}
\mathcal{Y}(t)=  &  e^{itK^{(r)}_{U(0)}K^{(l)}_{U(0)}} \mathcal{Y}%
(0)e^{-it\mathcal{L}(0) }\\
=  &  e^{itK^{(r)}_{U(0)}K^{(l)}_{U(0)}}\mathbf{W}(0)^{\ast} S^{\ast}
\mathbf{W}(0) \mathcal{L}(0)e^{-it\mathcal{L}(0)}\\
=  &  e^{itK^{(r)}_{U(0)}K^{(l)}_{U(0)}} S^{\ast} e^{-it\mathcal{L}%
(0)}\mathcal{L}(0).
\end{align*}
Thus (\ref{identity-1}) obtains.

Next, by (\ref{unitary equivalences}) and (\ref{Cauchy ast}), then for any
constant matrix $P \in\mathbb{C}^{n \times n}$, we have
\[
H_{U}^{(r)}H_{U}^{(l)}(P)=T_{UU^{\ast}}^{(r)}(P)
\]
and
\begin{align*}
\frac{\mathrm{d}}{\mathrm{d}t} (\mathbf{W}(t)^{\ast}(P)) = i \mathbf{W}%
(t)^{\ast} \mathcal{L}(t)(P) = i \mathcal{L}(0)\mathbf{W}(t)^{\ast} (P).
\end{align*}

The solution is
\begin{align*}
\mathbf{W}(t)^{\ast}(P) = e^{it\mathcal{L}(0)}\mathbf{W}(0)^{\ast}(P) =
e^{it\mathcal{L}(0)}(P).
\end{align*}
Thus (\ref{identity-2}) obtains.

Finally, by (\ref{cubic matrix}) and (\ref{Cauchy ast}), then
\begin{align*}
\frac{d}{dt}\left(  \mathbf{W}(t)^{\ast}(U(t)) \right)   &  = \frac{d}%
{dt}\left(  \mathbf{W}(t)^{\ast}\right)  (U(t)) +\mathbf{W}(t)^{\ast}%
(\partial_{t}U(t))\\
&  = i \mathbf{W}(t)^{\ast} T_{U(t) U(t)^{\ast}}^{(r)}(U(t)) -i\mathbf{W}%
(t)^{\ast}\left(  T_{U(t) U(t)^{\ast}}^{(r)}(U(t))\right) \\
&  =0.
\end{align*}
Thus,
\[
\mathbf{W}(t)^{\ast}(U(t))= \mathbf{W}(0)^{\ast}(U(0))=U(0),
\]
which proves (\ref{identity-3}).
\end{proof}

We now provide the proof of Theorem \ref{Sun explicit formula},
following the argument in \cite{Sun}. 


\textit{\textbf{Proof of Theorem \ref{Sun explicit formula}.}} We prove only
the first explicit formula. For every $\mathcal{V}\in H^{2}_{n\times n} $,
\begin{align*}
\underline{\mathcal{V}}(z)  &  =\sum_{m=0}^{\infty}\widehat{\mathcal{V}%
}(m)z^{m} =\sum_{k,j=1}^{n}\sum_{m=0}^{\infty}z^{m}\langle\mathcal{V}%
,S^{m}\mathbb{E}_{kj}\rangle\mathbb{E}_{kj}\\
&  =\sum_{k,j}\sum_{m=0}^{\infty}z^{m}\langle(S^{\ast})^{m}\mathcal{V}%
,\mathbb{E}_{kj}\rangle\mathbb{E}_{kj}\\
&  =\Delta\left(  \sum_{m=0}^{\infty}\left( zS^{\ast}\right) ^{m}
\mathcal{V}\right) \\
&  =\Delta\left(  (I-zS^{\ast})^{-1}\mathcal{V}\right)  .
\end{align*}
Then
\begin{align}
\label{harmonic eqaution}\underline{\mathcal{V}}(z)=\Delta\left(  (I-zS^{\ast
})^{-1}\mathcal{V}\right)  , z\in\mathbb{D}.
\end{align}

We first consider the case where the initial datum $U(0)$ is
rational. Let $U\in C^{\infty}(\mathbb{R},\mathcal{H}_{+}^{s}(\mathbb{T}%
,\mathbb{C}^{n\times n})) $ be a solution of (\ref{cubic matrix}), $s>1/2$.
Applying (\ref{harmonic eqaution}) to $\mathcal{V} =U(t)$, $\underline
{\mathcal{V}}(z)=\underline{U}(t,z)$ and using the unitarity of $\mathbf{W}%
(t)$, we get
\begin{align*}
\underline{U}(t,z)  &  =\Delta\left(  (I-zS^{\ast})^{-1}U(t)\right) \\
&  = \sum_{k,j} \langle\mathbf{W}(t)^{\ast}(I-zS^{\ast})^{-1}U(t),\mathbf{W}%
(t)^{\ast}\mathbb{E}_{kj}\rangle_{H^{2}_{n \times n}}\mathbb{E}_{kj},
\end{align*}
which yield
\begin{align}
\label{unitary equality}\underline{U}(t,z)= \sum_{k,j} \langle(I-z\mathbf{W}%
(t)^{\ast}S^{\ast}\mathbf{W}(t))^{-1}\mathbf{W}(t)^{\ast}U(t),\mathbf{W}%
(t)^{\ast}\mathbb{E}_{kj}\rangle_{H^{2}_{n \times n}}\mathbb{E}_{kj}.
\end{align}
Indeed,
\begin{align*}
\underline{U}(t,z)  &  = \sum_{k,j} \langle\mathbf{W}(t)^{\ast}\sum
_{m=0}^{\infty}(zS^{\ast})^{m}U(t),\mathbf{W}(t)^{\ast}\mathbb{E}_{kj}%
\rangle\mathbb{E}_{kj}\\
&  = \sum_{k,j} \langle\sum_{m=0}^{\infty}z^{m}\mathbf{W}(t)^{\ast}(S^{\ast
})^{m}\mathbf{W}(t)\mathbf{W}(t)^{\ast}U(t),\mathbf{W}(t)^{\ast}%
\mathbb{E}_{kj}\rangle\mathbb{E}_{kj}\\
&  = \sum_{k,j} \langle\sum_{m=0}^{\infty}z^{m}\left(  \mathbf{W}(t)^{\ast
}S^{\ast}\mathbf{W}(t)\right)  ^{m} \mathbf{W}(t)^{\ast}U(t),\mathbf{W}%
(t)^{\ast}\mathbb{E}_{kj}\rangle\mathbb{E}_{kj}\\
&  = \sum_{k,j} \langle(I-z\mathbf{W}(t)^{\ast}S^{\ast}\mathbf{W}%
(t))^{-1}\mathbf{W}(t)^{\ast}U(t),\mathbf{W}(t)^{\ast}\mathbb{E}_{kj}%
\rangle\mathbb{E}_{kj}.
\end{align*}

By Proposition 3.5 in \cite{Sun}, which relies on Kronecker's theorem
for rational symbols, we obtain
\[
\operatorname{Ran}H^{(r)}_{U(0)} = \operatorname{Ran}\left( H^{(r)}%
_{U(0)}H^{(l)}_{U(0)}\right) .
\]
Since $U(0) \in\operatorname{Ran}H^{(r)}_{U(0)}$, we have $U(0) \in
\operatorname{Ran}\left( H^{(r)}_{U(0)}H^{(l)}_{U(0)}\right) $. By
(\ref{identity-1}), the following identity holds on the range of
$H^{(r)}_{U(0)}H^{(l)}_{U(0)}$,
\begin{align}
\label{identity-4}\mathbf{W}(t)^{\ast} S^{\ast} \mathbf{W}(t) =  &
e^{itK_{U(0)}^{(r)} K_{U(0)}^{(l)}} S^{\ast} e^{-itH_{U(0)}^{(r)}
H_{U(0)}^{(l)}}.
\end{align}
It remains to plug identities (\ref{identity-2}), (\ref{identity-3}),
(\ref{identity-4}) into (\ref{unitary equality}). We finally obtain
\begin{align*}
\underline{U}(t,z) =  &  \sum_{k,j=1}^{n} \langle(I-ze^{itK_{U(0)}^{(r)}
K_{U(0)}^{(l)}} S^{\ast} e^{-itH_{U(0)}^{(r)} H_{U(0)}^{(l)}})^{-1}(U(0)),
e^{itH_{U(0)}^{(r)} H_{U(0)}^{(l)}} (\mathbb{E}_{kj})\rangle\mathbb{E}_{kj}\\
=  &  \Delta\left(  (I-ze^{-itH^{(r)}_{U(0)}H^{(l)}_{U(0)}}e^{itK^{(r)}%
_{U(0)}K^{(l)}_{U(0)}}S^{\ast})^{-1}e^{-itH^{(r)}_{U(0)}H^{(l)}_{U(0)}%
}(U(0))\right)  ,
\end{align*}
which proves the formula for rational initial data. The general case
$U(0)\in\mathcal{H}_{+}^{s}(\mathbb{T},\mathbb{C}^{n\times n})$, $s>1/2$,
follows by the approximation argument and the continuity of the flow map used
in the proof of Theorem 1.6 in \cite{Sun}.

\bigskip\noindent\textbf{Acknowledgement}: P. Ma was supported by NNSF of
China (Grant numbers 12171484,12571140), the Natural Science Foundation of
Hunan Province (Grant number 2023JJ20056), the Science and Technology
Innovation Program of Hunan Province (Grant number 2023RC3028), and Central
South University Innovation-Driven Research Programme (Grant number
2023CXQD032), Hunan Basic Science Research Center for Mathematical Analysis
(2024JC2002). M. Li was also supported by the Graduate Research and Innovation
Projects of Hunan Province (Grant number CX20250159).


\bigskip

\begin{thebibliography}{99}                                                                                               %

\bibitem {Conway}J. Conway, A course in functional analysis, 2nd ed., Graduate
Texts in Mathematics, Vol. 96, Springer, New York, 1990.

\bibitem {Gerard06}P. G\'{e}rard, Nonlinear Schr\"{o}dinger equations in
inhomogeneous media: wellposedness and illposedness of the Cauchy problem,
\textsl{Proc. Int. Congr. Math. (ICM)} , Madrid, Eur. Math. Soc. vol. 3
(2006), 156--182.

\bibitem {Gerard26}P. G\'{e}rard, Hardy spaces of holomorphic functions and
explicit formulae for a class of integrable partial differential equations,
\textsl{Proc. Int. Congr. Math. (ICM)}, Philadelphia, SIAM, vol. 2 (2026), 81--101.

\bibitem {Gerard08}P. G\'{e}rard, S. Grellier, L'\'{e}quation de Szeg\H{o}
cubique, S\'{e}minaire Goulaouic-Schwartz (2008), 1--19.

\bibitem {Gerard10}P. G\'{e}rard, S. Grellier, The cubic Szeg\H{o} equation.
\textsl{Ann. Scient. \'{E}c. Norm. Sup.} \textbf{43} (2010), 761--810.

\bibitem {Gerard12}P. G\'{e}rard, S. Grellier, Invariant tori for the cubic
Szeg\H{o} equation, \textsl{Invent. Math.} \textbf{187} (2012), no. 3, 707--754.

\bibitem {Gerard12-Anal}P. G\'{e}rard, S. Grellier, Effective integrable
dynamics for a certain nonlinear wave equation, \textsl{Anal. PDE,}
\textbf{5}(5) (2012), 1139--1155.

\bibitem {Gerard15}P. G\'{e}rard, S. Grellier, An explicit formula for the
cubic Szeg\H{o} equation. \textsl{Trans. Amer. Math. Soc.} \textbf{367}
(2015), 2979--2995.

\bibitem {Gerard17}P. G\'{e}rard, S. Grellier, The cubic Szeg\H{o} equation
and Hankel operators, \textsl{Ast\'{e}risque.} \textbf{389}, (2017): 114.

\bibitem {Gerard23}P. G\'{e}rard, A. Pushnitski, Unbounded Hankel operators
and the flow of the cubic Szeg\H{o} equation. \textsl{Invent. Math.}
\textbf{232} (2023), 995--1026.

\bibitem {Gerard-Pu24}P. G\'{e}rard, A. Pushnitski, The Cubic Szeg\H{o}
Equation on the Real Line: Explicit Formula and Well-Posedness on the Hardy
Class. \textsl{Comm. Math. Phys.} \textbf{405}, (2024), no. 7, Paper No. 167,
31 pp.

\bibitem {GerardJEMS}P. G\'{e}rard, A. Pushnitski, An inverse problem for
Hankel operators and turbulent solutions of the cubic Szeg\H{o} equation on
the line. \textsl{J. Eur. Math. Soc.} (JEMS) \textbf{27} (2025), no. 11, 4591--4648.

\bibitem {GuLuo}C. Gu, S. Luo, Invariant subspaces of the direct sum of
forward and backward shifts on vector-valued Hardy spaces. \textsl{J. Funct.
Anal.} \textbf{282} (2022), no. 9, Paper No.109419, 31 pp.

\bibitem {Halmos61}P. Halmos, Shifts on Hilbert spaces, \textsl{J. Reine
Angew. Math.} \textbf{208} (1961), 102--112.

\bibitem {Kato66}T. Kato, Perturbation theory for linear operators.
Springer-Verlag, New York, 1966.

\bibitem {Lax68}P. Lax, Integrals of nonlinear equations of evolution and
solitary waves, \textsl{Comm. Pure Appl. Math.} \textbf{21}(5) (1968), 467--490.

\bibitem {mats}Y. Matsuno, A direct approach for solving the cubic Szeg\H{o}
equation. \textsl{Stud. Appl. Math.} 153 (2024), no. 4, Paper No. e12770, 33 pp.

\bibitem {Nagy10}B. Sz.-Nagy, C. Foia\c{s}, H. Bercovici, L. K\'{e}rchy,
Harmonic Analysis of Operators on Hilbert Space, Springer, 2010.

\bibitem {Peller}V. Peller, Hankel operators and their applications. Springer,
New York, 2003.

\bibitem {Pocovnicu11}O. Pocovnicu, Traveling waves for the cubic Szeg\H{o}
equation on the real line, \textsl{Anal. PDE} \textbf{4} (2011),no. 3, 379--404.

\bibitem {Pocovnicu11-2}O. Pocovnicu, Explicit formula for the solutions of
the cubic Szeg\H{o} equation on the real line and applications,
\textsl{Discrete Contin. Dyn. Syst. A} \textbf{31} (2011) no. 3, 607--649.

\bibitem {Reed80}M. Reed, B. Simon, Methods of modern mathematical physics. I:
Functional analysis, Revised and enlarged edition, Academic Press (1980).

\bibitem {Reed802}M. Reed, B. Simon, Methods of modern mathematical physics.
II: Fourier analysis, Self-adjointness, Academic Press (1980).

\bibitem {Rodn04}I. Rodnianski, W. Schlag, Time decay for solutions of
Schr\"{o}dinger equations with rough and time-dependent potentials.
\textsl{Invent. Math.} \textbf{155} (2004), no. 3, 451--513.

\bibitem {Rose}M. Rosenblum, The absolute continuity of Toeplitz's matrices,
\textsl{Pacific J. Math.} \textbf{10} (1960), 987--996.

\bibitem {Sobolev}A. Sobolev, D. Yafaev, Multichannel scattering theory for
Toeplitz operators with piecewise continuous symbols. \textsl{Anal. PDE},
\textbf{15}(6), (2022), 1457--1486.

\bibitem {Sun}R. Sun, The matrix Szeg\H{o} equation, preprint, available on arXiv:2309.12136

\bibitem {Sun2}R. Sun, Matrix solutions of the cubic Szeg\H{o} equation on the
real line. \textsl{Math. Phys. Anal. Geom.} \textbf{28} (2025), no. 1, Paper
No. 4, 27 pp.

\bibitem {Timotin}D. Timotin, The invariant subspaces of $S\oplus S^{\ast}$,
\textsl{Concr. Oper.} \textbf{7}(1) (2020), 116--123.

\bibitem {Zakharov72}V. E. Zakharov, A. B. Shabat, Exact theory of
two-dimensional self-focusing and one-dimensional self-modulation of waves in
nonlinear media. \textsl{Sov. Phys. JETP} \textbf{34}(1) (1972), 62--69.
\end{thebibliography}
\end{document}